\documentclass{amsart}

\RequirePackage[OT1]{fontenc}
\usepackage{amsfonts, amsmath, amssymb, amsthm,url,dsfont, gensymb,multicol,bbm}
\usepackage[margin=3cm,footskip=1cm]{geometry}
\usepackage[english]{babel}
\RequirePackage[colorlinks,citecolor=blue,urlcolor=blue]{hyperref}
\usepackage{booktabs}

\RequirePackage[numbers]{natbib}
\usepackage[above,below,section]{placeins}

\usepackage[shortlabels]{enumitem}
\setlist{
	listparindent=\parindent,
	parsep=0pt,
}

\newtheorem{thm}{Theorem}

\newtheorem{defi}[thm]{Definition}

\theoremstyle{definition}
\theoremstyle{remark} 
\newtheorem{rem}[thm]{Remark}

\usepackage[dvipsnames]{xcolor}
\usepackage{tikz}
\usepackage{amsmath,bm,bbm,amsthm, amssymb}
\usepackage{fullpage}
\usepackage{color}
\usepackage{dsfont}
\usepackage{animate}
\usepackage{etoolbox}
\usepackage{comment}
\usepackage{pgf}
\usetikzlibrary{calc}
\usetikzlibrary{patterns}
\usetikzlibrary{arrows}
\usetikzlibrary{decorations.pathreplacing}
\usepackage[utf8]{inputenc}
\usepackage{pgfplots}

\usepackage{tcolorbox}
\usepackage{adjustbox}
\usepackage[final]{pdfpages}
\usepackage[ruled,vlined,linesnumbered,algosection,resetcount]{algorithm2e}
\usepackage{blkarray}
\usepackage{arydshln}
\usepackage{wrapfig}
\usepackage[font=scriptsize, labelfont=sc]{caption}

\usetikzlibrary{calc}
\theoremstyle{plain}
\newtheorem{theorem}{Theorem}
\newtheorem{proposition}[theorem]{Proposition}
\newtheorem{corollary}[theorem]{Corollary}
\newtheorem{lemma}[theorem]{Lemma}

\theoremstyle{definition}
\newtheorem{definition}[theorem]{Definition}

\newtheorem{example}[theorem]{Example}

\theoremstyle{remark}
\newtheorem{remark}[theorem]{Remark}

\numberwithin{theorem}{section}

\def\mc{\mathcal}
\def\ms{\mathsf}
\DeclareMathOperator{\Cov}{Cov}
\DeclareMathOperator{\Var}{\mathsf{Var}}

\newcommand{\X}{\mc X}

\DeclareMathOperator{\diam}{diam}
\newcommand{\E}{\mathbb{E}}
\def\P{\mathbb P}
\newcommand{\FF}{\mc F}
\newcommand{\PP}{\mc P}

\newcommand{\N}{\mathbb N}
\newcommand{\Z}{\mathbb{Z}}
\newcommand{\R}{\mathbb{R}}

\newcommand{\one}{\mathds{1}}
\newcommand{\Q}{\mathbb Q}
\newcommand{\eps}{\varepsilon}
\DeclareMathOperator{\de}{deg}

\DeclareMathOperator{\su}{sup}

\def\d{{\mathrm{d} }}

\newcommand{\GG}{\mc G}
\newcommand{\B}{\mc B}
\def\xx{\boldsymbol x}

\definecolor{wiasblue}   {cmyk}{1.0, 0.60, 0, 0}
\definecolor{mlugreen}{RGB}{172,6,52}
\definecolor{darmstadt}{RGB}{135,206,250}

\def\ni{\noindent}
\def\E{\mathbb E}
\def\M{\mathbb M}

\def\P{\mathbb P}
\def\Q{\mathsf{Leb}_{[0, 1]}}
\def\R{\mathbb R}
\def\Xx{\mathbb X}
\def\Uu{\mathbb U}
\def\Z{\mathbb Z}
\def\mc{\mathcal}
\def\ms{\mathsf}

\def\lrsa{\leftrightsquigarrow}
\def\la{\lambda}
\def\La{\Lambda}
\def\a{\alpha}

\def\s{\sigma}
\def\su{\subseteq}
\def\bs{\boldsymbol}

\def\b{\beta}
\def\de{\delta}

\def\es{\varnothing}
\def\one{\mathbbmss{1}}

\def\De{\Delta}

\def\G{\Gamma}

\def\ff{\infty}
\def\tff{\uparrow\infty}

\def\vp{\varphi}

\def\Leb{\ms{Leb}}
\def\d{{\rm d}}
\def\k{\kappa}

\def\GG{\mc G}
\def\GGG{\mc G^{\ms{expl}}}
\def\GGd{\mc G^{\ms{dc}}}
\def\GGc{\mc G^{\ms{cl}}}
\def\YY{\mc Y}

\def\NN{\mc N}
\def\BB{\mc B}
\def\PP{\mc P}
\def\FF{\mc F}
\def\CC{\mc C}

\def\XX{\mc X}

\def\f{\frac}

\def\r{\rho}

\def\im{\item}
\def\sm{\setminus}

\def\bep{\begin{proof}}
	\def\enp{\end{proof}}
\def\bepr{\begin{proposition}}
	\def\enpr{\end{proposition}}
\def\bec{\begin{corollary}}
	\def\enc{\end{corollary}}
\def\bea{\begin{align}}
	\newcommand\eea{\end{align}}
\def\beas{\begin{align*}}
	\def\eeas{\end{align*}}
\def\bet{\begin{theorem}}
	\def\ent{\end{theorem}}
\def\bee{\begin{example}}
	\def\ene{\end{example}}
\def\zz{\bs z}

\def\bede{\begin{definition}}
	\def\ende{\end{definition}}
\def\ber{\begin{remark}}
	\def\enr{\end{remark}}
\def\beca{\begin{cases}}
	\def\enca{\end{cases}}
\def\bel{\begin{lemma}}
	\def\enl{\end{lemma}}
\def\been{\begin{enumerate}}
	\def\enen{\end{enumerate}}

\def\S{\mathbb S}

\def\beit{\begin{itemize}}
	\def\enit{\end{itemize}}
\def\befr{\begin{frame}}
	\def\enfr{\end{frame}}
\def\ti{\times}

\def\Var{\ms{Var}}
\def\Cov{\ms{Cov}}

\def\diam{\ms{diam}}

\renewcommand\le{\leqslant}
\renewcommand\ge{\geqslant}

\def\x{x}

\def\bN{{\mathbf N}}

\def\xx{{\bs x}}
\def\ss{{\bs s}}
\def\rr{{\bs r}}
\def\ww{{\bs w}}

\def\y{y}
\def\yy{{\bs y}}

\def\Nlf{\mathbf N}

\def\dd{\,{\rm d}}

\def\var{\Var}

\def\becbb{\begin{center}\begin{tcolorbox}[{colback=Dandelion!20}]}
		\def\encbb{\end{tcolorbox}\end{center}}
\def\becb{\begin{center}\begin{tcbox}[{colback=Dandelion!20}]}
		\def\encb{\end{tcbox}\end{center}}

\def\bef{\begin{figure}[!h]}
	\def\enf{\end{figure}}
\def\bels{\begin{lstlisting}}
	\def\enls{\end{lstlisting}}

\def\betp{\begin{tikzpicture}}
	\def\entp{\end{tikzpicture}}
\def\co{\colon}
\def\endo{\end{document}}

\def\ZZ{\mc Z}

\def\N{\mathbb N}

\def\dtv{d_{\ms{TV}}}

\def\Trr{T^{\ms{RCM}}}
\def\Tf{T^{\ms{RCM, \ff}}}
\def\Trd{T^{\ms{RCM}, \ms{dc}}}
\def\Trc{T^{\ms{RCM}, \ms{cl}}}

\def\Pd{\PP^{}}
\def\Pds{\PP^{*}}

\def\X{\XX}

\def\dk{d_{\ms K}}

\def\tdc{T^{\ms{dc}}}
\def\Ep{E_{\ms{perc}}}
\def\Ed{E_{\ms{dens}}}

\begin{document}
	
	\title{
	Normal and Poisson approximation for Gibbs point processes with pair potentials 
	}
	
	\begin{abstract}
		We provide a Poisson approximation result for dependent thinnings of Gibbs point processes as well as qualitative and quantitative central limit theorems for geometric functionals of Gibbs point processes in increasing observation windows. The present paper extends prior work on finite-range Gibbs processes to processes with repulsive pairwise interaction of unbounded interaction range as well as processes on marked Euclidean space. The proofs rely on coupling different Gibbs processes using the disagreement coupling technique, which we generalize to infinite-volume domains under a suitable non-percolation condition. For the case of repulsive pairwise interactions, we introduce a version of disagreement coupling that constructs the Gibbs process by thinning a random connection model thus making previous approximation methods more flexible.
	\end{abstract}
\author{Christian Hirsch}
\address[Christian Hirsch]{Department of Mathematics, Aarhus University, Ny Munkegade 118, 8000 Aarhus C, Denmark}
\email{hirsch@math.au.dk}
\address[Christian Hirsch]{DIGIT Center, Aarhus University, Finlandsgade 22, 8200 Aarhus N, Denmark}

\author{Moritz Otto}
\address[Moritz Otto]{Mathematical Institute, Leiden University, Netherlands}
\email{m.f.p.otto@math.leidenuniv.nl}

\author{Anne Marie Svane}
\address[Anne Marie Svane]{Department of Mathematics\\
Aalborg University \\
Aalborg, 8000, Denmark}
\email{annemarie@math.aau.dk}	

\maketitle
	
	%
	%
	\section{Introduction}
\label{sec:int}
Point processes are important  in statistics as models  for random point patterns such as locations of plant or animal species,  cells in biological tissue, or disease occurrences \cite{Baddeley:Rubak:Wolf:15}.   The simplest point process model is the Poisson point process, which models complete spatial randomness, see \cite{LP} for a comprehensive overview. However, in many instances such a model is too simplistic due to spatial dependencies between points. This makes it essential to consider models for dependent point patterns. In the context of spatial statistics, a particularly popular class of models is Gibbs point processes. These are motivated by statistical physics as models for interacting particle systems \cite{ruelle}. The present paper considers  functionals of Gibbs point processes observed in large windows. We provide two types of approximate distribution results.
\been
	\item \textbf{Dependent thinnings of Gibbs point processes.} In such a thinning, each point of the Gibbs process is either kept or removed depending on the configuration of the surrounding points. We provide bounds on the Kantorovich-Rubinstein distance between the thinned process and a Poisson process. Results like this have recently been established for both Poisson processes \cite{BSY202, ober} and also for Gibbs processes \cite{dp}. In the absence of complete spatial independence, one needs to impose mixing conditions on the Gibbs process. In the case of pair-wise interaction, this can be obtained by assuming an exponential decay condition on the diameter of connected clusters in a dominating random connection model. This relaxes the non-percolation assumption made in \cite{dp}. As a second step of generalization, we consider thinning functions that depend on the point configuration within a random radius of exponentially decaying size instead of a deterministic radius as in \cite{dp}.
\item \textbf{Geometric functionals of point patterns.} We consider functionals of a point pattern such as the average nearest neighbor distance, Voronoi cell characteristics, and topological summaries called persistent Betti numbers \cite{shirai}. We provide normal approximation results for such functionals when applied to a Gibbs point process observed in a large window. Results of this type have a long history for Poisson processes, see e.g.\ \cite{Mehler,penrose}, and more recently, results have also appeared for Gibbs point processes \cite{benes,CX22,hos,gibbs_limit,gibbsCLT}. We show (i) a qualitative central limit theorem (CLT) requiring minimal assumptions on the stabilization of the functional in increasing windows and on the cluster sizes of a dominating percolation model and (ii) a quantitative CLT measured in the Kolmogorov distance under an additional exponential stabilization condition and assuming exponential tails of cluster diameters in the dominating percolation model. 
\enen

For both types of problems, the main contribution of this paper is to generalize existing results to a larger class of Gibbs point processes. Since the asymptotic theory of Gibbs point process is a vigorous research stream, we briefly review the existing results in literature. A key challenge for deriving asymptotic results is handling the spatial correlations of the point process. Mainly two techniques have been used in the literature, namely a {perfect simulation technique} and {disagreement coupling}. Both seek to control spatial correlations by the clusters of an associated percolation model.
\been
\im {\bf Perfect simulation.} The idea is to represent the Gibbs process as a thinning of a free birth-death process \cite{fernandez,ferrari}. The application to CLTs for Gibbs processes was pioneered in \cite{gibbs_limit,gibbsCLT}. Later, quantitative Wasserstein bounds were derived for point processes with sufficiently fast correlation decay \cite{CX22}, which could then be verified for the Gibbs processes from \cite{gibbs_limit,gibbsCLT}. 
While the perfect simulation technique is easy to work with, the arguments are only valid for short-range interactions  satisfying that a certain branching process derived from a space-time percolation model is sub-critical. 
\im {\bf Disagreement coupling.} This technique was originally introduced for discrete models  \cite{maes} and only recently transferred to the continuum case \cite{HT,HHT,dp}. It provides explicit couplings between different Gibbs models by thinning a common dominating Poisson process allowing to control the long range effects of locally modifying the model.  Disagreement coupling was used in \cite{dp} for Poisson approximation and some early CLTs for U-statistics were obtained in \cite{benes}.  More general normal approximation results were obtained in \cite{hos}. There, the Gibbs point processes were restricted to live in Euclidean space and have a finite interaction range satisfying a non-percolation condition on a dominating Boolean model. This relaxes the branching process condition for the perfect simulation technique, but, as we will see in this paper, it is still  unnecessarily strict in many cases. 
\enen
This paper continues the disagreement coupling studies from \cite{hos} and  considers two generalizations.
\been
	\item {\bf Models with local interactions.} We consider Gibbs processes in marked Euclidean space satisfying a local interaction property and further assume non-percolation of an associated random graph.  Building on the general disagreement coupling framework of \cite{dp} in finite-volume domains, we generalize the constructions of \cite{hos} to show how Gibbs processes in infinite-volume domains can be explicitly constructed by Poisson thinnings and we give quantitative bounds on the effect of locally modifying the model. As a new result compared to \cite{hos}, we give explicit couplings between infinite-volume processes with differing boundary conditions.  
	\item {\bf Models with pair potential.} To weaken the percolation restriction and to allow models with infinite interaction range, we take up an idea from \cite{betsch} to replace the associated percolation model by a random connection model (RCM), where the connection probabilities depend on the strength of the interaction between pairs of points. In \cite{betsch}, this is done by approximating the RCM by models of the form in \cite{dp}. However, for our proofs, we need to set up a version of disagreement coupling acting directly by thinning the clusters of an RCM in both finite- and infinite-volume domains.  
\enen
To summarize, the main contributions of the present paper is that it is the first to give both quantitative and qualitative Poisson and normal approximation results for a general class of pair potentials satisfying a suitable non-percolation condition, see \eqref{cond:perc_RCM} below.

The paper is structured as follows. In Section \ref{sec:mod}, we introduce the Gibbs models we consider, and the main theorems are stated in Section \ref{sec:results} with examples of applications given in Section \ref{sec:exa}. In Section \ref{sec:const}, we generalize the disagreement coupling results from \cite{hos} to standard Borel spaces, and in Section \ref{sec:RCM}, we set up a disagreement coupling scheme for the RCM case. In both sections,  proofs are deferred to Section \ref{sec:proofs_const} and~\ref{sec:proofsRCM}, respectively, allowing the reader to skip them upon first reading. We first give the proof of Theorem \ref{thm:3} in Section \ref{sec:weak} since this is conceptually simpler. Section \ref{sec:pac} is devoted to some common Palm coupling framework for the proofs of Theorems \ref{thm:1} and \ref{thm:2}, which are then given in Sections \ref{sec:papx} and \ref{sec:malliavin_stein}, respectively. Finally, Section~\ref{sec:proofthm2} contains a proof of a technical lemma used in the proof of Section \ref{sec:malliavin_stein}.

%
%
\section{Model and main results}
\label{sec:mod}

%
%
\subsection{Gibbs point processes}
\label{ss:Gibbs}

Let $\Xx$ be a standard Borel space with Borel sigma-algebra $\B$ equipped with a diffuse measure $\la$. We assume that there is a sequence $W_1\su W_2 \su \dots \su \Xx$ such that 
\begin{align}
	\label{eq:la}
	\Xx =\bigcup_n W_n \quad \text{ and } \quad \la(W_n)<\ff, \qquad n\ge 1.
\end{align}
Let $\B_0$ denote the set of \emph{bounded Borel sets}, i.e.\ sets $B\in \B$ with $B\su W_n$ for some $n\ge 1$. Let $\Nlf$ denote the space of \emph{locally finite} subsets of $\Xx$, i.e., sets $\vp \su \Xx$ such that $\vp \cap W_n$ is finite for all $n\ge 1$.  We sometimes identify elements $\vp \in \Nlf$ with their associated counting measures and  write $\vp(B) = |\vp \cap B|$ for the cardinality of $\vp\cap B$.
We equip $\Nlf$ with the smallest sigma-algebra $\NN$ making the maps $\vp \mapsto \vp(B)$ measurable for
all $B\in\BB_0$. Let $\Nlf_0\su \Nlf$ denote the subspace of finite sets. For $B\in \BB$,  we let $\Nlf_B\su \Nlf$ denote the subspace of locally finite sets contained in $B$ and for $\vp \in \Nlf$ we write $\vp_B:=\vp\cap B\in \Nlf_B$.

A {\em simple point process} on $\Xx$ is a measurable map
 $(\Omega,\FF)\to(\Nlf,\NN)$ where $(\Omega,\FF,\P)$ is some probability space. A Gibbs point process is defined in terms of a {\em Papangelou intensity} (PI) which is a function $\k : \Xx \ti \Nlf \to [0,\infty)$ satisfying the so-called \emph{cocycle condition}
\begin{equation}\label{cond:cocycle}
	\k(\xx,\vp)\k(\yy,\vp \cup \{\xx\} ) = \k(\yy,\vp)\k(\xx,\vp \cup \{\yy\}). 
\end{equation}

\bede A point process $\X$ is a {\em Gibbs point process on $\Xx$} with PI $\k$ if the {\em GNZ-equation} \cite{Georgii76,NgZe79} 
\begin{equation}\label{cond:GNZ}
	\E\Big[\int_{\Xx} f(\xx,\X) \X(\d \xx)\Big] =  \int_{\Xx}\E\big[ f(\xx,\X\cup\{\xx\}) \k(\xx,\X)\big] \la(\d \xx)
\end{equation}
holds for any measurable $f: \Xx\ti \Nlf \to [0,\infty)$.
If $Q\su \Xx$ and $\psi\in \Nlf_{\Xx\sm Q}$, then \emph{the Gibbs process $\XX(Q,\psi)$ on $Q$ with  boundary conditions $\psi$} and PI $\k$ is a point process on $Q$ satisfying
\begin{equation}\label{cond:GNZ_boundary}
	\E\Big[\int_Q f(\xx,\X(Q,\psi)) \X(Q,\psi)(\d \xx) \Big]= \int_Q  \E\big[f(\xx,\X(Q,\psi)\cup \{\xx\}) \k(\xx,\X(Q,\psi) \cup \psi)\big] \la(\d \xx).
\end{equation}
\ende
For a more detailed introduction to Gibbs point processes, we refer the reader to  \cite{betsch2,dereudre}.
Note that $\XX(Q,\psi)$ is the same as a Gibbs process on $\Xx$ with PI $\k'(\xx,\vp)=\k(\xx, \vp \cup \psi)\one\{\xx \in Q\}$.
We will assume throughout that $\k$ is \emph{locally stable}, i.e.\
\begin{equation}\label{cond:PIbound}
	\sup_{\vp \in \Nlf}\k(\xx,\vp) \le 1\quad \text{ for all }\xx \in \Xx.
\end{equation}

\ber\label{rem:activity}
In the literature, local stability is usually formulated as 
	$\sup_{\vp \in \Nlf}\k(\xx,\vp) \le \a(\xx)$ for all $\xx \in \Xx$, for some $\a:  \Xx \to [0 ,\infty)$ with 
$\int_{W_n} \a(\xx)\la(\d \xx) <\infty $ for all $ n\ge 1$. However, changing the reference measure to $\a \la$ and the PI to $\k/\a$ (with $0/0:=0$), we may assume $\a=1$. 
\enr

Local stability ensures that the Gibbs point process is stochastically dominated by a Poisson process of intensity $\la$.
Throughout, we let $\Pd$ denote a Poisson point process on $\Xx$ with intensity measure $\la$ and we let $\Pd_Q$ denote its restriction to $Q\in \BB$. We refer to \cite{LP} for a comprehensive introduction to Poisson point processes. Suppose $Q\in \BB_0$ and $\psi \in \Nlf_{\Xx\sm Q}$. Define $\tilde{\k}(\cdot ,\psi):\Nlf_Q \to [0,\infty)$ by  $\tilde{\k}(\es,\psi)=1$ and
$$\tilde{\k}(\{\xx_1,\dots,\xx_m\},\psi) = \prod_{i=1}^m \k(\xx_i, \{\xx_1,\dots,\xx_{i-1}\}\cup \psi )$$
where $\xx_1,\dots,\xx_m\in Q$ are pairwise different points. The cocycle condition \eqref{cond:cocycle} ensures that $\tilde \k$ does not depend on the ordering of $\xx_1,\dots,\xx_m$. 
The \emph{partition function} is given by
$$Z_Q(\psi) = \E[\tilde{\k}(\Pd_Q,\psi)]$$
 As noted in \cite{dp}, the assumptions \eqref{cond:PIbound} and $Q\in \BB_0$  ensure $0<Z_Q(\psi ) <\infty$. In this case, $\X(Q, \psi)$ has a density with respect to $\Pd_Q$ given by $\frac{1}{Z_Q(\psi)}\tilde{\k}(\cdot,\psi)$. In particular, $\X(Q,\psi) $ exists and its distribution is unique. 

Without the assumption $Q\in \BB_0$, the Gibbs point process is not guaranteed to exist. Moreover,  the Gibbs distribution is not necessarily unique as demonstrated for the Widom-Rowlinson model in \cite{houdebert}. A phase transition between uniqueness and non-uniqueness is believed to occur for most models. Thus, we  make some further assumptions to ensure uniqueness. The first is that the PI has \emph{local interactions}, that is, 
we assume that there is a symmetric relation $\sim$ on $\Xx$, which is measurable in the sense that $(\xx,\yy)\mapsto \one\{\xx\sim \yy\}$ is measurable. Let $N(\xx)=\{\yy\in \Xx\mid \xx\sim \yy\}$ denote all neighbors of $\xx$ with respect to $\sim$.  Then, we say that the model has \emph{local interactions} if
\begin{equation}\label{cond:PI_sim}
	\k(\xx,\vp) = \k(\xx,\vp\cap N(\xx)).
\end{equation}
We assume $\la(N(\xx)) <\infty $ for all $\xx \in \Xx$. Sometimes, we require this to hold uniformly in the sense that
%
%
\begin{equation}\label{cond:unif_sim}
\int_{W_n } \la(N(\xx) ) \la(\dd \xx) <\ff, \qquad n \ge 1.
\end{equation}
More generally, we write $N(B)=\bigcup_{\xx \in B} N(\xx)$ for $B\in \BB$.

For any $\vp \in \Nlf$, the relation $\sim$  defines a graph $G(\vp) $ with vertex set $\vp$ and an edge between $\xx,\yy\in \vp$ whenever $\xx\sim \yy$. For $\xx \in\Xx$, we let $C(\xx,\vp)$ be the connected component of $\xx$ in $G_{\xx}(\vp):=G(\vp\cup\{\xx\})$. We refer to the vertices of a connected component as a \emph{cluster}.
 We write $A\lrsa_\vp B$ if there is a point in $\vp\cap A$ and a point in $\vp\cap B$ that are in the same connected component of $G(\vp)$.
Intuitively, $C(\xx,\vp)$ captures how far the interactions with $\xx$ propagate. This motivates the following non-percolation assumption: For $\la$-almost all $\xx \in \Xx$,
\begin{equation}\label{cond:perc}
	\P(|C(\xx,\Pd)|<\ff ) =1.
\end{equation}
It will often be convenient to reformulate \eqref{cond:perc} as
\begin{equation}\label{cond:perc2}
	\lim_{n\to \infty} \P\big( C(\xx,\Pd) \not\su W_n \big) =0.
\end{equation}
This ensures that the interactions with $\xx$ do not propagate infinitely far out. Note that $C(\xx,\Pd)$ is the connected component of $\xx$ in the graph $G(\Pd_{\xx})$ where $\Pd_{\xx}=\Pd\cup \{\xx\}$ denotes the Palm distribution of $\Pd$ at $\xx$.
With this assumption, we have the following result from \cite[Thm. 9.1]{betsch2}.

%
%
\bet[Existence and uniqueness of Gibbs processes]\label{thm:unique1}
A Gibbs distribution $\X(\Xx,\es)$ with PI satisfying \eqref{cond:PIbound},  \eqref{cond:PI_sim} and \eqref{cond:perc} is guaranteed to exist and be unique.
\ent

While the above setup is very general, our results will be stated for a more specific setting. We will consider two main scenarios: 
\been
\im[] {\bf Model (L)}. Gibbs models with \emph{local interactions} on marked Euclidean spaces, and
\im[] {\bf Model  (P)}. Gibbs processes with \emph{pair potentials} on marked Euclidean spaces.
\enen 
Below, we define the two models more precisely.

%
%
\subsection{Model (L): Gibbs models with local interactions on marked Euclidean spaces}
\label{ss:marked}

In the following, we say that $\Xx$ is a \emph{marked Euclidean space} if  $\Xx=\R^d \ti \M$, where the mark space $\M$ is some  Polish space equipped with a probability measure $\mu$ and $\R^d$ is equipped with the Lebesgue measure $\Leb$. In this setting, we choose the sequence  $W_1\su W_2 \su \cdots \su \R^d \ti \M$ to be marked $d$-cubes $W_n=[-n/2,n/2 ]^d\ti \M$  of side length $n$. We write points $\xx\in \Xx$ as $\xx=(x,m)$.

\begin{defi}
 Model (L) is a Gibbs point process on marked Euclidean space with reference measure $\la=\a \Leb \otimes \mu$ where  $\a > 0$ is the \emph{activity}. The PI is assumed to satisfy \eqref{cond:PIbound},  \eqref{cond:PI_sim}, \eqref{cond:unif_sim} and \eqref{cond:perc2}. 
\end{defi}
 
For $A\su \R^d\ti \M$ and $y\in \R^d$, we define the set translation $A+y=\{(x+y,m)\in \R^d \times \M \mid (x,m)\in A\}$. We say that a PI $\k$ on marked Euclidean space is \emph{translation-invariant} if 
$$\k\big((x,m),{\vp}\big) = \k\big((x+y,m), {\vp}+y\big)$$ 
for all $(x,m)\in \R^d\ti \M$, ${\vp}\in\Nlf_{\R^d\ti \M } $, and $y\in \R^d$.  
 We sometimes need a stronger version of the non-percolation condition \eqref{cond:perc2} requiring exponential tails of the spatial diameter of clusters. More precisely, we assume that there are $c_{\ms D, 1},c_{\ms D, 2}>0 $ such that
\begin{equation}\label{cond:sharp}
        \sup_{x \in\R^d}\int_{\M} \P\big(\diam(C((x,m),\Pd)) > r\big) \mu(\d m) \le c_{\ms D, 1} \exp(-c_{\ms D, 2} r).
\end{equation}
Here, $\diam (A)$ refers to the spatial diameter of  the projection to $\R^d$ of a subset $A\su \R^d \ti \M$.

Sometimes, we will also need to assume exponential tails of the neighborhood of a point in the sense that there are constants $c_{\ms{v},1},c_{\ms{v},2}>0$ such that
 \begin{equation}\label{cond:exponential_sim}
\sup_{x  \in \R^d} \int_{\M} \int_{\Xx\sm B_r(\xx)} \one \{(x,m) \sim (y,m')\} \la(\d(y,m'))\mu(\d m )\le c_{\ms{v},1}\exp({-c_{\ms{v},2} r}), \qquad r>0,
\end{equation}
where $B_r(x)$ is the Euclidean ball of radius $r$ around $x\in \R^d$ and for $\xx=(x,m)\in \R^d \ti \M$ we write $B_r(\xx):=B_r(x)\ti \M$. Note that \eqref{cond:exponential_sim} implies \eqref{cond:unif_sim}.

We now discuss specific examples of point processes that satisfy the conditions.
\bee[Euclidean point processes with finite interaction range]\label{ex:hos}
A Gibbs point process on $\Xx = \R^d$ has \emph{finite interaction range} $r_0>0$ if
$$\k(\xx,\vp) = \k\big(\xx,\vp \cap B_{r_0}(\xx)\big).$$
This corresponds to a symmetric relation given by $ \xx_1 \sim \xx_2 $ if and only if $\|\xx_1 - \xx_2\|\le r_0$. Thus,  the condition \eqref{cond:perc} is equivalent to non-percolation of the Boolean model with intensity $\a$ and balls of radius $r_0/2$. Moreover, the criterion \eqref{cond:perc2} implies \eqref{cond:sharp}, see e.g.\ \cite{raoufi_sub,ziesche}. This was the setting studied in \cite{hos}.
\ene

%
%
\bee[Interacting particle models]
\label{ex:particles}
Another relevant model class is interacting particle processes studied in e.g.\ \cite{benes,betsch2}. Here, $\M$ is the family of compact sets having center of gravity at the origin. A pair $(x,C)\in \R^d \ti \M$ is identified with the particle $x+C\su \R^d$. In a model of non-overlapping particles, the PI is given by 
$$\k\big((x,C),\vp\big) = \one\bigg\{(x+C)\cap\bigcup_{(y,C')\in \vp} (y+C') = \es \bigg\}. $$
Here, we can take the symmetric relation to be $(x,C)\sim (y,C')$ if and only if $(x+C)\cap (y+C')\neq \es$. The percolation condition \eqref{cond:perc} corresponds to non-percolation of a germ-grain model with grain distribution $\mu$.
It immediately implies \eqref{cond:sharp} when grain diameters are bounded \cite{ziesche}. In the case where the grains are balls of random radii and the radius distribution has exponential tails, \eqref{cond:perc} implies \eqref{cond:sharp} by \cite{raoufi_sub}. However, if the tail probabilities of the radii decay slower than exponentially, \eqref{cond:sharp} cannot hold, but the non-percolation condition \eqref{cond:perc} may still hold for low intensities $\a$ \cite{gouere}. 
Other interacting particle models considered in the  literature include the  Widom-Rowlinson  model \cite{widom}, segment processes \cite{betsch2}, and the general models considered in \cite{benes,betsch2}.
\ene


\subsection{Model (P): Gibbs models with pairwise interaction on marked Euclidean space}\label{ss:pair}

Some of the most well-studied Gibbs models are models with pairwise interaction, i.e., the PI has the form
\begin{equation*}
	\k(\xx,\vp) = \exp\bigg(-\sum_{\yy\in \vp\sm \{\xx\}} v(\xx,\yy) \bigg),
\end{equation*}
where the \emph{pair potential} $v: \Xx \ti \Xx \to [0,\infty] $ is symmetric in its two arguments. The assumption $v\ge 0$  is necessary to ensure \eqref{cond:PIbound} and implies that the process is repulsive at all distances. We remark that there are physically relevant models like the Lennard-Jones potential that do not satisfy this.  In Model (P), we always  require that
\begin{equation}\label{cond:int_v}
	\int_{\Xx} \pi(\xx, \yy) \la( \d \yy )<\ff, \qquad \xx \in \Xx,
\end{equation}
 where
 \begin{equation}\label{eq:def_pi}
 	\pi(\xx,\yy) = 1-e^{-v(\xx,\yy)}.
 \end{equation}
  This ensures existence of the Gibbs process \cite{gibbsexist}. 
For the proofs, we further need that \eqref{cond:int_v} is locally integrable in the sense that
\begin{equation}\label{cond:local_int_v}
	\int_{W_n}\int_{\Xx}\pi(\xx, \yy)\la( \d \yy) \la(\d \xx) < \ff, \qquad n \ge 1.
\end{equation}

If the pair potential $v(\xx,\cdot)$ has bounded support for every fixed $\xx$, one could define a symmetric relation as $\xx\sim \yy$ if $v(\xx,\yy) > 0$ and be in the setting of Model (L). However, the associated non-percolation  criterion \eqref{cond:perc2} would be unnecessarily strict. 
Following \cite{betsch}, one may instead relate the Gibbs process to a random connection model (RCM) where, rather than adding deterministic edges between $\xx,\yy\in \Pd$ whenever $v(\xx,\yy)>0$, edges are only added with a certain probability depending on $v$. This leads to a larger percolation threshold and, moreover, does not require any restrictions on the support of $v$.
To make this idea precise, we first define the RCM.

%
%
\bede[RCM $\G(Q, \psi)$]
\label{def:RCM}
For $\vp,\psi\in \Nlf$, define the random graph $\G(\vp,\psi)$ as the simple graph 
 with vertex set $\vp \cup \psi$ and edges added independently with probability  $\pi(\xx, \yy) $
of having an edge between $\xx \in \vp$ and $\yy\in \vp \cup\psi$.  The \emph{random connection model} $\G(\PP_Q, \psi)$ is the  graph with vertex set consisting of the Poisson process $\PP_Q $ together with a set of fixed points $\psi$. Conditionally on $\PP_Q=\vp$ it has distribution $\G(\vp,\psi)$. The distribution of $\G(\PP_Q, \psi)$ is denoted $\G(Q,\psi)$.
\ende

For $Q\in \BB_0$ and $\psi\in \Nlf_{\Xx\sm Q}$, the density of $\X(Q,\psi)$ with respect to $\Pd_Q$ is given, up to a normalizing constant, by 
\begin{equation*}
	\tilde{\k}(\{\xx_1,\dots,\xx_l\},\psi) = \prod_{i=1}^l \exp\bigg(-\sum_{\yy\in \{\xx_1,\dots,\xx_{i-1}\}\cup \psi} v(\xx_i,\yy) \bigg) = \prod_{i=1}^l \prod_{\yy\in  \{\xx_1,\dots,\xx_{i-1}\}\cup\psi} e^{-v(\xx_i,\yy)}.
\end{equation*}
In particular, $\tilde \k \big(\{\xx_1,\dots,\xx_l\},\psi\big)$ equals the probability that $\G\big(\{\xx_1,\dots,\xx_l\},\psi\big)$ has no edges. Thus, one may think of $\X(Q,\psi)$  as having a density with respect to the RCM $\G(\PP_Q,\psi)$. The density is proportional to the indicator $\one\{\G(Q,\psi)\text{ has no edges}\}$  with normalizing constant ${Z}_Q(\psi) = \P(\G(Q,\psi)\text{ has no edges}).$
This suggests a close relation between the Gibbs process and the RCM, and in fact, in Section \ref{sec:RCM} we will construct $\X(Q,\psi)$ as a subset of the vertices of $\G(Q,\psi)$ that have no edges between them.

For $\xx \in\Xx$, let $\G_{\xx}(\Pd, \es) := \G(\Pd_{\xx}, \es)$ be the RCM of the Palm process at $\xx$ and let $C\big(\xx,\G_{\xx}(\Pd,\es)\big)$ denote the connected component of $\xx$ in  $\G_{\xx}(\Pd,\es)$. The equivalent of the non-percolation criterion \eqref{cond:perc2} for the RCM is that for $\la$-almost all $\xx \in \Xx$,
\begin{equation}\label{cond:perc_RCM}
	\lim_{n\to \ff} \P_{\xx}\Big( C(\xx,\G_{\xx}(\Pd,\es)) \cap (\Xx \sm W_n) \ne \es \Big) =0.
\end{equation}

\bet[Existence and uniqueness of Gibbs processes with pair potential  \cite{betsch,gibbsexist}]\label{thm:unique2}
A Gibbs distribution on $\Xx$ with non-negative pair-potential satisfying \eqref{cond:int_v} and \eqref{cond:perc_RCM} exists and is unique.
\ent

Our main results will again be concerned with models on marked Euclidean space, hence we define Model (P) as follows.

\begin{defi}
	Model (P) is a Gibbs point process on marked Euclidean space with reference measure $\la = \a \Leb \otimes \mu$, $\a>0$, having a non-negative pair potential satisfying \eqref{cond:int_v}, \eqref{cond:local_int_v} and \eqref{cond:perc_RCM}.
\end{defi}

In the main Theorems \ref{thm:1} and~\ref{thm:2} below, we further require exponential decay of cluster sizes of the RCM analogous to \eqref{cond:sharp}. That is, there are $c_{\ms D, 1},c_{\ms D, 2}>0$ such that
\begin{equation}
	\label{cond:sharp_vv}
	\sup_{x \in\R^d}\int_{\M}\P\Big(\diam\big(C((x,m),\G_{(x,m)}(\Pd,\es))\big) > r\Big) \mu(\d m) \le c_{\ms D, 1} \exp(-c_{\ms D, 2} r).
\end{equation}
Conditions for this to hold can be found in \cite{mikhail}. 
We may also assume exponential decay of $\pi$, i.e., there exist $c_{\ms v, 1}, c_{\ms v, 2}>0$ such that 
\begin{align}
	\label{cond:sharp_v}
	\int_{\M} \int_{\M} \pi((x,m), (y,m')) \mu(\d m) \mu(\d m')\le c_{\ms v, 1} \exp\big(-c_{\ms v, 2} \|x-y\|\big) , \qquad x,y \in \R^d. 
\end{align}

 We will sometimes write $\XX^\la(Q,\psi)$ and $\G^\la(Q,\psi)$ to emphasize that the reference measure is $\la$ and $\G(Q,\psi,\pi)$ if we wish to emphasize the connection probabilities $\pi$. Note that $\XX^\la(Q,\psi)$ has the same distribution as $\XX^{\rho_\psi\la}(Q,\es)$ with  reference measure $\rho_\psi\la$ where 
	\begin{align}
		\label{eq:rp}
		\rho_\psi(\xx) := \prod_{\yy\in \psi} \big(1-\pi(\xx,\yy)\big),
	\end{align}
	since the same GNZ equation is satisfied.

\bee(Strauss model)
A popular model in spatial statistics is the \emph{Strauss process} \cite{strauss}, which is a point process on $\R^d$ with  pair potential $v((x,m),(y,m'))= -\ln(\gamma ) \one\{\|x-y\|\le R\}$, where $\gamma \in [0,1)$,  and reference measure $\beta\Leb$. The special case $\gamma=0$ is the hard-sphere model. 
\ene

The non-overlapping particle model from Example \ref{ex:particles} is another example with pairwise interaction. The percolation models arising from considering it as a model of type (L) or (P) coincide. An example of a model that does not have pairwise interactions is the Widom-Rowlinson model. See \cite{rasmus} for more examples.

\begin{rem}
	A different condition for uniqueness of a Gibbs model with non-negative pair potential satisfying \eqref{cond:int_v} was given in \cite{michelen2}. This was based on suitably counting self-avoiding walks in the RCM. The relation to the non-percolation criterion \eqref{cond:sharp_vv} is not fully understood, but the authors of \cite{michelen2} argue that at least in high dimensions, their criterion is better. Under the same criterion, they show a strong spatial mixing property in \cite{michelen}. It is possible that this mixing condition suffice to derive some CLTs, however, the strong results presented here crucially rely on the coupling constructions  in Section \ref{sec:const} and \ref{sec:RCM}  available only under \eqref{cond:sharp_vv}.
\end{rem}

\section{Main results and conditions}\label{sec:results}
In this section, we present the three main results of the paper: a Poisson approximation result in Theorem~\ref{thm:1}, a quantitative CLT in Theorem \ref{thm:2}, and a qualitative CLT in Theorem \ref{thm:3}. All results hold for Gibbs models of type (L) and (P). In particular, $\Xx=\R^d\ti \M$. Each result will have some additional assumptions on the Gibbs process as well as on the functional. Examples of functionals where  the conditions can be verified are given in Section \ref{sec:exa}.

Both Theorem \ref{thm:1} and \ref{thm:2} involve sums of stabilizing scores. We first give the common definitions. A \emph{score function} is a measurable map $g:\Xx \times \bN \to [0, \ff)$. Given a score function, we will consider the random measure
\begin{align}
	\Xi := \Xi[\XX]& :=\sum_{\xx \in \XX } g(\xx, \XX)\,\delta_\xx, \label{eqn:poixidefXi1}
\end{align}
where $\de_\xx$ denotes the Dirac measure at $\xx$. The  \emph{intensity measure} of $\Xi $ is
\begin{align}
	\label{eq:La}
	\La(A) := \E[\Xi(A)], \qquad A \in \B.
\end{align}

 We assume that $g$ is \emph{locally stabilizing}, i.e.\ that there is a \emph{stabilization radius}  $R(\xx,\vp)$ such that for all $\vp \in \Nlf $ and $\xx \in \vp$, we have
$$g(\xx, \vp) = g\big(\xx, \vp \cap B_{R(\xx, \vp)}(\xx) \big).$$
The stabilization radius must satisfy the stopping property  
\begin{equation} 
	\label{eq:stop}
	\big\{R(\xx,\varphi)\le r\big\} = \big\{R(\xx,\varphi \cap B_r(\xx_{}))\le r\big\}
\end{equation}
for all $r>0$.
Finally, we impose that the score function is \emph{exponentially stabilizing} when applied to the Gibbs process $\XX$, i.e.\ the stabilization radius has exponential tails in the sense that there are $c_{\ms{es}, 1},c_{\ms{es}, 2}> 0 $ such that
\begin{align}
	\sup_{\substack{1\le k\le 5\\x_1, \dots, x_k\in \R^d}} & \int_{\M^k }\hspace{-.2cm} \P\big(R((x_1,m_1),\XX\cup \{(x_1,m_1),\ldots,(x_k,m_k)\}) > r\big) \mu(\d m_1) \cdots \mu(\d m_k)
  \le c_{\ms{es}, 1}\exp(-rc_{\ms{es}, 2}).\label{eq:s2}
\end{align}
We further assume that $g$ satisfies the following \emph{hereditary} property. For all $\xx,\yy\in \Xx$ and  $\vp \in \Nlf$,
\begin{align} \label{eq:here}
g(\xx,\vp\cup \{\xx\})=0 \quad \Longrightarrow \quad g\big(\xx,\vp \cup \{\xx,\yy\}\big)=0. 
\end{align}

 We will assume that $\Xi$ is observed in a window $Q_n := [-n/2, n/2]^d\times \M$ for some $n\ge 1$ such that $\la(Q_n) = \a n^d$. This window shape is chosen for convenience. The proof techniques can easily be adapted to other window shapes, see \cite{hos,penrose}, but one has to be more careful about boundary effects. 

%
%
\subsection{Poisson approximation}
\label{ss:thm1}

If the score function $g$ takes values in $\{0,1\}$, then the random measure $\Xi$ in \eqref{eqn:poixidefXi1} becomes a new point process, which is a thinning of $\X$. 
Our Poisson approximation result gives an upper bound on the distance between $\Xi\cap Q_n$ and a Poisson process on $Q_n$ having the same intensity measure $\La$ as $\Xi$. 
We formulate our result in the \emph{Kantorovich-Rubinstein distance} between two point processes $\xi$ and $\zeta$ on $Q_n$, defined by
\begin{align}
	\label{eq:ddkr}
	{d_{\mathsf{KR}}}(\xi,\zeta):=\sup_{h \in \text{Lip}} \big| \E[h(\xi)]- \E[h(\zeta)]\big|.
\end{align}
Here, $\text{Lip}$ is the class of functions $\Nlf_{Q_n} \to \R$ which are 1-Lipschitz with respect to the total variation distance on $\Nlf_{Q_n}$ (considered as a space of counting measures) given by
$$\dtv(\vp,\psi) = \sup_{A\in\BB,A\su Q_n} |\vp(A)-\psi(A)|,\qquad \vp,\psi\in \Nlf_{Q_n}.$$
Having introduced this terminology, we now state the Poisson approximation result.

%
%
\bet[Poisson approximation for dependent thinning of Gibbs processes] 
\label{thm:1} 
Let $\XX$ be an infinite-volume Gibbs point process of Model (L) or (P) satisfying conditions \eqref{cond:sharp} and \eqref{cond:exponential_sim} or \eqref{cond:sharp_vv} and \eqref{cond:sharp_v}, respectively. Assume $\M$ is locally compact and
 that $g_n$ is a $\{0, 1\}$-valued exponentially stabilizing score function, i.e.\ satisfying  \eqref{eq:stop},  \eqref{eq:s2}, and \eqref{eq:here}. Let $\Xi_n$ be the corresponding point process \eqref{eqn:poixidefXi1} and let $\mc M_n$ be a Poisson process on $\Xx$ with the same intensity measure $\La_n$ given by \eqref{eq:La}.  
Then,
\begin{align*}
&d_{KR}(\Xi_n \cap Q_n, \mc M_n \cap Q_n) \le c_{\ms{poi}} \la(Q_n)^{-1}+ 2\int_{Q_n}\int_{Q_n} \one \big\{S_{n,x} \cap S_{n,y} \ne \es\big\} \La_n(\d x) \La_n( \d y)\\
&\,+2\int_{Q_n}\int_{Q_n} \one \big\{S_{n,\xx} \cap S_{n,\yy} \ne \es\big\}  \E \Big[g_n(\xx,\XX\cup\{\xx,\yy\}) g_n(\yy,\XX\cup\{\yy\}) \k(\xx,\XX\cup\{\yy\}) \k(\yy,\XX)\Big] \la(\d \xx)\la (\d \yy),\nonumber
\end{align*}
where $S_{n,\xx}=B_{3c_{\ms St}\log n}(\xx)$ with $c_{\ms St} :=14\max(c_{\ms{es},2},c_{\ms{v},2},c_{\ms{D},2})$ and  $c_{\ms{poi}}=c_{\ms{poi}}(d, \a,  c_{\ms{es},1},c_{\ms{v},1},c_{\ms{D},1})$. 
\ent

While the right-hand side of the bound in Theorem \ref{thm:1} is more involved than the corresponding result for Poisson approximation of Poisson functionals (see \cite{BSY202}), it exhibits a similar structure. The first term puts a bound on long range effects due to exponential stabilization and spatial interactions in the Gibbs process. Noting that $\La_n(\d \xx) =\E[g_n(\xx,\XX\cup\{\xx\})\k(\xx,\XX)]\la(\d \xx)$, the first integral on the right-hand side controls the first-order contributions of the score function $g_n$, while the last term controls its local second-order contributions. More specifically, if the second integral vanishes asymptotically as $n \to \infty$, this implies that the local influence of an atom $\xx \in \XX$ with $g_n(\xx,\XX)=1$ on the score value $g_n(\yy,\XX)$ of other atoms $\yy \in \XX$ with $S_{n,\xx} \cap S_{n,\yy} \neq \es$ is asymptotically negligible.

The proof reveals that the $c_{\ms{poi}} \la(Q_n)^{-1}$ term in Theorem \ref{thm:1} can be made to be of order $\la(Q_n)^{-m}$, $m\ge 1$, by choosing $c_{\ms St}$ larger. However, it is usually the two integrals that determine the rate when $n\to \ff$.
In Section \ref{ss:pois}, we discuss an example where the rate of the integrals can be determined. Theorem \ref{thm:1} generalizes the Poisson approximation results from \cite{dp} in several key ways. First, we extend the class of Gibbs models where, instead of assuming subcriticality of an underlying random geometric graph as in Model (L), we impose the weaker condition of exponential decay for the underlying RCM, as specified in  \eqref{cond:sharp_vv}. Second, we relax the stabilization condition on $g_n$, which in \cite{dp} was required to have a deterministic stabilization radius.

%
%
\subsection{Quantitative CLT for exponentially stabilizing functionals}
\label{ss:thm2}
Our second main result is a quantitative normal approximation result in the spirit of \cite{chen} for functionals of the form $\Xi(Q_n)$. In this section,  score functions are no longer restricted to be $\{0, 1\}$-valued. 
We impose the following {bounded moment condition}. 
\begin{align} 
	\label{eq:m2}
	\sup_{n \ge 1}\sup_{k\le 5} \sup_{x_1, \dots, x_k \in \R^d}	 \E\int_{\M^k} \big[g_n((x_1, m_1),\XX\cup \{(x_1,m_1), \dots, (x_k, m_k)\})^6\big]\mu(\d m_1)\cdots \mu(\d m_k)=:c_{\ms m}<\ff.
\end{align}
The results are given in terms of the \emph{Kolmogorov distance} which is given for two real-valued stochastic variables $X$ and $Y$ as
$$\dk(X,Y) = \sup_{u\in \R} |\P(X\le u) - \P(Y\le u)|.$$

%
%
\begin{theorem}[Quantitative normal approximation of Gibbsian score sums]
\label{thm:2}
	Let $\XX$ be an infinite-volume Gibbs point process from Model (L) satisfying \eqref{cond:sharp} and \eqref{cond:exponential_sim} or from Model (P) satisfying \eqref{cond:sharp_vv}  and \eqref{cond:sharp_v}. Assume that $g_n$ is an exponentially stabilizing score function satisfying \eqref{eq:stop}, \eqref{eq:s2},  \eqref{eq:here}, and \eqref{eq:m2}.  Then, for $n\ge 1$, $\Xi_n(Q_n)$ defined at \eqref{eqn:poixidefXi1} satisfies
\begin{align*}
	\dk\bigg({\f{\Xi_n(Q_n)-\La_n(Q_n)}{\sqrt{\Var(\Xi_n(Q_n))}} }, N(0,1) \bigg) \le c_{\ms{norm}}\f{\la(Q_n)(\log \la(Q_n))^{2d}}{\Var({\Xi_n(Q_n)})^{3/2}},
\end{align*}
where $N(0,1)$ is a standard normal variable and $c_{\ms{norm}}=c_{\ms{norm}}(d, \a, c_{\ms m}, c_{\ms{es},1},c_{\ms{es},2},c_{\ms{v},1},c_{\ms{v},2},c_{\ms{D},1},c_{\ms{D},2})$. 
\end{theorem}
This result is comparable to \cite[Thm. 2.4]{hos} with the main novelty being the generalization to larger classes of Gibbs processes. The logarithmic factor in the bound is usually found in the Gibbs literature  \cite{gibbs_limit,gibbsCLT}, but we remark that it can be avoided in the Poisson case \cite{Mehler} and  it is  open whether this is also true in the Gibbs case. 
%
%
\subsection{Qualitative CLT for weakly stabilizing functionals}
\label{ss:thm3}

We finally give a qualitative CLT for a Gibbs process $\XX$ on $\R^d \times \M$ with translation-invariant PI. We consider a  functional $H: \Nlf_0 \to \R$ which is
 \emph{translation-invariant}, i.e.\ $H(\vp+y)=H(\vp)$ for all $y\in \R^d$, $\vp \in \Nlf_0$.  For $n\ge1$, we let $H_n$ denote the functional
\begin{equation}\label{eq:H_n_Gibbs}
	H_n(\vp) = H(\vp \cap Q_n), \qquad \vp \in \Nlf.
\end{equation}

 The CLT requires a \emph{weak stabilization} condition. Let $Q_{z,n}:=Q_n+z$ for $z\in \R^d$. We require that almost surely,
\begin{equation}\label{cond:weak}
	H(\XX)-H(\XX\sm Q_l) := \lim_{n\to \infty} \big(H(\XX\cap Q_{w_n,m_n})-H((\XX\cap Q_{w_n,m_n})\sm Q_l)\big)
\end{equation}
exists for all $l\ge 0$ and all sequences $m_n\in \N$ and $w_n\in \Z^d$ such that $\R^d=\bigcup_{k\ge 1} \bigcap_{n\ge k} Q_{w_n,m_n}$. Note that the limit is necessarily independent of the sequence $Q_{w_n,m_n}$.

We also need a moment condition. Since the proof relies on a general CLT for sums of martingale differences, this condition is stated in terms of these martingale differences.  Let $Z_{n+\sqrt{n}}$ be the set of lattice points $z\in \Z^d$ such that $Q_{z,1}$ intersects the window $Q_{n+\sqrt{n}}$. Then, $Q_{n+\sqrt{n}}\su \bigcup_{z\in Z_{n+\sqrt{n}}} Q_{z,1}$. Order $Z_{n+\sqrt{n}}$ lexicographically as $z_1,\dots,z_{k_n}$. Let $\FF_{0,n}$ be the trivial $\s$-algebra and $\FF_{k_n,n}=\sigma(\X)$. Moreover,  define 
$$ \FF_{i,n}=\FF_{z_i}:= \s\Big( \X \cap \bigcup_{{z\in \Z^d, z\preceq z_i}} Q_{z,1}\Big),\qquad 0<i<k_n,$$
 where $\preceq$ denotes lexicographical ordering on $\Z^d$. Then, $\FF_{i,0} \su \FF_{i,1}\su \dotsm \su \FF_{k_n,n}$ is an increasing sequence of $\s$-algebras containing information about $\X$ in increasing domains of $\Xx$. In particular, $\E\big[H_n(\XX) | \FF_{i,n}\big]$ define a martingale for fixed $n$. We require the following  moment condition on the martingale differences.
\begin{equation}\label{cond:cond_moment}
		\sup_{n\ge 1}\sup_{z_i\in Z_{n+\sqrt n}} \E\big[\big(\E\big[H_n(\XX) | \FF_{i,n}\big] - \E\big[H_n(\XX) |\FF_{i-1,n}\big] \big)^4\big] < \infty.
\end{equation}

%
%
\bet[CLT for translation-invariant functionals]
\label{thm:3}
Let $\XX$ be an infinite-volume Gibbs point process   of Model (L)  or Model (P)  with translation-invariant PI. Let $H$ be a translation-invariant functional and define $H_n$ as in \eqref{eq:H_n_Gibbs}. Assume that $H$ is weakly stabilizing \eqref{cond:weak} and the conditional moment condition \eqref{cond:cond_moment} holds. 
Then, the limit
\begin{equation*}
	\De_0 = \lim_{n\to \infty } \left(\E\big[ H_n(\X )| \X\cap \FF_{0} \big] - \E\big[ H_n(\X ) | \X\cap \FF_{0_-}\big]\right)
\end{equation*}
exists in $L^2$, where $0_-$ is the point before 0 in the lexicographical ordering of $\Z^d$. Let $\s^2:= \E[\De_0^2]$. Then,
$
\la(Q_n)^{-1}\Var(H_n(\X) ) \to \s^2
$
and
\begin{equation*}
	\la(Q_n)^{-1/2}\left(H_n(\X) - \E[H_n (\X)] \right) \to N(0,\s^2)
\end{equation*}
in distribution.
\ent

The weak stabilization condition is generally weaker than exponential stabilization, see the discussion in \cite{hos}. Moreover, note that we do not rely on any assumptions on the tails of the spatial cluster size distribution in the associated percolation model. The price for this is that the moment condition is more involved. To verify it, it will often be convenient to make assumptions on the spatial cluster sizes. More easily checkable conditions that imply \eqref{cond:cond_moment} are given in \cite{hos} under the additional assumption of exponential tails of cluster sizes \eqref{cond:sharp} or \eqref{cond:sharp_vv}, see Section \ref{sec:bet} for a discussion of  the general case. Positivity of the limiting variance $\sigma^2$ is not guaranteed, see \cite{hos} for an example of how to verify this.


	%
	%
\section{Examples}
\label{sec:exa}
In this section, we present specific examples of functionals to which our main results, Theorems \ref{thm:1}, \ref{thm:2}, and \ref{thm:3}, apply. In \cite{hos}, we discussed a number of examples for the case of finite-range interaction. Since we succeeded in proving Theorems \ref{thm:1}, \ref{thm:2} and \ref{thm:3} under essentially the same stabilization and moment conditions as in the finite-range case from \cite{hos},  also the examples carry over. To avoid redundancy, we keep the presentation here short and refer to \cite{hos} for details. One important difference is that \cite{hos} does not investigate the question of Poisson approximation, i.e., our Theorem \ref{thm:1}. Hence, we will be more detailed here. 

An additional assumption in many of the examples in \cite{hos} is {\em Poisson-likeness}. Following \cite{gibbs_limit}, a Gibbs process $\XX$ on $\R^d$ is called Poisson-like if  (i) it is stochastically dominated by a homogeneous Poisson point process, and (ii) there exist $c_{\ms{PL}}, r_1>0$ such that for all $r \ge r_1,\, x \in \R^d$ and $ \vp \in \Nlf_{B_r(x)^c}$, it holds that
\begin{equation}\label{eq:PLdef}
\P(\XX_{B_r(x)}=\es \mid \XX_{B_r(x)^c}=\vp) \le e^{-c_{\ms{PL}}r^d}.
\end{equation}
Poisson-likeness was shown in \cite[Proposition A.3]{hos} under a very mild condition for processes with finite interaction range. For Model (P) with infinite interaction range, Poisson-likeness does not hold in general because the left hand side of \eqref{eq:PLdef} can typically be made arbitrarily small by adding more points to $\vp$. This is solved in \cite{gibbs_limit} by assuming a hard-core radius $r_0>0$ such that $v(x,y)=\ff$ whenever $\|x-y\|\le r_0$. This ensures that $\vp$ cannot be too dense, and Poisson-likeness was shown in \cite[Lemma 3.3]{gibbs_limit}  under this condition. A full study of when Poisson-likeness holds is beyond the scope of this article, and we will leave it as an assumption in the examples below.

%
%
\subsection{Example for Theorem \ref{thm:1}: Isolated points}
\label{ss:pois}
Let $\XX$ be a Poisson-like infinite-volume Gibbs point process of Model (P) on $\R^d$. We consider the thinning $\Xi_n$ of  $\XX\cap Q_n$   consisting of those points $x \in \XX\cap  Q_n$  for which the distance to the nearest neighbor in $\XX$ exceeds $u_n(x)$ for some measurable function $u_n:Q_n\to [0,\ff)$. The process $\Xi_n$ can also be understood as a Mat\'ern type I thinning of $\Xi_n$ (see \cite[Section 10]{dp}). We can write $\Xi_n$ in the form \eqref{eqn:poixidefXi1} by letting $$g_n(x,\XX)=\one\{(\XX \sm \{x\}) \cap B_{u_n}(x)=\es\}.$$
Choose $u_n$ so large that
\begin{align} \label{exam:nn}
 \sup_{x \in Q_n} \P(\XX \cap B_{u_n(x)}(x)=\es) \le \frac{\b_2}{\la(Q_n)}.
\end{align}
for some positive constant $ \b_2$. This is clearly possible because of Poisson-likeness \eqref{eq:PLdef}.

The hereditary property \eqref{eq:here} obviously holds. Moreover, $R(x,\XX)= d(x,\X\sm\{x\})$, where is $d(x,\X\sm\{x\})$ is the distance from $x$ to the closest point in $\X\sm \{x\}$, is a  stabilization radius independent of $n$ satisfying the stopping condition  \eqref{eq:stop}.  The stabilization condition  \eqref{eq:s2} follows from Poisson-likeness of  $\XX$, since
$$
\P(R(x,\X)> r) \le \P(\XX \cap B_{r}(x)=\es) \le e^{-c_{\ms{PL}} r^d},\quad x \in Q_n.
$$

We now apply Theorem~\ref{thm:1}. Note that \eqref{exam:nn} implies
$$
\E \Big[ g_n(x,\XX \cup \{x\})  \Big] =\P(\XX \cap B_{u_n(x)}(x)=\es) \le \frac{\beta_2}{\la(Q_n)},\quad x \in Q_n.
$$
 Since $\k \le 1$, the first integral term on the right-hand side in Theorem \ref{thm:1} is bounded by
\begin{align*}
& \iint \one \big\{S_{n,x} \cap S_{n,y} \ne \es\big\} \La_n^2(\d x,\d y) \\
&\quad =\a^2 \int_{Q_n} \int_{Q_n} \one \big\{S_{n, x} \cap S_{n, y} \ne \es\big\} \E \Big[ g_n(x,\XX \cup \{x\}) \k(x,\XX)  \Big] \E \Big[ g_n(y,\XX \cup \{y\}) \k(y,\XX)  \Big] \d y \d x\\
&\quad \le \f{\a^2\b_2^2}{\la(Q_n)^{2}} \int_{Q_n} \int_{Q_n} \one \big\{S_{n,x} \cap S_{n,y} \ne \es\big\} \d x \d y \in O\Big( \f{(\log \la(Q_n))^d}{\la(Q_n)}\Big).
\end{align*}

Before bounding the second integral term from Theorem \ref{thm:1}, we note that
\begin{align*}
&\E \Big[g_n(x,\XX\cup\{x,y\}) g_n(y,\XX\cup\{ y\}) \k(x,\XX\cup\{y\}) \k(y,\XX)\Big]\\
&\quad \le \P\Big(\XX \cap B_{u_n(x)}(x)=\es, \XX \cap B_{u_n(y)}(y)=\es\Big)\one\{ \|x-y\|\ge \max\{u_n(x),u_n(y)\} \}\\
&\quad \le  \P\Big(\XX \cap B_{u_n(x)}(x)=\es, \XX \cap B_{\f{u_n(y)}{2}}\Big(y+\f{u_n(y)}{2} \f{y-x}{\|y-x\|}\Big)=\es\Big)\\
& \quad =  \P\Big(\XX \cap B_{u_n(x)}(x)=\es \Big) \P\Big(\XX \cap B_{\f{u_n(y)}{2}}\Big(y+\f{u_n(y)}{2} \f{y-x}{\|y-x\|}\Big)=\es \Big| \XX \cap B_{u_n(x)}(x)=\es\Big).
\end{align*}
Since the two balls considered in the last line are disjoint, we obtain from Poisson-likeness  that the conditional probability in the last display is bounded by 
$\exp(-c_{\ms{PL}}(u_n(y)/2)^d)$. 
Since $\X$ can be realized as a thinning of $\PP$ (see e.g. Proposition \ref{prop:RCM_ff} below),
$$
e^{-\a\k_du_n(y)^d} = \P(\PP \cap B_{u_n(y)}(y)=\es) \le  \P(\XX \cap B_{u_n(y)}(y)=\es) \le \frac{\beta_2}{\la(Q_n)}. 
$$
It follows that
$\exp(-c_{\ms{PL}}(u_n(y)/2)^d)\le (\b_2/\la(Q_n))^{c_{\ms{PL}}/(2^d\a\k_d)}.$
This gives
\begin{align*}
&\int_{Q_n} \int_{Q_n} \one \big\{S_{n,x} \cap S_{n,y} \ne \es\big\}  \E \Big[g_n(x,\XX\cup\{x,y\}) g_n(y,\XX\cup\{x,y\}) \k(x,\XX\cup\{y\}) \k(y,\XX\cup \{x\})\Big] \d x \d y\\
&\quad \le \Big(\f{\b_2}{\la(Q_n)}\Big)^{\frac{c_{\ms{PL}}}{2^d\a\k_d}} \int_{Q_n} \int_{Q_n} \one \big\{S_{n,x} \cap S_{n,y} \ne \es\big\}  \P(\XX \cap B_{u_n(x)}(x)=\es) \d x \d y \in O\Big(\f{(\log \la(Q_n))^d}{\la(Q_n)^{c_{\ms{PL}}/(2^d\a \k_d)}}\Big).
\end{align*}
Hence, we conclude from Theorem \ref{thm:1} that
$$
{d_{\ms{KR}}}(\Xi_n \cap Q_n, \mc M_n \cap Q_n) \in O\Big(\f{(\log \la(Q_n))^d}{\la(Q_n)^{c_{\ms{PL}}/(2^d\a \k_d)}}\Big).
$$

%
%
\subsection{Example for Theorem \ref{thm:2}: Edge lengths in $k$-nearest neighbor graph}
\label{ss:knn}
In the $k$-nearest neighbor (kNN) graph with vertices in $\R^d$, connections are drawn from any given vertex to the $k$ other vertices that are closest in Euclidean distance. More precisely, given a locally finite vertex set $\vp \su \R^d$ an edge is drawn between two points $x, y \in \vp$ if $x$ is one of the $k$-nearest neighbors of $y$ or the other way around. Normal approximation for the total edge length of kNN graphs with vertices on a Gibbs point process  was considered in \cite[Section 3.1]{hos} and \cite[Theorem 4.1]{CX22}.  Formally, we can put this example in the setting of score functions by setting $g(x, \vp)$ to be one half of the length of all kNN edges incident to $x$. Then, we consider the total length given as 
\begin{align*}
	H^n(\vp) := \sum_{x \in \vp \cap Q_n} g(x,\vp).
\end{align*}
To apply our normal approximation result Theorem \ref{thm:2}, we need to verify the stabilization and moment conditions   \eqref{eq:s2} and \eqref{eq:m2}. These conditions are identical to the finite-range case from \cite{hos}, so assuming Poisson-likeness, no additional work is needed. The exponential decay of the percolation probabilities was automatic in the setting of \cite{hos}, but  this is an assumption in general. In the case of Model (P), we refer to \cite{mikhail} for a more detailed discussion.

%
%
\subsection{Example for Theorem \ref{thm:2} with marks: Edge counts in Gilbert graphs}
\label{ss:gilb}
Example \ref{ss:knn} considers a situation of quantitative normal approximation involving unmarked Gibbs processes. Now, we discuss an example illustrating the marked case. Here, the mark space is given by $\M := [0, \infty)$ and the marks are interpreted as random radii. Given a locally finite set $\varphi \su \Xx$, we connect $\xx = (x, m)$ and $\xx' = (x', m')$ by an edge if $\|x - y\| \le m + m'$. Then, we consider the functional given by the number of edges in the associated graph. Similarly as in Example \ref{ss:knn},  we can define the score function $g(\xx, \vp)$ to be the number of neighbors in the graph that have a smaller mark than $\xx$. Then, we consider the total edge count given as 
\begin{align*}
	H^n(\vp) := \sum_{\xx \in \vp \cap Q_n} g(\xx,\vp).
\end{align*}
The stabilization radius can be chosen as twice the mark, i.e., $R\big((x, m), \vp\big) = 2m$. 
For Gibbs processes whose PI does not depend on the marks, the resulting Gibbs process will be an independently marked Gibbs process. Hence, we can treat the percolation conditions in the same way as in Example \ref{ss:knn}. Moreover, for the stabilization condition \eqref{eq:s2}, it is sufficient to assume an exponential tail of the original mark distribution. Finally, the moment condition \eqref{eq:m2} follows since the Gibbs process is dominated by a Poisson point process which has finite moments of any order.

\subsection{Example for Theorem \ref{thm:3}: Persistent Betti numbers}
\label{sec:bet}
Persistent Betti numbers $\beta_q^{r,s}$, $0<r<s$ and $q\in \N_0$,  are summaries of point processes coming from topological data analysis. We will not give the details here, but see \cite{edHar} for a general introduction to topological data analysis and \cite{shirai} for the definition of persistent Betti numbers. CLTs for Betti numbers of Poisson processes were first given in \cite{shirai,yogeshAdler2}. In \cite[Sec. 3.1]{hos}, this was extended to the Gibbs processes in Example \ref{ex:hos}. We can generalize these results to other processes of Model (L)  as well as Model (P). The weak stabilization condition \eqref{cond:weak} holds deterministically \cite{shirai}, so  Theorem \ref{thm:3} applies if we can verify the moment condition \eqref{cond:cond_moment}. This was done in \cite{hos} using only the exponential tails of the cluster diameters for the Boolean model. Carefully inspecting the proof, one finds that it suffices  that $\P_o(o\lrsa Q_n^c) $ is of order $O(n^{-20d-5-\eps})$, where $o$ denotes the origin. Hence, the proof also holds in the case of pairwise interactions when the clusters of the RCM satisfy the analogous criterion. Poisson-likeness was not an assumption, but the variance positivity proof in \cite{hos} assumed finite interaction range.

	%
%

\section{Disagreement coupling}
\label{sec:const}
 The key technical tool of our proofs is disagreement coupling. This technique explicitly constructs a Gibbs process as a thinning of its dominating Poisson process. Moreover, it allows Gibbs processes with different boundary conditions to be coupled to the same Poisson process in such a way that they disagree only on specific clusters of the corresponding percolation model. This section is devoted to setting up the framework for Model (L) while Section \ref{sec:RCM} considers the generalization to Model (P).
 
We rely on the basic disagreement coupling framework from \cite{dp} for Gibbs processes on bounded domains, which we recall in Section \ref{ss:fw}.   Section \ref{ss:pert} gives bounds on the disagreement between Gibbs processes with locally perturbed PIs. In Section \ref{ss:emb}, we show how to extend disagreement coupling to unbounded domains. Most proofs are very similar to those from \cite{hos} for the models in Example \ref{ex:hos} with the only difference being that $\R^d$ must be replaced by $\Xx$. Hence, the proof for those direct analogs are omitted. 
Proofs of new results and proofs where the arguments substantially differ from \cite{hos} are deferred to Section \ref{sec:proofs_const}. While our applications are to marked Euclidean spaces, the theory is set up in an arbitrary Polish space $\Xx$  and specialized to the marked Euclidean case only when necessary.

%
%

%
%
\subsection{Couplings with a Poisson process}
\label{ss:fw}
%
%

In this section, we consider the construction of a Gibbs process $\XX(Q,\psi)$ on a bounded domain $Q\in \BB_0$ with boundary conditions $\psi \in \Nlf_{\Xx\sm Q}$ given by thinning  a Poisson process with intensity measure $\lambda$.  We first review the basic version of the algorithm from \cite[Sec. 5]{dp} in Section \ref{sss:spe}. As in \cite{hos}, this will be the building block for the more flexible disagreement coupling algorithm  introduced in Section \ref{sss:dc}.

Henceforth, for  $B\in\BB$, we let $B^*=B\ti [0,1]$. We let $\Pds$ be a  Poisson point process on $\mathbb{X}^* $ with intensity measure $ \lambda \times \Leb_{[0,1]}$, where $\Leb_{[0,1]}$ denotes Lebesgue measure on $[0,1]$. We write $\Pd$ for the projection of $\Pds$ onto the $\mathbb{X}$ coordinate, which is a Poisson process on $\Xx$ with intensity measure $\la$. More generally, we use $\vp^*$ to denote elements of $\Nlf_{\Xx^*}$ and $\vp$ denotes the projection of $\vp^*$ to the $\Xx$-coordinate. For $B\in\BB$, we let $\vp^*_B=\vp^*\cap B^*$ and $\vp_B=\vp\cap B$.

%
%
\subsubsection{The standard Poisson embedding}
\label{sss:spe}

Let $Q \in \BB_0 $ and $\psi \in \Nlf_{\Xx \sm Q}$. In particular, $Q$ has finite volume $\lambda(Q)<\ff$. We recall the definition of the thinning map  $T_{Q, \psi}: \Nlf_{Q^*} \to \Nlf_Q$ from \cite{dp} in Definition \ref{def:standard} below. For $\vp^* \in  \Nlf_{Q^*}$, the set $T_{Q, \psi}(\vp^*)$ will be a thinning of $\vp$ in the sense that  $T_{Q, \psi}(\vp^*)\su \vp$. The last coordinate of a point in $\vp^* $ is only used to decide whether or not to keep the point in the thinning. 
The key property of $T_{Q, \psi}$ is that when applying it to a Poisson process $\Pds_Q$ we obtain the Gibbs process $\X(Q,\psi)$, that is, $T_{Q, \psi}(\Pds_Q)\stackrel{d}{=} \X(Q,\psi)$, where $\stackrel{d}{=}$ denotes equality in distribution. We refer to $T_{Q, \psi}(\Pds_Q)$ as the \emph{standard Poisson embedding}, since it provides an embedding $\XX(Q, \psi)\su \Pd_Q$. We call it standard because it will be the main ingredient for other Poisson embeddings introduced in later sections under the name disagreement couplings. 

To define $T_{Q, \psi}$, we need an injective measurable map  $\iota\co Q \to (0,1)$, which is guaranteed to exist in any standard Borel space. We  write $T_{Q, \psi, \iota}$ when we wish to emphasize that $T_{Q, \psi}$ depends on $\iota$. The map $\iota$ defines a total ordering $\le_{\iota}$ on $Q$. We sometimes refer to the map $\iota$ itself as an ordering.
To construct $T_{Q, \psi}(\vp^*)$, we go through the points of $\vp$ in the order given by ${\iota}$ and decide for each point whether or not to keep it in the thinned process. For  $\xx \in  Q $, let
\begin{align*}
&	Q_{(-\ff,\xx)} {}:=\{ \yy \in Q \co\yy <_\iota\xx\},\quad Q_{(\xx, \ff)} {}:=\{ \yy \in Q \co\yy >_\iota\xx\}, \\ 	& Q_{(-\ff,\xx]} {}:=\{ \yy \in Q \co \yy \le_\iota\xx\}, \quad Q_{[\xx, \ff)} {}:=\{ \yy \in Q \co \yy \ge_\iota\xx\}.
\end{align*}
For $\xx \in  Q$ and $\psi \su Q_{(-\ff,\xx)}$, define the retention probability
\begin{align}
	p(\xx, Q, \psi) := \k(\xx, \psi) \f{Z_{Q_{(\xx, \ff)}}\big(\psi \cup \{\xx\}\big)}{Z_{Q_{(\xx, \ff)}}(\psi)}.\label{def:p}
\end{align}
It is shown in \cite[Lem. 4.1]{dp} that under Assumption \eqref{cond:PIbound}, $0\le p(\xx, Q, \psi) \le 1$. 

We now have all the ingredients needed for an algorithmic definition of the standard Poisson embedding.

%
%
\bede[Standard Poisson embedding]
\label{def:standard}
Let $\vp^*\in\Nlf_{Q^*}$. Then,  $T_{Q, \psi}(\vp^*)= \{\xx_1,\xx_2, \dots\}\su \vp$ is the point set obtained by the following algorithm.
\begin{itemize}
	\item[$\bullet$] To define the first point $\xx_1 :=\xx_1(\vp^*)$, consider the set 
	\begin{equation}\label{eq:defx1}
	\big\{\xx\in Q \co (\xx, u) \in \vp^*  \text{ and } u \le p(\xx, Q, \psi) \big\}.
	\end{equation}
	We choose $\xx_1$ to be the smallest element in \eqref{eq:defx1} with respect to the $\iota$-ordering. If \eqref{eq:defx1} is empty, then the construction terminates with $T_{Q, \psi}(\vp^*) = \es$. 
	\item[$\bullet$] Otherwise, we proceed inductively. Suppose, we have defined the points $\xx_1,\ldots,\xx_k $. This will be the restriction of $T_{Q, \psi}(\vp^*)$ to $Q_{(- \ff,\xx_k]}$. It remains to define $T_{Q, \psi}(\vp^*)$ on $Q_{(\xx_k, \ff)} $. We define $\xx_{k+1}=\xx_{k + 1}(\vp^*)$ to be the $\iota$-smallest element of the set 
	\begin{equation}\label{eq:admissible}
		\big\{\xx\in Q_{(\xx_k, \ff)}  \co (\xx, u) \in \vp^*  \text{ and } u \le p(\xx, Q, \psi \cup \{\xx_1, \dots,\xx_k\})\big\}.
	\end{equation}
	\item[$\bullet$] If \eqref{eq:admissible} is empty, the construction terminates  with $T_{Q, \psi}(\vp^*)=\{\xx_1,\xx_2, \ldots,\xx_k\}$. 
\end{itemize}
\ende

The key result of \cite[Thm. 5.1]{dp} is that the thinned process $T_{Q, \psi}(\Pds_Q) $ is the Gibbs process $\XX(Q, \psi)$. 

%
%
\bepr[Correctness of the Poisson embedding \cite{dp}]
\label{pr:emb}
Assume that the PI satisfies \eqref{cond:PIbound}. Let $Q \in \B_0$ and $\psi \in \Nlf_{\Xx\sm Q}$. Then, $T_{Q, \psi}(\Pds_Q)\su \Pd_Q$ is distributed as the Gibbs process $\XX(Q, \psi)$. 
\enpr

%
%
\subsubsection{Disagreement coupling in finite-volume domains}\label{sss:dc}

The standard Poisson embedding defines a coupling between $\XX(Q,\psi)$ and $\Pd_Q$ for any choice of boundary conditions $\psi$, but it is not clear how different choices of $\psi$ affect the embeddings of $\XX(Q,\psi)$. In this section, we modify the standard Poisson embedding to define \emph{disagreement coupling}, which provides more flexible ways of constructing Gibbs processes as Poisson thinnings. The name disagreement coupling refers to one of the main uses, namely obtaining a coupling between Gibbs processes $\XX(Q, \psi)$ and $\XX(Q, \psi')$ with differing boundary conditions $\psi,\psi' \in \Nlf_{\Xx\sm Q}$ such that the two Gibbs processes only  disagree near $\psi\Delta \psi'$, where $\Delta $ denotes the symmetric set difference.  We first describe the general framework for constructing disagreement couplings in Definition \ref{def:dc} and then introduce the specific version we are going to use in Definition \ref{ex:dc1}. The construction is the immediate generalization of the one in \cite{hos} to general spaces $\Xx$.

For the definition, we need the notion of a \emph{stopping set}. A stopping set $S$ on a Borel set $Q\su \Xx$ assigns to each  $\vp \in \Nlf_Q$ a Borel set  $S(\vp)\su Q$. We require that the map $(\xx,\vp) \mapsto \one_{\{\xx \in S(\vp)\}}$ is measurable and that $S$ has the property that $S(\vp)$ depends only on $\vp\cap S(\vp)$, i.e.\  $S(\vp) =  S(\vp\cap S(\vp)) $ for all $\vp \in \Nlf_Q$. For instance, $S(\vp)=Q_{(-\ff,\xx ]}$, where $\xx $ is the $\iota$-smallest point of $\vp$, is a stopping set. We write $S=S(\vp)$ when $\vp$ is clear from the context. A key property of a stopping set is the following, see e.g. \cite[Thm. A.3]{stoppingset}: If $S$ is a stopping set on $Q$ and $\PP_Q$ is a Poisson processes on $Q$ with intensity measure $\lambda$, then 
\begin{equation}\label{eq:stopping}
	\text{conditionally on } \PP_{S(\PP_Q)}, \,	\PP_{Q\sm S(\PP_Q)} \text{ is a Poisson process on $Q\sm S(\PP_Q)$ with intensity $\lambda_{Q\sm S(\PP_Q)}$.}
\end{equation}

%
%
To define disagreement coupling, we need the following ingredients:
\begin{itemize}
	\item[$\bullet$] A sequence $S_0^*, S_1^*, S_2^*, \dots $ of stopping sets on $Q^*$ of the form $S_n^*=S_n\times [0,1]$ such that for any $\vp^* \in \Nlf_{Q^*}$,  
	$$S_0^*(\vp^*)=\es\su S_1^*(\vp^*)\su S_2^*(\vp^*) \su \cdots  \quad \text{ and } \quad \bigcup_n S_n^*(\vp^*) = Q^*.$$ 
	\item[$\bullet$] A family of injective orderings  $\iota_n:=\iota_n^{\vp^*}:Q\sm S_n(\vp^*) \to \R $ defined for $n\ge 0$ and $\vp^*\in \Nlf_{Q^*}$ such that $\iota_n^{\vp^*}$ depends only on $\vp^*$ via $\vp^* \cap S_n^*(\vp^*)$. Moreover, we require
	\begin{equation}\label{eq:iota_n}
		\iota_n(\xx) < \iota_n(\yy )\quad \text{ for all } \quad\xx\in S_{n+1}\sm S_n,\, \yy \in Q\sm S_{n+1}. 
	\end{equation}
		We extend $\iota_n$ to a map $\tilde{\iota}_n:=\tilde\iota_n^{\vp^*}:Q  \to \R\cup\{-\ff\} $  by setting $\tilde{\iota}_n^{\vp^*}(\xx)=-\ff$ if $\xx \in S_n(\vp^*)$ and require that the map $(\xx, \vp^*) \mapsto \tilde{\iota}^{\vp^*}_n(\xx) $ is measurable as a function $Q \times \Nlf_{Q^*} \to  \R\cup\{-\ff\}$.
\end{itemize}
We think of the stopping sets as providing a step-wise exploration of $Q$, where the set  $S_{n}(\vp^*)\sm S_{n-1}(\vp^*)$ is explored in step $n$. 
In disagreement coupling, the point pattern   $\vp\in \Nlf_Q$ is thinned during this exploration, that is, the points $\vp \cap (S_{n}(\vp^*)\sm S_{n-1}(\vp^*))$ are thinned in step $n$. For the definition, recall that we  write $T_{Q, \psi, \iota}$ to highlight the dependence of the thinning map $T_{Q, \psi}$ on the ordering $\iota$.
\bede[Disagreement coupling]\label{def:dc}
Let $\psi \in \Nlf_{\Xx\sm Q}$. Disagreement coupling is a map $\tdc_{Q,\psi}: \Nlf_{Q^*} \to \Nlf_{Q} $.  The thinning $\tdc_{Q,\psi}(\vp^*) \su \vp$ of a point pattern $\vp^*$ is constructed inductively in $n$ by constructing the thinned point pattern $\xi_n:=\xi_n(\vp^*)$ on the set $S_n\sm S_{n-1}$ as follows: 
\begin{itemize}
	\item[1.] First, set  $\chi_1:=\chi_1(\vp^*):=T_{Q,\psi, \iota_0} (\vp^*)$. Let $\xi_1:=\chi_1\cap S_1$. 
	\item[2.] Assume inductively that we have constructed $\xi_i$ on $S_i\sm S_{i-1}$ for $i\le n$. The thinning of $\vp \cap S_n$ is $\xi_1\cup\ldots \cup \xi_n$. Next construct the thinning on $S_{n+1}\sm S_n$. For this, add $\xi_1\cup\ldots \cup \xi_n$ to the boundary conditions $\psi$ and apply the standard thinning map $ T_{Q\sm S_n,\psi \cup \xi_1 \cup \ldots \cup \xi_n,\iota_n}$ on $Q\sm S_n$  to obtain  $ \chi_{n+1}:= \chi_{n+1}(\vp^*):= T_{Q\sm S_n,\psi \cup \xi_1 \cup \ldots \cup \xi_n,\iota_n} (\vp^*_{Q\sm S_n})$. Set $\xi_{n+1}:=\chi_{n+1} \cap (S_{n+1}\sm S_n)$. 
	\item[3.] Finally, define $\tdc_{Q,\psi}(\vp^*) :=\bigcup_{n\ge 1} \xi_n $. 
\end{itemize}
\ende

Taking $S_n(\vp^*) =Q_{(-\ff,\xx_n]}$, where $\xx_n$ is the $n$th smallest point of $\vp$, we recover the standard embedding map from Definition \ref{def:standard}.  The following proposition shows that applying $\tdc_{Q,\psi}$ to a Poisson process  yields a Gibbs process generalizing Proposition \ref{pr:emb}. 

%
%
\bepr[Correctness of disagreement coupling]
\label{th:dc_prop}
Suppose the PI satisfies \eqref{cond:PIbound}. Let $Q \in \B_0$  and let $\psi \in \Nlf_{\Xx\sm Q}$. Then, $\tdc_{Q,\psi}(\Pds_Q)\stackrel{d}{=}\XX(Q, \psi)$. 
\enpr

\begin{proof}
	The proof is identical to the proof of \cite[Prop. 4.4]{hos} after replacing $\R^d$ by $\Xx$.
\end{proof}

 We want to relate the disagreement probabilities of Gibbs processes with different boundary conditions to percolation clusters. For this, we introduce the following cluster-based version of disagreement coupling, which is the version mainly used in this paper.

%
%
\bede\label{ex:dc1}[Cluster-based thinning] 
Assume the existence of a symmetric relation $\sim$ on $\Xx$ such that  \eqref{cond:PI_sim} holds. 
	Let $B\su \Xx \sm Q$ and $\psi\in \Nlf_{\Xx \sm Q}$. The \emph{cluster-based thinning map}, denoted by $T_{Q,B, \psi}^{\ms{cl}} : \Nlf_{Q^*} \to \Nlf_{Q}$, is the version of disagreement coupling arising from the following choice of stopping sets and orderings: 	Choose an initial ordering  $\iota: Q\to (0,1)$.
	 Let $S_1:=N(B)\cap Q$ be the deterministic set consisting of all points in $Q$ that are neighbors to a point in $B$ under the relation $\sim$.   Suppose $S_n$  is defined. We continue with one of the two cases below. 
	\begin{itemize}
\item[1.]	If $N(\vp \cap S_n)\cap (Q\sm S_n) \neq \es $, we take $S_{n+1} := S_n \cup N(\vp\cap S_n)\cap Q $. 
	\item[2.] Otherwise, let $S_{n+1}=S_n\cup Q_{(-\ff,\xx]}$ where $\xx$ is the $\iota$-smallest point in $(Q\sm S_n)\cap \vp$.  If no such point exists, we set $S_{n+1}=Q$.  
	\end{itemize}
	
	Define the orderings $\iota_n: Q \sm S_n \to (0,1)$ satisfying \eqref{eq:iota_n} by
	\begin{align*}&\iota_0(\xx)=\tfrac1{2} \iota(\xx) \one{\{\xx \in S_1 \}} + ( \tfrac1{2} \iota(\xx) + \tfrac1{2} ) \one{\{\xx\in Q\sm S_1 \}},\\
		&\iota_n(\xx)=\tfrac1{2} \iota(\xx) \one{\{\xx \in N(\vp \cap S_n)\cap ( Q \sm S_n) \}} + ( \tfrac1{2} \iota(\xx) + \tfrac1{2} ) \one{\{\xx\in Q\sm (S_n\cup N(\vp \cap S_n)  ) \}}, \qquad n\ge 1. 
	\end{align*}
\ende	
	
	In the cluster-based thinning, we thus take $S_{n+1}$ to be all neighbors of the so far explored points $\vp\cap S_n$. If all these neighbors have already been explored, we instead take $S_{n+1}\sm S_n$ to consist of all so far unexplored points up to and including the smallest point in $\S_n\sm \vp$ in the $\iota$-ordering. The stopping sets thus explore the clusters of $\GG(\vp)$ starting from the points of $\vp$ that are neighbors to a point in $B$. As long as we are in case 1., we explore the same clusters, while in case 2., we jump to a new cluster.  The order in which the remaining clusters are visited is determined by the $\iota$-smallest point in the cluster. 
	
	The advantage of the cluster-based thinning is that when we jump to a new cluster in 2., it means that $(\vp \cap S_n) \cup (\psi \cap B)$ is not neighbor to any point in $Q\sm S_n$. Hence, by Assumption \eqref{cond:PI_sim}, the retention probabilities \eqref{def:p} used in the remainder of the algorithm are not influenced by any of the points in $(\vp \cap S_n) \cup (\psi \cap B)$. Therefore, the boundary points $\psi \cap B$ only affect the thinning of clusters connected to $B$, and the  thinning of each new cluster does not depend on the thinning decisions made on the previous clusters.  
This means that when we change the boundary conditions locally in $B$, it only affects the thinning of the clusters of $\vp$ connected to a point in $B$. 

We summarize the main properties of the cluster-based thinning in Proposition~\ref{pr:rad_prop}.
The last statement of the proposition bounds the \emph{total variation distance} between two  Gibbs processes. The total variation distance between two random variables $X, Y$ with values in a common measurable space $S$ is defined by
\begin{align}
	\label{eq:dtv}
	\dtv(X, Y) := \sup_{f\co S \to [0, 1]}\big|\E[f(X)] - \E[f(Y)]\big|,
\end{align}
where the supremum is taken over all measurable functions $f\co S \to [0, 1]$.
%
%
\bepr[Properties of the cluster-based thinning]\label{pr:rad_prop}
Suppose \eqref{cond:PIbound} and \eqref{cond:PI_sim} holds. Let $Q \su \Xx$ with $Q\in\BB_0$ and let $B\su  \Xx \sm Q$ be a countable set.  The cluster-based thinning has the property that if $\psi,\psi'\in \Nlf_{\Xx\sm Q}$ with $\psi \De \psi' \su B $, then $T^{\ms{cl}}_{Q,B, \psi}(\Pds_Q)$ and $T^{\ms{cl}}_{Q,B, \psi'}(\Pds_Q)$ differ only on the clusters of $\GG(\Pd_Q) $  that have a point  that is neighbor to a point in $B$. In particular, for $A\su Q$
\begin{align*}
	\P(T^{\ms{cl}}_{Q,B, \psi}(\Pds_Q)\cap A \neq T^{\ms{cl}}_{Q,B, \psi'}(\Pds_Q)\cap A)& \leq \P(\Pd_A \lrsa_{\Pd_Q\cup B} B), \\
	\dtv(\XX(Q,\psi)\cap A , \XX(Q,\psi')\cap A) &\leq \P(\Pd_A\lrsa_{\Pd_Q \cup B } B).
\end{align*}

\enpr

We conclude this section by noting that disagreement coupling can also be used to show that \eqref{cond:perc2} implies uniqueness in distribution of infinite-volume Gibbs processes, see \cite[Thm 9.1]{betsch2}.

%
%
\subsection{Local perturbation of PIs}
\label{ss:pert}
In this section, we show a bound on how much local changes in the PI affect the Poisson embedding.  Several important consequences will be drawn from this lemma in later sections.  In the following, let $\XX(Q,\psi)$ and $\XX'(Q,\psi')$ denote Gibbs processes with PIs $\k$,  $\k'$ and boundary conditions $\psi,\psi'\in \Nlf_{\Xx \sm Q}$, respectively. A "prime" will be used to denote quantities associated with $\kappa'$, e.g. $T_{Q,\psi'}'$ denotes the standard thinning using the PI $\k'$. 

Lemma \ref{lem:distv} below bounds the disagreement probabilities between the first points in the respective Poisson embeddings in terms of the {total variation distance} betwen the corresponding Gibbs processes. 
Since the total variation distance has  an interpretation as the infimum of disagreement probabilities over all couplings, the lemma and its corollary below show that the Poisson embedding is not too far from optimal.

\bel [Disagreement probability of first point]
\label{lem:distv}
Let $A\su Q$ be elements of $\B_0$, $\psi,\psi'\in \Nlf_{\Xx\sm Q}$, and $\iota:Q\to (0,1)$ be injective with $\iota(\xx)<\iota(\yy)$ whenever $\xx\in A$ and $\yy \in Q\sm A$. Suppose that $\k, \k'$ are PIs satisfying \eqref{cond:PIbound}. 
Then,
$$
\P\big(\inf{}_{\iota} (T_{Q,\psi}(\Pds_Q) \cap A) \ne \inf{}_{\iota} (T'_{Q,\psi'}(\Pds_Q)\cap A)\big)\le (2+\lambda(A))\,\dtv(\XX(Q,\psi)\cap A,\XX'(Q,\psi')\cap A),
$$
where 
$\inf_\iota$  denotes the smallest point with respect to the ordering $\le_\iota$.
\enl

\bep
The proof is almost identical to the proof of \cite[Lem. 4.21]{hos} after replacing $\R^d$ by $\Xx$ and is hence omitted. 
\enp

For the next extension, assume that $\k$ and $\k'$ are closely related on some $P\su Q$ via a common PI $\k_0$ satisfying \eqref{cond:PIbound}, i.e.\ there are  $\eta,\eta' \in \Nlf_{\Xx\sm P}$ such that for all $\xx \in P$ and $\vp \in \Nlf$,	
\begin{align}\nonumber
	\k(\xx,\vp)&=\k_0(\xx,(\vp\cap P_1)\cup \eta)\\
	\k'(\xx,\vp) &= \k_0(\xx,(\vp\cap P_2) \cup \eta') 	\label{eq:k,k'}
\end{align}
where $P\su P_1,P_2\su Q$. In the statement below, there are no boundary conditions, as these can always be put as part of $\eta$ and $\eta'$, respectively, as noted in Section \ref{ss:Gibbs}. 

\bec[Disagreement probabilities for standard Poisson embedding]\label{cor:poiss_emb}
Suppose that $\k, \k'$ are of the form \eqref{eq:k,k'} and $\eta,\eta' \in \Nlf_{\Xx\sm Q}$. Let $A\su P$ be a Borel set. Let $\iota:Q\to (0,1)$ be an injective map such that $\iota(\xx)<\iota(\yy)$ whenever $\xx \in A$ and $\yy \in Q\sm A$.  Then, with $\bar\eta :=  \eta\Delta \eta'$,
\begin{equation*}
	\P\big(T_{Q,\es}(\Pds_Q)\cap A\ne T'_{Q,\es}(\Pds_Q) \cap A\big) \le 2\int_A \P\big(\xx\lrsa_{\Pd_\xx \cup \bar\eta} \Pd_{Q\sm P} \cup  \bar \eta \big)\lambda(\d\xx).
\end{equation*}
\enc

\bep
See Section \ref{sec:proofs_const}.
\enp

We can apply Corollary \ref{cor:poiss_emb} step-wise to give disagreement probabilities for perturbed PIs when using the cluster-based thinning. We will assume $\Xx = \R^d \times \M$ is marked Euclidean space and $Q=Q^0\times \M$. Define the initial ordering is as follows. We divide $\R^d$ into  cubes of diameter at most 1. We enumerate the cubes  intersecting $Q^0$ as $C_l^0$, $l=1,2,\ldots,L$ and let $C_l=C_l^0\times \M$. On $C_l$ we choose an ordering $\iota^l: C_l\cap Q \to (2^{-l-1},2^{-l})$ and let $\iota^0:\Xx \sm Q \to (0,2^{-L-1})$. We piece these orderings together to form an ordering $\iota:\Xx \to (0,1)$. We call this a \emph{cube-wise} ordering. 

The result is expressed in terms of probabilities of  two error events, where we use the notation for $A\su \Xx$ and $s>0$, 
$$B_s(A)= \{(x,m)\in \R^d\times \M \mid \|x -y\|\le s \text{ for some } (y,m')\in A\}.$$ 
\been
\im  For $Q\in \BB_0$, let
$\Ed := \big\{\Pd(Q) > 4\la(Q)\big\}$
denote the event that $\Pd(Q)$ substantially exceeds its expectation. 
 %
%
\im Let $A\su Q$ be Borel. Then, we define
$\Ep (A,s) := \big\{\Pd_A \lrsa_{\Pd_Q} \Pd_{Q\sm B_{2s}(A)} \big\}$.
\enen

We obtain the following bound on the disagreement probabilities similar to  \cite[Thm. 4.18]{hos}. The result here is only stated in the generality needed for our main theorems.

\bet\label{thm:disagree}(Disagreement probabilities for cluster-based thinning)
Let $A\su P=P_0\ti \M\su Q=Q^0\times \M \in \B_0$.
Fix $s>1$ and assume $B_{2s}(A) \su P \su Q$. Let  $\eta,\eta'\in \Nlf$ and let $\iota$ be a cube-wise ordering. Suppose that $\k, \k'$ are PIs of the form  \eqref{eq:k,k'} with a common PI $\k_0$ satisfying \eqref{cond:PIbound}, \eqref{cond:PI_sim}, and \eqref{cond:unif_sim}. 
Then, with $\bar{\eta} =  \eta \De \eta'$,
\begin{align*}
&\P\big( T^{\ms{cl}}_{Q,\es,\es}(\Pds_Q) \cap A \ne T^{\prime,\ms{cl}}_{Q,\es,\es}(\Pds_Q)\cap A \big)\le \P(\Ed ) + \P(\Ep(A,s))\\
& + 4\int_{ B_{2s}(A)} \P(\xx\lrsa_{\Pd_{\xx}\cup\bar{\eta}} \Pd_{Q\sm P}\cup \bar{\eta})\la(\d\xx)+ 2\int_{Q\sm B_{2s}(A)}\int_{B_s(A)}\one\{\xx\sim\yy\} \la(\d\xx) \la(\d\yy).
\end{align*}
\ent

\bep
See Section \ref{sec:proofs_const}.
\enp

%
%
\subsection{Thinning constructions of infinite-volume Gibbs processes}
\label{ss:emb}

The cluster-based thinning  in Section  \ref{ss:fw} was only defined in a bounded window $Q\in \BB_0$. We now extend this to a thinning algorithm on arbitrary $U\su \Xx$ by taking  limits of cluster-based thinnings in bounded domains $U_n=U\cap W_n$. 

%
%
\bepr[Convergence of cluster-based thinning]
\label{pr:xxff}
Assume that the PI satisfies \eqref{cond:PIbound}, \eqref{cond:PI_sim} and \eqref{cond:perc2} and that the ordering $\iota: \Xx \to (0,1)$ satisfies $\iota(W_{n+1}\sm W_n)\su (2^{-(n+1)},2^{-n})$. Let $A\su U \su \Xx$ with $A\in \B_0$ and $\psi \in \Nlf_{\Xx\sm U}$. Set $\xi_n=\Pd_{U\sm W_n}$. Then, on the event $\{\Pd_A  \not \lrsa_{\Pd_{U }} \xi_n\}$, it holds for all $m\geq n$,
\begin{equation}\label{e:eqonA}
	T^{\ms{cl}}_{U_n,\xi_n, \psi }(\Pds_{U_n}) \cap A = T^{\ms{cl}}_{U_{m} ,\xi_m,  \psi }(\Pds_{U_m })\cap A.
\end{equation}
In particular,  there is a limiting process on $U$,
\begin{align}
	\label{eq:xxff}
	 T^{\ff}_{U,\psi}(\Pds_U):= \bigcup_{n_0 \ge 1} \bigcap_{n \ge n_0} T^{\ms{cl}}_{U_n,\xi_n, \psi}(\Pds_{U_n})
\end{align}
such that for all  $n\ge 1$,
$$ \P\big( T^{\ff}_{U,\psi}(\Pds_U)  \cap A\ne T^{\ms{cl}}_{U_{m}, \psi}(\Pds_{U_{m} })\cap A \text{ for some $m \ge n$}\big) \leq \P(\Pd_A   \lrsa_{\Pd_U} \xi_n).$$ 
The process $T^\ff_{U,\psi}(\Pds_U)$ is the infinite-volume Gibbs process $\X(U,\psi)$ with PI $\k$.
\enpr

%
%

The theorem below shows how we can also construct the Gibbs process on $\Xx$ in two steps. The proof is essentially  the DLR-equations \cite{dobrushin,lanford}. However, since these are usually only stated when $U\in \B_0$, we state the result below.

\bepr
\label{pr:chimera}
Let $U \su \Xx$ be Borel and assume \eqref{cond:PIbound},  \eqref{cond:PI_sim}, and \eqref{cond:perc2} hold.  Let $\X$ be an infinite-volume Gibbs process on $\Xx$, let $\xi = \X \cap (\Xx \sm U)$, and define $\XX'=\xi \cup T^{\ff}_{U,\xi }({\PP}^*_U)$, where ${\PP}^*$ is independent of $\X$. Then, $\XX'$ is again an infinite-volume Gibbs process on $\Xx$.
\enpr

Proofs of Proposition \ref{pr:xxff} and \ref{pr:chimera} are given in Section \ref{sec:proofs_const}.
The above construction of the Gibbs process in unbounded domains does not allow us to control what happens if we locally change the boundary conditions $\psi$ inside a bounded set $B\su \Xx \sm U$. This requires a generalization of Proposition \ref{pr:xxff} which will be given later in Proposition \ref{prop:2_ff} in Section \ref{sss:dcff}.


\section{Disagreement coupling for Gibbs processes with pair potential}\label{sec:RCM}
The close relation between Gibbs point processes with pair potential and the RCM  indicated in Section~\ref{ss:pair} motivates the idea  that disagreement coupling could be constructed by working on the clusters of an RCM rather on the graph with edges given deterministically by the relation $\sim$. This poses some extra challenges since the edges of the RCM are random. In \cite{betsch}, this was solved by approximating the RCM by a sequence of graphs  with deterministic edges in some enlarged space. However, for the proofs of the main Theorems \ref{thm:1}, \ref{thm:2}, and \ref{thm:3}, we need to construct the Gibbs process directly as a thinning of an RCM. Therefore, we extend the disagreement coupling framework to the RCM in this section. The proofs are deferred to Section \ref{sec:proofsRCM}, but we remark that the main proofs rely on the approximation idea of \cite{betsch}, which allows us to draw on the machinery set up in  Section \ref{sec:const}.

%
%
\subsection{Standard embedding in the RCM in finite-volume domains}
\label{ss:setup} We first give a version of the standard embedding where the Gibbs process on some $Q\in \BB_0$ is constructed by thinning an RCM by visiting the vertices in the order given by some injective ordering map $\iota:\Xx \to (0,1)$. We will assume that $\sup \iota(\Xx\sm Q) \le \inf \iota(Q)$.

%
%
\subsubsection{Graph construction from marked points}
\label{sss:exact}
By its definition, the RCM requires additional randomness to decide whether edges are inserted in the model or not. In this section, the additional randomness is added to the model as random marks on the vertices. From these marked points, the RCM can be constructed by a deterministic construction rule.
More generally, we construct graphs on $Q$ from  elements $\bar\vp^*\in \Nlf_{\bar{Q}^*}$, where $\bar {Q}^* := Q^* \times [0,1]^\Z=Q \times [0,1] \times [0,1]^\Z$ is equipped with the product measure $\lambda\otimes \Q \otimes \Q^\Z$. The first coordinate of a point $(\xx,u,\rr)\in \bar\vp^*$ corresponds to a vertex of the graph,  while the third coordinate is used to define the edges. The second coordinate is only used later when we construct a thinning of the graph and is sometimes omitted from the notation. When we project the points in $\bar\vp^*$ onto one or more coordinates, we obtain the point sets  $\vp^*\in \Nlf_{Q^*}$,  $\vp \in \Nlf_{Q}$, and $\bar\vp \in \Nlf_{Q\times [0,1]^\Z}$, respectively. If two points in $\bar\vp^*$ project to the same point in $\vp$, we say that  that $\vp$ has \emph{multiple points}.

Let $\bar\vp^*\in \Nlf_{\bar Q^* }$ and let $\psi \in\Nlf$ be a fixed set of points. Assume that $\vp$ does not have multiple points. For $(\xx,u,\rr)\in \bar\vp^*$, we associate some of the coordinates in the sequence $\rr = (r_i)_{i\in \Z}$ with points in $\psi \cup \vp$. For this, let $\tilde{\iota}$ be an ordering such that  $\sup \tilde{\iota} (W_{n+1} \sm W_n  ) \le \inf \tilde{\iota}(W_n)$. This induces a total ordering on $\vp\cup \psi$.  We associate the $i$th coordinate $r_i$ of $\rr$,  with the $i$th point  after $\xx$ in the $\tilde{\iota}$-ordering of $\psi\cup \vp$. If $\yy$ is the $i$th point after $\xx$, we write $r_i=r_{\xx,\yy}$. In particular, $r_0=r_{\xx,\xx}$. Often, we write $(\xx,u,\rr)\in \bar\vp^*$ as $(\xx,u,(r_{\xx,\yy})_{\yy\in \psi\cup\vp})$. The coordinates of $\rr$ that are not matched to a point in $\vp \cup \psi$ will never be used. The ordering $\tilde{\iota}$ is fixed throughout and is only used to make this assignment. We construct a graph from an element $\bar\vp^*\in \Nlf_{\bar{Q}^*}$ as follows.

%
%
\bede[Graph construction from $\Nlf_{\bar{Q}^*}$]
\label{def:qrcm}
Let $\bar\vp^*\in \Nlf_{\bar{Q}^*}$ and $\psi\in\Nlf_{\Xx\sm Q}$. If $\vp$ does not have multiple points, we form a graph $\Gamma(\bar\vp^*,\psi)$ with vertex set $\vp\cup \psi$ and an edge between $\xx\in \vp$ and $\yy\in \vp \cup \psi$ if $\yy<_\iota \xx$ and $r_{\xx,\yy}\le \pi(\xx,\yy)$. If $\vp$ has multiple points, we set $\Gamma(\bar\vp^*,\psi)$ to be the empty graph.
\ende

Let $\psi\in\Nlf_{\Xx\sm Q}$. Let $\Phi^*_{Q,\psi}$ be a Poisson process on $\bar Q^*$, which we will think of as 
a  Poisson process $\Pds_Q$ on $Q^*$ with points independently marked by a sequence of i.i.d.\ uniform points in $[0,1]$. The $\psi$ in the notation $\Phi_{Q,\psi}^*$ indicates that for each point $(\xx,u,\rr)\in \Phi^*_{Q,\psi}$, we  think of  $\rr$ as a sequence $\rr=(r_{\xx,\yy})_{\yy\in \psi \cup \Pd_{Q}}$ as explained above. 
Applying the graph construction in Definition \ref{def:qrcm} to $\Phi^*_{Q,\psi}$, we get that $\Gamma(\Phi^*_{Q,\psi},\psi)$ is a random graph with vertex set $\Pd_{Q} \cup \psi$ and given $\Pd_Q$, there is an edge between $\xx\in \Pd_{Q}$ and $\yy\in \Pd_{Q}\cup \psi$ if $\yy<_{\iota} \xx$ and $r_{\xx,\yy}\le \pi(\xx,\yy)$. As a result, there is an edge between $\xx$ and $\yy$ with probability $\pi(\xx,\yy)$ and  all edges are independent of each other. Thus, the distribution of $\Gamma(\Phi^*_{Q,\psi},\psi)$ is that of the RCM $\Gamma(\Pd_Q,\psi)$ introduced in Definition \ref{def:RCM}. Note that there are two marks, $r_{\xx,\yy}$ and $r_{\yy,\xx}$, that could be used for deciding if there is an edge between $\xx$ and $\yy$. Here, we used $r_{\xx,\yy} $ when $\yy<_\iota \xx$. We will later see different rules for choosing between $r_{\xx,\yy}$ and $r_{\yy,\xx}$ that lead to different constructions of the RCM.

\subsubsection{Standard embedding in the RCM} \label{sss:se_RCM}
We now introduce a version of the standard Poisson embedding, where the Gibbs process is embedded in the vertex set of the  RCM $\Gamma(\Phi^*_{Q,\psi},\psi)$ in such a way that there are no edges between the Gibbs points. Given the boundary conditions $\psi\in \Nlf_{\Xx\sm Q}$, we define  a thinning map  $\Trr_{Q,\psi}: \Nlf_{\bar Q^* } \to \Nlf_{Q}$ as follows. 
Let $\bar\vp^*\in \Nlf_{\bar Q^*}$. If $\vp$ contains multiple points, we set $\Trr_{Q,\psi}(\bar\vp^*)=\es$.  Otherwise, let $\Gamma =\Gamma(\bar\vp^*,\psi)$.  Then, $\Trr_{Q,\psi}(\bar\vp^*)$ is defined by following the thinning scheme from Definition~\ref{def:standard}  on $\vp^*$, but using the retention probabilities 
\begin{equation}\label{eq:ret_RCM}
p^{\ms{RCM}}(\xx,Q,\phi,\Gamma):= p(\xx,Q,\phi) \frac{\one\{\xx\nsim_{\Gamma} \phi\}}{\prod_{\yy\in\phi}(1-\pi(\xx,\yy))}=\one\{\xx\nsim_{\Gamma} \phi\} \frac{Z_{Q_{(\xx,\ff)}}(\phi\cup \{\xx\})}{Z_{Q_{(\xx,\ff)}}(\phi)},
\end{equation}
where $\xx\in\vp$,  $\phi\su (\vp \cup \psi) \cap Q_{(-\ff,\xx)}$ and $\xx \nsim_{\Gamma} \phi$ means that there are no direct edges from $\xx$ to any point $\phi$ in $\Gamma$. If the denominator is zero, then $\xx $ is connected to a point in $ \phi$ with probability 1 and we use the convention $0/0:=0$.
Note that the construction only depends on $\bar\vp^*$ via $\Gamma$ and $\vp^*$. Each time we consider a point for thinning, it is automatically removed if it has an edge to a point previously kept. Otherwise, it is only removed with a certain probability. Thus, the subgraph of $\Gamma$ with vertex set $\Trr_{Q,\psi}(\bar\vp^*)$ will have no edges. Applied to $\Phi_{Q,\psi}^*$ and the associated RCM $\Gamma =\Gamma(\Phi_{Q,\psi}^*,\psi)$, we show in Proposition \ref{prop:pois_emb_RCM} that $\Trr_{Q,\psi}(\Phi_{Q,\psi}^*)\stackrel d=\X(Q,\psi)$.

\bepr\label{prop:pois_emb_RCM}
Suppose the PI is given in terms of a non-negative pair-potential.
Let $Q\in \BB_0$ and $\psi\in \Nlf_{\Xx\sm Q}$. 
 Then, $\Trr_{Q,\psi} (\Phi^*_{Q,\psi})\stackrel d= \X(Q,\psi)$. 
\enpr

\bep
See Section \ref{ss:proofs_se_RCM}.
\enp

\subsection{Disagreement coupling in finite-volume domains for the RCM}\label{ss:dc_RCM}

The aim of this section is to construct flexible disagreement couplings allowing us to construct the Gibbs process by cluster-wise thinning of an RCM similar to Definition \ref{ex:dc1}. We first give the construction of a single step of the algorithm in Section \ref{sss:fixed}. An important difference from Definition \ref{ex:dc1} is that the RCM is not fixed in advance. Instead, it is constructed by a step-wise exploration procedure.  Moreover, rather than exploring the underlying Poisson process in increasing parts of space as in Definition \ref{ex:dc1}, we divide it into a sequence of point processes, each of which is Poisson given the previous ones, and explore one at a time. A similar sequential  construction  was successfully used earlier in literature \cite[Proposition 2]{rcm}.  The construction of the RCM is given in Section \ref{sss:expl}. The thinning takes place during  this iterative construction of the RCM. This is explained in Section \ref{sss:dc_rcm}.  

\subsubsection{Poisson embedding from fixed points}\label{sss:fixed}
In this section, we construct a thinning map $\Trr_{Q,\nu,\psi}: \Nlf_{\bar{Q}^*} \to \Nlf_Q$, where $\nu \in \Nlf_{\Xx\sm Q}$ is a fixed set of points  (not necessarily disjoint from $\psi$) playing the role of $B$ in Section \ref{sec:const}.  Applying $\Trr_{Q,\nu,\psi}$ to the process $\Phi^{*}_{Q,\psi\cup \nu}$, we obtain a process with distribution $\X(Q,\psi)$, but the thinning now starts by considering the vertices of $\Gamma(\Phi^{*}_{Q,\psi\cup \nu},\psi\cup \nu)$ connected by a direct edge to a point in $\nu$. The thinning $\Trr_{Q,\nu,\psi}$ will be used to form a single exploration step in the exploration-based thinning in Section \ref{sss:expl}.  
 We first give another  algorithm for constructing an RCM from $\Phi^*_{Q,\psi\cup \nu}$ that is different from $\Gamma(\Phi^*_{Q,\psi\cup \nu},\psi)$ on a realization level but has the same distribution.

\def\GGG{\Gamma^{\ms{Init}}}
\def\GGB{\Gamma^{\ms{Bulk}}}
\bede[Long-range percolation model $\Gamma_\nu(\bar\vp^*,\psi)$ from a finite set $\nu$]
\label{def:rcm_nu}
Let $\bar\vp^*\in\Nlf_{\bar{Q}^*}$ and $\psi,\nu\in \Nlf_{\Xx\sm Q}$ (not necessarily disjoint). Then, we let
\beit 
\im $V:=V_\nu(\bar\vp^*,\psi)=\{\xx\in \vp \mid \xx \sim_{\Gamma(\bar\vp^*,\psi\cup \nu)} \nu\}$  be the vertices of $\Gamma(\bar\vp^*,\psi\cup \nu)$ having a direct edge to a point in $\nu$; 
\im  $\GGG_\nu(\bar\vp^*,\psi)$ denote the subgraph of $\Gamma(\bar\vp^*,\psi\cup \nu)$ with vertices $V\cup \psi $ and all induced edges; 
\im $\GGB_\nu\big(\bar\vp^*,V\cup(\psi\sm \nu)\big)$ be the graph with vertex set $\vp \cup (\psi\sm \nu)$ where there is an edge between  $\xx\in \vp\sm V$ and $y\in V\cup (\psi\sm \nu)$ if $r_{\xx,\yy}\le \pi(\xx,\yy)$ and an edge between $\xx,\yy\in \vp \sm V$ if $\xx >_{\iota} \yy$ and $r_{\xx,\yy} \le \pi(\xx,\yy)$;
\im $\Gamma_\nu(\bar\vp^*,\psi) := \GGG_\nu\big(\bar\vp^*,\psi\big) \cup \GGB_\nu\big(\bar\vp^*,V\cup(\psi\sm \nu)\big),$
where the union of graphs is given by taking union of vertices and edges.
\enit
\ende

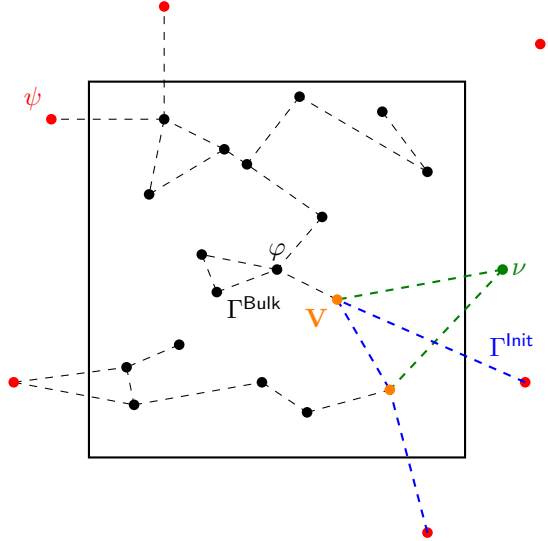
\begin{figure}[!htpb]
\begin{tikzpicture}
    \draw[thick] (0,0) rectangle (5,5);

    \coordinate (A) at (0.5,1.2);
    \coordinate (B) at (1.0,4.5);
    \coordinate (C) at (1.7,2.2);
    \coordinate (D) at (2.3,1.0);
    \coordinate (E) at (2.8,4.8);
    \coordinate (F) at (3.1,3.2);
    \coordinate (G) at (0.8,3.5);
    \coordinate (H) at (1.5,2.7);
    \coordinate (I) at (2.1,3.9);
    \coordinate (J) at (2.9,0.6);
    \coordinate (K) at (1.2,1.5);
    \coordinate (L) at (0.6,0.7);
    \coordinate (M) at (1.8,4.1);
    \coordinate (N) at (3.9,4.6);
    \coordinate (O) at (4.5,3.8);
    \coordinate (P) at (2.5,2.5);

    \coordinate (Q) at (3.3,2.1);
    \coordinate (R) at (4.0,0.9);

    \foreach \pt in {A,B,C,D,E,F,G,H,I,J,K,L,M,N,O,P} {
        \fill[black] (\pt) circle (2pt);
    }

\draw[dashed] (-1,1) -- (0.6,0.7);
\draw[dashed] (-1,1) -- (0.5,1.2);
\draw[dashed] (1,4.5) -- (1,6);
\draw[dashed] (1,4.5) -- (-0.5,4.5);

    \node at (P) [above] {$\varphi$};

    \foreach \x/\y in {
        -1/1, -0.5/4.5, 1/6, 4.5/-1, 6/5.5, 5.8/1
    } {
        \fill[red] (\x,\y) circle (2pt);
    }

    \definecolor{darkgreen}{rgb}{0.0,0.5,0.0}
    \fill[darkgreen] (5.5,2.5) circle (2pt);
    \node[darkgreen] at (5.5,2.5) [right] {$\nu$};

    \node[red] at (-0.5,4.5) [above left] {$\psi$};

    \draw[darkgreen, thick, dashed] (5.5,2.5) -- (Q);
    \draw[darkgreen, thick, dashed] (5.5,2.5) -- (R);

\draw[blue, thick, dashed] (5.8,1) -- (Q);
\draw[blue, thick, dashed] (4.5,-1) -- (R);
\draw[blue, thick, dashed] (Q) -- (R);

    \draw[dashed] (A) -- (L);
    \draw[dashed] (B) -- (G);
    \draw[dashed] (C) -- (P) node[xshift=-0.3cm, yshift=-0.5cm]{$\GGB$} ;
    \draw[dashed] (D) -- (J);
    \draw[dashed] (E) -- (I);
    \draw[dashed] (F) -- (P);
    \draw[dashed] (G) -- (M);
    \draw[dashed] (H) -- (C);
    \draw[dashed] (I) -- (M);
    \draw[dashed] (J) -- (R); 
    \draw[dashed] (K) -- (A);
    \draw[dashed] (L) -- (D);
    \draw[dashed] (M) -- (B);
    \draw[dashed] (N) -- (O);
    \draw[dashed] (O) -- (E);
    \draw[dashed] (P) -- (H);
    \draw[dashed] (Q) -- (P); 
    \draw[dashed] (F) -- (I);

    \coordinate[label=0:{\textcolor{blue}{{$\GGG$}}}] (V) at (5.2,1.5);

\fill[orange] (Q) circle (2pt);
\fill[orange] (R) circle (2pt);
\node[orange] at (Q) [below left] {$\textbf{V}$};

\end{tikzpicture}
\caption{Illustration of the construction of $\Gamma_\nu(\bar\vp^*,\psi)$. The green point is $\nu$, and the red points are $\psi$. The orange points are $V$, which consists of the points connected to $\nu$. $\GGG$ is shown in blue as the subgraph with vertex set $V\cup \psi$ and $\GGB$ is the graph shown in black. }
	\label{fig:rcm_nu}
\end{figure}

The construction of the long-range percolation model  $\Gamma_\nu$ is illustrated in Figure \ref{fig:rcm_nu}. Note that the green edges are not part of $\Gamma_\nu$ since $\nu$ is disjoint from $\psi$ in this example and hence $\nu$ is not part of the vertex set of $\Gamma_{\nu}$. The green edges are only used for defining $V$. Note that $\Gamma_\nu(\bar\vp^*,\psi)$ and $\Gamma(\bar\vp^*,\psi)$ have the same edges between points in $V\cup \psi$ and between points in $(\vp \sm V)\cup \psi$, but the edges connecting $V$ to $\vp \sm V$ (edges between orange and black points in Figure \ref{fig:rcm_nu}) may be different because different marks may have been used to decide if edges are present. The next lemma shows that applying $\Gamma_\nu$ to $\Phi_{Q,\psi\cup \nu}^*$, we obtain an RCM.

\bel[The graph $\Gamma_\nu(\Phi^*_{Q,\psi\cup\nu},\psi)$ is an RCM] \label{lem:RCM_it}
The point processes $V=V_\nu(\Phi^*_{Q,\psi})$ and  $\Pd_Q\sm V$ are independent Poisson processes on $Q$ with intensity measures $(1-\rho_\nu)\la$ and  $\rho_\nu\la$, respectively, where 
\begin{align}
	\label{eq:rn}
\rho_\nu(\xx):=\prod_{\yy\in \nu}(1-\pi(\xx,\yy)).
\end{align}
Moreover, $\Gamma_\nu(\Phi^*_{Q,\psi\cup\nu},\psi)\stackrel d= \Gamma(Q,\psi)$ and conditionally on $\GGG_\nu(\Phi^*_{Q,\psi\cup\nu},\psi)$, we have $\GGB(\Phi^*_{Q,\psi\cup\nu} ,V\cup (\psi\sm \nu))$ is an RCM $\Gamma^{\rho_\nu\la}(Q,V\cup (\psi\sm \nu) )$ on $Q$ with underlying Poisson intensity $\rho_\nu\lambda$.
\enl

\bep
See Section \ref{ss:proofs_dc_RCM}.
\enp

Next, we give a different way of constructing the Gibbs process $\XX(Q,\psi)$ as a thinning of $\Gamma_\nu(\Phi^*_{Q,\psi\cup \nu},\psi\cup \nu)$.  For this we define a thinning map $\Trr_{Q,\nu,\psi}: \Nlf_{\bar{Q}^*} \to \Nlf_Q$. For $\bar\vp^*\in \Nlf_{\bar{Q}^*} $, we let $\Trr_{Q,\nu,\psi}(\bar\vp^*)=\es $ if $\vp$ has multiple points. Otherwise, let $\Gamma=\Gamma(\bar\vp^*,\psi\cup \nu)$ and $V=V_\nu(\bar\vp^*) $. 
The thinning proceeds by first thinning $V$ and given the result $\xi\su V$, construct a thinning of $\vp \sm V$.  

\begin{itemize}
\item For the thinning of $V$, we visit the points of $V$ in the order of $\iota$. Suppose that $\xx_1,\dots,\xx_k\in V\cap Q_{(-\ff,\xx)}$ are the points so far kept in the thinning when considering the thinning of $ (\xx,u)\in V^*$. Let $\phi=\psi\cup \{\xx_1,\ldots,\xx_k\}$. 
Let $\Gamma^{\rho_\nu^{\xx}\lambda}(Q,\psi\cup \nu)$ be an RCM with vertices on a Poisson process $\PP_{Q}^{\rho_\nu^{\xx}\lambda}$ of intensity $\rho_\nu^{\xx}\lambda$ where 
$$\rho_\nu^{\xx}(\yy) =  \rho_\nu(\yy)\one\{\yy\in Q_{(-\ff,\xx )}\} + \one\{\yy\in Q_{(\xx,\ff)}\}.$$
Then, $\xx$ is kept if $u\le p(\xx,Q,\nu,\phi,\Gamma)$ where
\begin{equation}\label{eq:p_nu}
 p(\xx,Q,\nu,\phi,\Gamma)=\one\{\xx\nsim_{\Gamma} \phi \} \frac{Z_Q(\phi\cup\{\xx\},\nu,\xx)}{Z_Q(\phi,\nu,\xx)}
 \end{equation}
with 
\begin{align}\nonumber
Z_Q(\phi,\nu,\xx)={\E}\bigg( \prod_{\yy\in \PP_{Q_{(-\ff,\xx)}}^{\rho_\nu^{\xx}\lambda}} \one\{ \yy \nsim_{\Gamma^{\rho_\nu^{\xx}\lambda}(Q,\phi\cup \nu)} \phi\sm \nu\}\prod_{\yy\in \PP_{Q_{(\xx,\ff)}}^{\rho_\nu^{\xx}\lambda}}\one\{ \yy \nsim_{\Gamma^{\rho_\nu^{\xx}\lambda}(Q,\phi\cup \nu)}\phi\}\\
 \times\prod_{\yy_1,\yy_2\in \PP_{Q}^{\rho_\nu^{\xx}\lambda }} \one\{\yy_1 \nsim_{\Gamma^{\rho_\nu^{\xx}\lambda}(Q,\phi\cup \nu)} \yy_2\} \bigg).\label{eq:Z_bar_def}
 \end{align}
\item Given the result $\xi \su V$ of thinning $V$, we follow the Poisson embedding procedure on $\vp^*\sm V^*$ with the retention probabilities \eqref{eq:ret_RCM} but where $\Gamma$ is now the graph $\GGB(\bar\vp^*,V \cup (\psi\sm\nu))$, the boundary conditions are $\xi \cup (\psi\sm \nu)$,  the reference measure is $\rho_\nu\lambda$ when computing partition functions, and the ordering is 
$$\iota_{\text{new}}(\xx) = \tfrac{1}{2}\iota(\xx)\one\{\xx \notin Q\sm V\} + \big(\tfrac{1}{2} + \tfrac{1}{2}\iota(\xx))\one\{\xx \in Q\sm V\}.$$ 
The latter is a technical detail, since we always assume the boundary points to come first in the ordering.
Denote the thinning of $\vp\sm V$ by $\xi'$.
\item The resulting thinned process is denoted $\Trr_{Q,\nu,\psi}(\bar\vp^*)=\xi \cup \xi'$.
\end{itemize}

Applying $\Trr_{Q,\nu,\psi}$ to $\Phi_{Q,\psi\cup \nu}^*$, we obtain a Gibbs process. 

%
%
\bel[Correctness of thinning construction starting from $\nu$]
\label{lem:TRCM_V}
Let $\k$ be a PI with non-negative pair potential.
Let $Q\in \BB_0$ and $\nu,\psi \in \Nlf_{\Xx\sm Q}$. Then,
 $\Trr_{Q,\nu,\psi}(\Phi^*_{Q,\psi \cup\nu})\stackrel d=\X(Q,\psi)$. Conditionally on $\xi =\Trr_{Q,\nu,\psi}(\Phi^*_{Q,\psi \cup\nu})\cap V$, the distribution of $\xi'=\Trr_{Q,\nu,\psi}(\Phi^*_{Q,\psi\cup \nu})\sm \xi$ is $  \X^{\rho_\nu\la}(Q,\xi \cup (\psi\sm \nu))$.
\enl

\bep
The proof can be found in Section \ref{ss:proofs_dc_RCM}.
\enp

\subsubsection{An exploration-based construction of the RCM}
\label{sss:expl}
We next give an exploration-based construction of the RCM by iterative use of  Definition \ref{def:rcm_nu}. 
Let $\bar\vp^*\in \Nlf_{\bar Q^*}$ and $\nu,\psi\in \Nlf_{\Xx\sm Q}$.
Then, we construct the graph $\GGd(\bar\vp^*,\psi)$ by a step-wise procedure. 
At the beginning of step $k+1$, we have 
\begin{itemize}
\im a set $S_{k}\su \vp \cup \psi$ of so far explored points,
\im a graph $\GGd_k(\bar\vp^*,\psi)$ with vertex set $S_k$.
\end{itemize}
During the step $k+1$, we further introduce a reference point $\xx_{k + 1}$ which will be the starting point of the cluster explored in step $k+1$.
The details are given in Definition \ref{def:GGd} below. We write $S_k^*$ and $\bar S_k^*$ when the marks from $\bar\vp^*$ are kept and write $(\GGd_k)^*$ for the graph corresponding to $\GGd_k$ but having vertex set $S_k^* \cup \psi$. For $\bar\vp^*\in \Nlf_{\bar Q^*}$, we use a slight abuse of  notation writing $ \bar\vp^*\sm \bar A^* $ where for $(\xx,u,\rr)\in \bar\vp^*\sm \bar A^* $, the sequence $\rr$ is still identified with $(r_{\xx,\yy})_{\yy\in \vp\cup\psi\cup\nu}$ in the same way as the corresponding point in $\bar\vp^*$.  

\bede[Exploration based construction of $\GGd$]
\label{def:GGd}
Let $\bar\vp^*\in \Nlf_{\bar Q^*}$ and $\nu,\psi\in \Nlf_{\Xx\sm Q}$.
If $\vp $ has multiple points, we set $\GGd=\es$. Otherwise, initialize the construction by setting  $\nu_0=\es$, $S_0=V_0=\psi$, $\xx_0=-\ff$, $\xi_0=\es$, $\iota_0=\iota$,  and  $\GGd_{0}$ the graph with vertex set $\psi$ and no edges.  Suppose step $k$ is completed and $\vp\sm S_k\ne \es$.  
Then, we first determine a set $V_{k + 1}$ of points to be explored in the next step as follows.
\been
\item Define the window of unexplored points $Q_{k+1}=Q_{(\xx_k,\ff)}\sm S_k.$
\item 
Define the ordering
$$\iota_{k+1}(\xx) = \tfrac12\iota_k(\xx) \one\{\xx\in \Xx \sm Q_{k+1} \} + (\tfrac12  + \tfrac12\iota_k(\xx))\one\{\xx\in Q_{k+1}\}$$
that reorders the explored points  to come before $Q_{k+1}$.
\item Determine a subset $\nu_{k+1}$ of $(\nu\cup S_k)\sm (\cup_{l\le k}\nu_l)$ by some deterministic rule only depending on $(\GGd_k)^*$ and $\cup_{l\le k}\nu_l$.
\item 
Construct the graph $\Gamma_{\nu_{k+1}}(\bar\vp^*\sm \bar  S_k^*,S_k\sm \cup_{l\le k} \nu_l) $  as in Definition \ref{def:rcm_nu} using the ordering $\iota_{k+1}$.
\item Let $V_{k+1} = V_{\nu_{k+1}}(\bar\vp^*\sm \bar S_k^*,S_k\sm \cup_{l\le k} \nu_l)$ be the set of points in $\vp\sm S_k$   connected to  $\nu_{k+1}$ in $\Gamma_{\nu_{k+1}}(\bar\vp^*\sm \bar S_k^*,S_k\sm \cup_{l\le k} \nu_l)$. Those are the so far unexplored points connected to $\nu_{k+1}$. 
	\enen
	
To determine the next step, we distinguish two cases depending on whether or not the exploration step was successful. That is, whether $V_{k + 1} \ne \es$ or $V_{k + 1}= \es$.

{\bf Case 1 -- Successful exploration $\bs{V_{k + 1} \ne \es}$.}

\begin{figure}
	\begin{center}
\begin{tikzpicture}
	\def\xmin{0}
	\def\xmax{8}
	\def\ymin{0}
	\def\ymax{2}
	
	\draw[thick] (\xmin,\ymin) rectangle (\xmax,\ymax);
	
	\coordinate (pt1) at (1.2, 1.1);
	\coordinate (pt3) at (3.3, 1.7);
	\coordinate (pt4) at (4.8, 1.2);   
	\coordinate (pt5) at (4.1, 0.3);   
	\coordinate (pt6) at (5.2, 1.5);
	\coordinate (pt7) at (5.7, 0.5);
	\coordinate (pt8) at (6.8, 1.0);
	\coordinate (pt9) at (7.1, 0.7);
	\coordinate (pt10) at (7.5, 1.3);
	
	\foreach \i in {1,3,4,5} {
		\fill[red] (pt\i) circle (1.5pt);
	}
	\foreach \i in {6,...,10} {
		\fill[black] (pt\i) circle (1.5pt);
	}
	
	\node[red] at (1.2,0.8) {$\xx_{k+1}$};
	
	\node[red] at (1,1.8) {$S_k$};
	
	\draw[dashed] (pt4) -- (pt6);
	\fill[blue!50] (pt6) circle (1.5pt);
	
	\draw[dashed] (pt4) -- (pt9);
	\fill[blue!50] (pt9) circle (1.5pt);
	
	\draw[dashed] (pt5) -- (pt7);
	\fill[blue!50] (pt7) circle (1.5pt);
	
	\draw[dashed] (pt5) -- (pt8);
	\fill[blue!50] (pt8) circle (1.5pt);
	
	\node[blue!50] at (6.5,1.6) {$V_{k + 1}$};
	
	\draw[red, dashed] (pt1) -- (pt3);
	\draw[red, dashed] (pt1) -- (pt4);
	\draw[red, dashed] (pt5) -- (pt1);
	
\end{tikzpicture}
\begin{tikzpicture}
    \def\xmin{0}
    \def\xmax{8}
    \def\ymin{0}
    \def\ymax{2}

    \draw[thick] (\xmin,\ymin) rectangle (\xmax,\ymax);

    \coordinate (pt1) at (1.2, 1.1);   
    \coordinate (pt3) at (3.3, 1.7);   
    \coordinate (pt4) at (4.8, 1.2);   
    \coordinate (pt5) at (4.1, 0.3);   

    \coordinate (pt6) at (5.2, 1.5);   
    \coordinate (pt7) at (5.7, 0.5);   
    \coordinate (pt8) at (6.8, 1.0);   
    \coordinate (pt9) at (7.1, 0.7);   
    \coordinate (pt10) at (7.5, 1.3);  

    \foreach \i in {1,3,4,5} {
        \fill[red] (pt\i) circle (1.5pt);
    }

    \foreach \i in {7,...,10} {
        \fill[black] (pt\i) circle (1.5pt);
    }

    \fill[blue] (pt6) circle (1.5pt);
	\node[blue] at ($(pt6) + (0.3, 0.3)$) {$\xx_{k + 1}$};

    \draw[red, dashed] (pt1) -- (pt3);
    \draw[red, dashed] (pt1) -- (pt4);
    \draw[red, dashed] (pt5) -- (pt1);

    \node[red] at (1,1.8) {$S_k$};

\end{tikzpicture}
\end{center}
\caption{Left: Successful exploration $V_{k+1}\ne  \es$. Right: Failed exploration $V_{k+1} = \es$. }\label{fig:explore}
\end{figure}
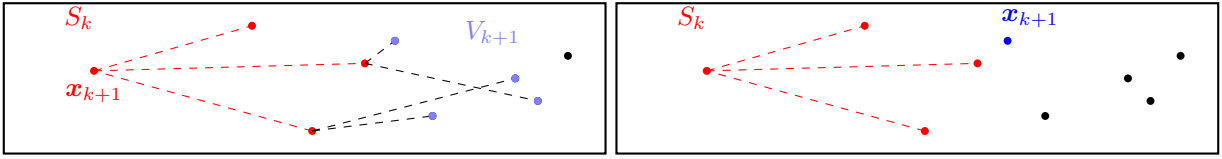

In this case, we keep exploring the same cluster, see Figure \ref{fig:explore} left.
	\been
\item Let   $S_{k+1}=S_k\cup V_{k+1}$ be the points explored after the $(k+1)$st step and keep the reference point, i.e., $\xx_{k + 1} = \xx_k$.
\item 
Define 
$\GGd_{k+1}=\GGd_k \cup \GGG_{\nu_{k+1}}(\bar\vp^*\sm \bar S_k^*,S_k\sm \cup_{l\le k} \nu_l)$ 
\item If $S_{k+1}=\vp\cup\psi$, the algorithm terminates with $\GGd = \cup_{l=0}^{k+1} \GGd_l$.
\enen

{\bf Case 2 -- Failed exploration $\bs{V_{k + 1} = \es}$.}

In this case, we start exploring a new cluster, see Figure \ref{fig:explore} right.

	\been
\item Let $\xx_{k + 1}=\inf_{\iota_{k+1}}(\vp\sm S_k)$ be the first point in the new cluster and $S_{k+1}=S_k\cup \{\xx_{k + 1}\}$ the points explored after step $k+1$.
\item 
Define $\GGd_{k+1}=\GGd_k \cup \GG_{\{\xx_{k + 1}\}}$ where $\GG_{\{\xx_{k + 1}\}}$ consists of the vertices $\{\xx_{k+1}\}\cup \psi$ and all edges in $\G_{\nu_{k+1}}(\bar\vp^*\sm \bar S_k^*,S_k\sm \cup_{l\le k} \nu_l)$ between $\xx_{k+1}$ and $\psi\sm \cup_{l\le k+1} \nu_l$.

\item If $S_{k+1}=\vp\cup \psi$, the algorithm terminates with $\GGd = \cup_{l=0}^{k+1} \GGd_l$.
\enen
In both cases, we let $V_{k+1}'=S_{k+1}\sm S_k$ be the set of points that are explored in step $k+1$.
\ende

Note that the graph $\GGd$ is not given from start, but is being built while exploring it.  The definition can be seen as a repeated use of the graph construction from fixed points in Definition \ref{def:rcm_nu} where in each step we only keep the $\GGG$ part of the graph. The only free choice in the algorithm is how we choose the sets $\nu_k$.

The main property of this construction will be that $\GGd(\Phi^*_{Q,\psi\cup \nu})$ is an RCM $\Gamma(Q,\psi)$. 

\bepr[Distributional correctness of the exploration based graph construction]
\label{prop:dc_RCM} Let $\psi,\nu\in \Nlf_{\Xx\sm Q}$.
Then, $\GGd(\Phi^*_{Q,\psi\cup \nu},\psi)\stackrel d= \Gamma(Q,\psi)$.
\enpr

\bep
See Section \ref{ss:proofs_dc_RCM}.
\enp

The choice for $\nu_k$ that we will mainly work with in the following is the one given below.

\bede[Cluster-based graph construction $\GGc_\nu$]\label{def:cl_RCM}
The  \emph{cluster-based graph construction} (denoted $\GGc_\nu$) is the graph resulting from Definition \ref{def:GGd} if we set $\nu_1 = \nu$ and $\nu_{k+1}=V_k'$ for $k\ge 1$.
\ende

The cluster-based graph $\GGc_\nu$ is built by first building the clusters connected to  $\nu$. In the first step, we explore all points directly connected to $\nu$, in the second step all points connected by two (but not one) edges and so on. The first time we encounter $V_k=\es$, these clusters are completed, and we start a new cluster from $\xx_{k+1}$. Each time $V_l =\es$, a cluster is completed and a new one is started from $\xx_{l+1}$. 
The point $\xx_l$ will be the $\iota$-smallest point in the cluster, and the sequence of starting points $\xx_l$ is increasing with respect to the $\iota$-ordering, meaning that the clusters are explored in the order of their $\iota$-smallest point. This is very similar to the cluster-based exploration from Definition \ref{ex:dc1}, where neighboring points were determined by the relation $\sim$. The difference is that the edges are not there from the beginning, but are being built from the marks $r_{\xx,\yy}$.

\subsubsection{Disagreement coupling for the RCM}\label{sss:dc_rcm}

Next, we define a version of disagreement coupling for the RCM. It assumes an exploration based graph construction. In the $k$th  exploration step, the thinning is applied to the explored set $V_k'$. The precise definition is given below.

%
%
\bede[Disagreement coupling for the RCM]
\label{def:trd_RCM}
Assume an exploration based graph construction $\GGd(\bar\vp^*,\psi)$ as in Definition \ref{def:GGd}.
If $\vp $ has multiple points, we set $\Trd_{Q,\nu,\psi}(\bar\vp^*)=\es$. Otherwise,  we construct thinnings $\xi_k\su V_k'$ step-wise as follows. Let $B_k$ denote a set of boundary conditions for step $k$ and initialize  $B_0=\psi$ and $\xi_0=\es$. Suppose step $k$ is completed and we have constructed $\xi_1,\ldots,\xi_k$. 

\been
\item Set $B_{k+1}= (B_k\sm \nu_k) \cup \xi_k$.
\item Set $\xi_{k+1}= \Trr_{Q_{k+1},\nu_{k+1},B_{k+1},\iota_{k+1}}(\bar\vp^*_{Q_{k+1}})\cap V_{k+1}'$.
\item If $S_{k+1}=\vp\cup\psi$, the algorithm terminates with $\Trd_{Q,\nu,\psi}(\bar\vp^*) =\cup_{l=1}^{k+1} \xi_l$.
\enen
\ende

The next theorem shows that applying this construction to $\Phi^*_{Q,\psi\cup\nu}$ yields a Gibbs process. 

\bet[Distributional correctness of disagreement coupling for the RCM]\label{thm:dc_RCM} Let $\k$ be a PI with non-negative pair potential.
Let $Q\in \BB_0$ and $\psi,\nu\in \Nlf_{\Xx\sm Q}$. Then,
 $\Trd_{Q,\nu,\psi}(\Phi^*_{Q,\psi\cup\nu})\stackrel d= \X(Q,\psi)$.
\ent

\bep
The proof can be found in Section \ref{ss:proofs_dc_RCM}.
\enp

Again, we will work with the following version.

\bede[Cluster-based thinning for the RCM]\label{def:cl_T_RCM}
In the case where the exploration based graph construction is the cluster-based construction $\GGc_\nu$ from Definition \ref{def:cl_RCM}, we call the corresponding thinning algorithm the \emph{cluster-based thinning of the RCM} and denote it by $\Trc_{Q,\nu,\psi}$.
\ende 

The advantage of the cluster-based thinning is that the part of the boundary conditions $\psi\cap \nu$  will only affect the clusters connected to $\nu$ allowing us to control the effect of changing boundary conditions inside $\nu$.

\bec
Suppose $\k$ has a non-negative pair potential.
Let $Q\in \BB_0 $ and $A\su Q$. Let $\psi,\psi'\in \Nlf_{\Xx\sm Q}$.  Then,
$$\P\Big(\Trc_{Q,\psi\Delta \psi',\psi} (\Phi^*_{Q,\psi\cup\psi'})\cap A \ne \Trc_{Q,\psi\Delta \psi',\psi'} (\Phi^*_{Q,\psi\cup\psi'})\cap A\Big)\le \P( \psi\Delta \psi' \lrsa_{\Gamma(Q,\psi\Delta \psi')} \Pd_A).$$
\enc

\subsection{Local perturbation of PIs}\label{ss:pert_RCM}

We now consider a slight modification of our models with pair potential, namely processes with PI of the form
\begin{equation}\label{eq:ktilde}
	\bar{\k}(\xx,\vp) = \bar\a(\xx) \exp\bigg(-\sum_{\yy\in \vp} \bar v(\xx,\yy) \bigg)
\end{equation}
where $\bar v$ is a pair potential and $\bar\a : Q \to [0,1]$ is an \emph{activity function}. The only difference from the models with pair potential considered so far is the activity function. 
We could handle the activity function by changing the reference measure to $\bar\a \la$, see Remark \ref{rem:activity}. However, we later consider a family of Gibbs models of the form \eqref{eq:ktilde} and want to have $\la$ as a common  reference measure. We will assume that there is a dominating PI $v$ such that $\bar v\le v$ for all the processes. We  construct them by thinning of a common RCM $\G^\la(Q,\psi\cup\nu,\pi)$, where the connection probabilities are $\pi = 1-e^{-v}$. We denote this thinning by $\bar T^{\ms{RCM,cl}}_{Q,\nu,\psi}$. It explores the clusters of $\GGc(\Phi_{Q,\nu\cup\psi}^*,\nu\cup\psi,\pi)$ in the same order as $ T^{\ms{RCM,cl}}_{Q,\nu,\psi}$. However, the retention probabilities will depend on $\bar \a$ and $\tilde v$.  We refer to Section \ref{ss:proofs_perturb} for the precise definition of $\bar T^{\ms{RCM,cl}}_{Q,\nu,\psi}$ and the proofs of  Theorem \ref{thm:Talpha} and  Lemma \ref{lem:pa1}  below.

In the following, we assume that $\Xx=\R^d \times \M$ is marked Euclidean space and $\k$ is a PI with pair potential $v$. Consider two PIs $\bar{\k}$ and $\bar{\k}'$ of the form \eqref{eq:ktilde}. We assume that there is a $P\su Q$ such that for all $\xx\in P$ and $\vp\in \Nlf$, we  have
\begin{align}\nonumber
	&\bar\k(\xx,\vp) =\a(\xx)\k(\xx,(\vp \cap P_1) \cup \eta )=\a(\xx) e^{-\sum_{\yy\in \eta}v(\xx,\yy)} e^{-\sum_{\yy\in \vp} v(\xx,\yy)\one\{\yy\in P_1\} }\\
	&\bar{\k}'(\xx,\vp) =\a(\xx)\k(\xx,(\vp \cap P_1') \cup \eta')=\a(\xx)e^{-\sum_{\yy\in \eta'}v(\xx,\yy)} e^{-\sum_{\yy\in \vp} v(\xx,\yy)\one\{\yy\in P_1'\} }\label{eq:tildekk'}.
\end{align}
We require $ P\su P_1\cap P_1'$, $P_1,P_1'\su Q$, and $\eta,\eta'\in \Nlf$ with $ \eta\De \eta'\in \Nlf_{\Xx\sm P}$ and $\a: P\to [0,1]$. Note that on $P$, $\bar{\k}$ has the form \eqref{eq:ktilde} with $\bar{\a}(\xx) = \a(\xx)e^{-\sum_{\yy\in \eta}v(\xx,\yy)} $ and $\bar v(\xx,\yy)=v(\xx,\yy)\one\{\yy\in P_1\}$ and similarly for $\bar{\k}'$. 
\def\EpRCM{E_{\ms{perc}}^{\mathsf{RCM}}}

The theorem below bounds the disagreement between the two processes in \eqref{eq:tildekk'} on $A\su P$. It uses the notation $\EpRCM(A,s)=\{\Pd_A\lrsa_{\G(Q,\es)} \Pd_{Q\sm B_{s}(A)}\}$ for  $s>0$.

\bet\label{thm:Talpha}(Disagreement for perturbed PIs with pair potential)
Let $\Xx=\R^d \times \M$ and assume $\iota$ is a cube-wise ordering as in Section \ref{ss:pert}.
Let $Q=Q_0\times \M \in \BB_0$ and $A=A_0\ti \M \su P=P_0\ti \M \su Q$. Let $\k$ be a PI with pair potential satisfying \eqref{cond:int_v}, \eqref{cond:local_int_v}, and  \eqref{cond:perc_RCM}. Let $\bar T^{\ms{RCM,cl}}_{Q,\es,\es}$ and $\bar T^{\prime,\ms{RCM,cl}}_{Q,\es,\es}$ denote thinnings  with respect to PIs $\bar{\k}$ and $\bar{\k}'$ of the form \eqref{eq:ktilde} that satisfy \eqref{eq:tildekk'} on $P$.  Let $s>1$ and assume  $B_{2s}(A) \su P $. Then,
\begin{align*}
	&\P(\bar T^{\ms{cl}}_{Q,\es,\es}(\Phi_{Q,\es}^{*})\cap A \ne \bar T^{\prime,\ms{cl}}_{Q,\es,\es}(\Phi_{Q,\es}^{*})\cap A) \le \P(\Ed)+\P(\EpRCM(A,s)) \\
	&+16\la(Q) \int_{B_{2s}(A)}\P(\xx\lrsa_{\G_\xx (Q,\bar{\eta})} \Pd_{Q\sm P} \cup \bar{\eta} ) \la(\d \xx) + 8\la(Q)\int_{Q \sm B_{2s}(A) }\int_{B_{s}(A)} \pi(\xx,\yy) \la(\d \xx) \la(\d \yy),
\end{align*}
where $\bar{\eta}=\eta \De \eta'$.
\ent

\begin{figure}[!htpb]
\centering
\begin{tikzpicture}[scale=0.9, every node/.style={font=\small}]

\begin{scope}

\draw[line width=1pt,fill,gray!30] (-2,0) rectangle (2,4);
\draw[line width=1pt] (-2,0) rectangle (2,4);

\fill (0,2) circle (1.2pt);
\node[anchor=west] at (0.1,2) {$x$};

\draw[fill,white] (0,2) circle (0.6);
\draw[dashed] (0,2) circle (0.6);
\draw[dashed] (0,2) circle (1.0);

\node at (-0.8, 0.5) {$ P_{1,-}$};

\end{scope}

\begin{scope}[xshift=6cm]

\draw[line width=1pt] (-2,0) rectangle (6,4);


\draw[dashed,fill,gray!20] (0,2) circle (1.6);
\draw[dashed,fill,white] (0,2) circle (0.6);
\draw[dashed] (0,2) circle (0.6);
\draw[dashed] (0,2) circle (1.1);
\draw[dashed] (0,2) circle (1.6);

\fill (0,2) circle (1.2pt);
\node[anchor=west] at (0,2) {$x$};

\node at (0,3.3) {$P_{1,-}$};
\node at (0,2.8) {$A$};
\node at (0,2.3) {$P_3$};

\fill (4,2) circle (1.2pt);
\node[anchor=west] at (4,2) {$y$};

\draw[dashed,fill,gray!20] (4,2) circle (1.6);
\draw[dashed,fill,white] (4,2) circle (0.6);
\draw[dashed] (4,2) circle (0.6);
\draw[dashed] (4,2) circle (1.6);

\fill (4,2) circle (1.2pt);
\node[anchor=west] at (4,2) {$y$};

\node at (4,3) {$P_{2,-}$};

\node at (2,0.5) {$P_1 = Q$};

\end{scope}

\end{tikzpicture}
\caption{
Left:  $P_1=Q$ is the whole square while $P_{1,-}=Q\setminus B_s(x)$ is the grey region, and $P_3=P_{2,-}=\varnothing$. The PI $\k'$ lives on all of $P_1=Q$, while $\k''$ vanishes on the ball and in the grey region, it ignores boundary points in the grey ball. Theorem \ref{thm:Talpha} bounds the probability that the two processes disagree outside the larger ball.
Right: $\k'$ vanishes on the smallest ball $P_3$ around $\xx$, while $\k''$ vanishes outside the two grey annuli. Theorem \ref{thm:Talpha} bounds the probability that the two processes disagree on the medium sized ball $A$ around $\xx$.
}
\label{fig:pert_RCM}
\end{figure}
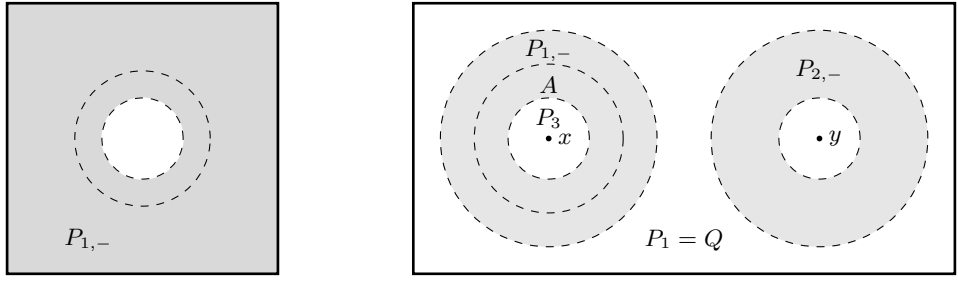

We will apply the theorem above to  PIs of the form given in the following lemma.  In Section \ref{sec:pac}, we will apply it in the concrete case where $Q$ is a window and $P_1=Q$, $P_{1,-} =Q\sm B_s(\xx)$, $P_3=P_{2,-}=\es$ and $A=Q\sm B_{4s}(\xx)$. This is illustrated in the left panel of Figure \ref{fig:pert_RCM}.  In Section \ref{sec:malliavin_stein}, we consider the case $P_1=Q$, $P_{1,-}=B_{r+3s}(\xx)$, $P_3= B_{s}(\xx)$, $P_{2,-}=B_{r + 3s}(\yy)\sm B_s(\yy)$,   $A=B_r(\xx)$ and $P=B_{r+3s}(\xx)$, see right panel of Figure \ref{fig:pert_RCM}.

\bel[Perturbed PIs]
\label{lem:pa1} 
Let $\k$ be a PI induced by a pair potential $v$ satisfying \eqref{cond:int_v}, \eqref{cond:local_int_v}, \eqref{cond:perc_RCM}. Let $Q=Q_0\ti \M$, $P_i =(P_i)_0\ti \M $ for $i=1,2,3$ and $P_{i, -}=(P_{i, -})_0\times \M$ for $i=1,2$ be such that $P_{i, -} \su P_i \su Q $ for $i= 1,2$,   $P_{1, -} \cap P_{2, -} =\es$ and  $P_3\su P_{1,-}$. Moreover, let $\psi\in \Nlf_{\Xx \sm P_1}, \psi_i \in \Nlf_{\Xx\sm P_{i,-}} $, $i\in \{ 1, 2\}$. Consider PIs $\k', \k''$ defined as 
\begin{align*}
	\k'(\xx,\vp) &= \k\big( \xx, (\vp\cap (P_1\sm P_3))\cup \psi\big) \one\{\xx\in P_1\sm P_3 \},\\
	\k''(\xx,\vp)&:=\begin{cases}
		 \k\big(\xx,(\vp\cap (P_{1,-}\sm P_3))\cup   \psi_1\big) & \xx \in P_{1, -}\sm P_3,\\
		\k\big(\xx,(\vp\cap P_{2,-})\cup   \psi_2\big) & \xx \in P_{2, -},\\
		0 &\text{otherwise}.
	\end{cases} 
\end{align*}
Then,
 $\k',\k''$ satisfy the assumptions of Theorem \ref{thm:Talpha} with $P=P_{1,-}$ and $\a(\xx)=\one\{\xx\notin P_3\}$.
\enl

%
%
\subsection{Extension to unbounded domains}\label{ss:ff_RCM}

\begin{figure}[!htpb]
\centering
\begin{tikzpicture}[scale=0.85, every node/.style={font=\small}]
    \tikzset{
  nuPoint/.style={circle, fill=black, inner sep=1.2pt}
}

\def\nx{5}
\def\ny{5}

\foreach \x/\y in {
  0.25/0.30, 0.70/0.65,
  1.40/0.20, 1.85/0.55,
  2.30/0.40, 2.75/0.70,
  3.35/0.25, 3.80/0.60,
  4.20/0.35
}{
  \node[nuPoint] at (\x,\y) {};
}

\foreach \i in {0,...,4} {
  \foreach \j in {0,...,4} {
    \ifnum\j<2
      \fill[gray!80] (\i,\j) rectangle ++(1,1);
    \else
      \ifnum\j=2
        \ifnum\i<3
          \fill[gray!80] (\i,2) rectangle ++(1,1);
        \fi
      \fi
    \fi
  }
}

\foreach \i in {0,...,4} {
  \fill[gray!80] (\i,0) rectangle ++(1,1);
}

\fill[gray!20] (0,3) rectangle ++(1,1);
\fill[gray!20] (0,4) rectangle ++(1,1);
\fill[gray!20] (1,4) rectangle ++(1,1);
\fill[gray!20] (2,4) rectangle ++(1,1);
\fill[gray!20] (3,4) rectangle ++(1,1);
\fill[gray!20] (4,4) rectangle ++(1,1);
\fill[gray!20] (4,3) rectangle ++(1,1);
\fill[gray!20] (4,2) rectangle ++(1,1);

\fill[white] (3,2) rectangle ++(1,1);
\node at (2.5,2.5) {$Q_{z,1}$};

\foreach \i in {0,...,4} {
  \foreach \j in {0,...,4} {
    \draw (\i,\j) rectangle ++(1,1);
  }
}

\foreach \x/\y in {
	0.3/3.7,
	1.1/4.5,
	2.5/4.8,
	3.7/4.2,
	4.5/4.4,
	4.2/2.8
} {
	\fill (\x,\y) circle (1.2pt);
}

\foreach \x/\y in {
  0.3/0.4,
  1.7/2.6,
  2.4/0.3,
  3.6/1.7,
  1.2/1.2,
  4.2/0.4
} {
  \fill (\x,\y) circle (1.2pt);
}

\foreach \x/\y in {
	1.5/3.6,
	2.9/3.3,
	3.5/2.7,
	3.6/1.7,
	1.2/1.2,
	4.2/0.4
} {
	\fill (\x,\y) circle (1.2pt);
}

\draw[very thick] (1,1)--(1,4) -- (4,4) -- (4,1)--(1,1);

\end{tikzpicture}
\caption{
The set $U$ is the union of the white and the light grey region consisting af all squares succeeding $Q_{z,1}$ lexicographically. The window $W_n$ is the middle $3\times 3$ square, so $U_n$ is the white region and $U\sm U_n$ is the light grey region. The points in the dark grey region are the boundary points $\psi$ while the points in the light grey region are $\nu_n$. The points in the white region is the Poisson process that we are thinning.
}
\label{fig:lex_order}
\end{figure}
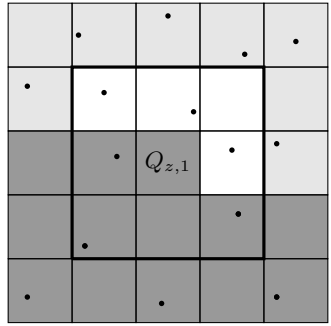

So far, we have constructed thinnings of the RCM only for Gibbs processes $\XX(Q,\psi)$ in bounded domains $Q\in \BB_0$. In this section, we consider the extension to the infinite-volume Gibbs process $\X(U,\psi)$ on an unbounded  Borel set $U\su \Xx$  with  boundary conditions $\psi\in \Nlf_{\Xx\sm U}$. We are going to use this for the proof of Theorem \ref{thm:3} in the setting where $U = \bigcup_{w \succ z} Q_{w,1}$ is the union of all cubes preceding a fixed cube $Q_{z,1}$ in the lexicographic order, see illustration in Figure \ref{fig:lex_order}. 
 In Section \ref{sss:ff_RCM}, we show that the cluster-based thinning in increasing windows converges to give a Gibbs process in unbounded domains. However, this construction does not control what happens when changing boundary conditions  $\psi$ locally. Thus, in Section \ref{sss:dcff}, we give a modified construction that is better able to control the effect of  changing $\psi$. In the proof of Theorem \ref{thm:3}, we will apply the result when $\psi$ is changed inside $Q_{z,1}$, see Figure~\ref{fig:lex_order}.

\subsubsection{A cluster-based thinning construction in infinite-volume domains}\label{sss:ff_RCM}

The idea is to construct $\X(U,\psi)$ as the limit of Gibbs processes $\XX(U_n,\psi)$, where $U_n=U\cap W_n$,  by letting $n\to\ff$. The processes are obtained by thinning a common process $\Phi^*_{U,\psi}$ with underlying Poisson process $\Pds_U$ and given $\Pds_U$, each $(\xx,\rr)\in \Phi^*_{U,\psi}$ consists of an $\xx\in \Pds_U$ equipped with a sequence $\rr=(r_{\xx,\yy})_{\yy\in \Pd_U\cup \psi}$ of i.i.d.\ uniform variables in $[0,1]$.  
Let $\nu_n=\Pd_{U\sm U_n}$. Then,  $\XX(U_n,\psi)$ can be constructed as $T^{\ms{RCM,cl}}_{U_n,\nu_n,\psi}(\Phi^*_{U_n,\nu_n\cup\psi})$, where $\Phi^*_{U_n,\nu_n\cup \psi}$ is the restriction of $\Phi^*_{U,\psi}$ to $U_n$ using the same association of marks with points $\rr=(r_i)_{i\in\Z}=(r_{\xx,\yy})_{\yy\in \Pd_U\cup \psi}$  as the original $\Phi^*_{U,\psi}$.  
The exploration-based graphs $\GGc (\Phi^*_{ U_n,\nu_n\cup\psi },\psi )$ converge when $n\to \ff$ to a well-defined graph $\GGc (\Phi^*_{ U,\psi},\psi)$,   which is distributed as an RCM $\Gamma(U,\psi)$.

\bel\label{lem:tilde_iotan}
Consider a PI with non-negative pair potential $\k$ satisfying \eqref{cond:int_v} and \eqref{cond:perc_RCM}. Let $U\su \Xx$ be a Borel set, $\psi\in \Nlf_{\Xx\sm U}$ and $\nu_n=\Pd_{U\sm U_n}$. Assume $\iota$ is such that $\iota(U_{n+1}\sm U_n)\su (2^{-(n+1)},2^{-n})$ and $\iota(\Xx\sm U) \su (-\ff,0)$.
The connected components of $ \GGc (\Phi^*_{ U_n,\nu_n\cup\psi},\nu_n\cup \psi ) $ not connected to $\nu_n$ are unchanged when increasing $n$. The limiting graph $$ \GGc (\Phi^*_{ U,\psi},\psi):= \bigcup_{n_0\ge 1} \bigcap_{n\ge n_0} \GGc (\Phi^*_{ U_n,\nu_n\cup \psi},\psi ) \stackrel d=\Gamma(U,\psi)$$ is the union of all these connected components. 
\enl

\bep
See Section \ref{ss:proffs_ff}
\enp

We next show that the Gibbs processes 
$T^{\ms{RCM,cl}}_{U_n,\nu_n,\psi}(\Phi^*_{U_n,\nu_n\cup\psi})$ converge to $\X(U,\psi)$.

\bepr\label{prop:RCM_ff}
Consider a PI with non-negative pair-potential $\k$ satisfying \eqref{cond:int_v} and \eqref{cond:perc_RCM}.
Let $A\su W_n$ be a Borel set. On the event $\{\Pd_{A\cap U} \not\lrsa_{\Gamma} \nu_n\}$, where $\Gamma=\GGc (\Phi^*_{ U,\psi},\psi)$,
$$\Trc_{U_m,\nu_m,\psi}(\Phi^*_{U_m,\psi\cup \nu_m})\cap A =\Trc_{U_n,\nu_n,\psi}(\Phi^*_{U_n,\psi\cup\nu_n})\cap A$$ 
almost surely whenever $m\ge n$. Let 
$$\Tf_{U,\psi}(\Phi^*_{U,\nu_n\cup \psi}) = \bigcup_{n_0\ge 1} \bigcap_{n\ge n_0} \Trc_{U_n,\nu_n,\psi}(\Phi^*_{U_n,\psi\cup\nu_n})$$
 denote the limiting process. Then, $\Tf_{\Xx,\es}(\Phi^*_{\Xx,\es})\stackrel d= \X(\Xx,\es)$ and if $\Psi\stackrel d=\X(\Xx,\es)\cap \sm U$, then  $\Psi\cup \Tf_{U,\Psi}(\Phi^*_{U,\Psi})\stackrel d= \X(\Xx,\es)$. 
\enpr

\bep
See Section \ref{ss:proffs_ff}
\enp

\subsubsection{Disagreement coupling for infinite-volume Gibbs processes}\label{sss:dcff}
In Section \ref{sec:weak}, we need to compare  infinite-volume Gibbs processes on $U\su \Xx$ with locally differing boundary conditions. For simplicity, we will assume that the boundary conditions disagree only on  $W_1$ and $W_1\cup U=\es$. Let $\psi\in \Nlf_{\Xx\sm (U\cup W_1)}$ and $\mu_1,\mu_2\in \Nlf_{W_1}$. We compare the two sets of  boundary conditions  $\psi\cup \mu_1$ and $\psi\cup \mu_2$. As in Section \ref{ss:ff_RCM}, we first construct Gibbs processes on finite domains $U_n=U\cap W_n$ and then let $n\to \ff$, but we use a modification of the cluster-based disagreement coupling. 
Assume $ \mu_1\cup \mu_2\su \mu \su W_1$ and $\nu_n = \Pd_{U\sm U_n}$.
We denote the modified  cluster-based disagreement coupling by $ T^{\ms{RCM,cl},2}_{U_n,\mu,  \nu_n,\psi\cup \mu_i}(\Phi^*_{ U_n,\mu\cup\psi\cup \nu_n})$. It starts from $\nu_1=\mu$. In each step, we take $\nu_{k+1}=V_k$ until at some point $V_{k_0}=\es$, i.e.\ we have explored all clusters connected to $\nu_1$.  In this case, we set $V_{k_0}'=\es$ and $\nu_{k_0+1}=\nu_n$. For $k>k_0$, we  set  $\nu_{k+1}=V_{k+1}'$.  Thus, we first construct all clusters connected to $\mu$. We denote the union of all these clusters  by $\mathcal{C}_\mu$. Then, we follow the usual cluster-based exploration scheme of Definition \ref{def:cl_RCM} and explore all clusters connected to $\nu_n$. Finally, the remaining clusters are explored in the order of their $\iota$-smallest point. The resulting graph is denoted $ \GG^{\ms{cl},2} (\Phi^*_{ U_n,\mu\cup\psi\cup \nu_n},\psi\cup \mu) $ and the associated thinning $ T^{\ms{RCM,cl},2}_{U_n,\mu,  \nu_n,\psi\cup \mu_i}(\Phi^*_{ U_n,\mu\cup\psi\cup \nu_n})$. The first failed exploration step slightly deviates from the set-up in Definition \ref{def:GGd}, however, it is easily verified that $ \GG^{\ms{cl},2} (\Phi^*_{ U_n,\mu\cup\psi\cup \nu_n},\psi\cup \mu) \stackrel d= \G(U_n,\psi \cup \mu \cup \nu_n)$ and $ T^{\ms{RCM,cl},2}_{U_n,\mu,  \nu_n,\psi\cup \mu_i}(\Phi^*_{ U_n,\mu\cup\psi\cup \nu_n}) \stackrel d= \X(U_n \psi\cup \mu_i)$.

We want to obtain a thinning  $ T^{\ms{RCM,\ff},2}_{U,\mu,\psi\cup \mu_i}$ as a limit of $ T^{\ms{RCM,cl},2}_{U_n,\mu ,\nu_n,\psi \cup\mu_i}$ when  $n\to 0$. To see that this is well-defined, we first note that by a  proof similar to Lemma \ref{lem:tilde_iotan}, we obtain the following result.
 
\bel
Assume that \eqref{cond:int_v} and \eqref{cond:perc_RCM} hold. Let $\psi\in \Nlf_{\Xx\sm (U\cup W_1)}$, $\mu_1,\mu_2\in \Nlf_{W_1}$ and $\mu_1\cup \mu_2 \su \mu$. Suppose $\CC_\mu\su W_n$.
The connected components of $ \GG^{\ms{cl},2} (\Phi^*_{ U_n,\mu\cup\psi\cup \nu_n},\psi\cup \mu\cup \nu_n) $ not connected to $\nu_n$ are unchanged when increasing $n$. The limiting graph $$ \GG^{\ms{cl},2} (\Phi^*_{ U,\psi\cup \mu},\psi\cup \mu):= \bigcup_{n_0\ge 1} \bigcap_{n\ge n_0} \GG^{\ms{cl},2} (\Phi^*_{ U_n,\mu \cup \psi\cup \nu_n},\psi\cup \mu ) \stackrel d=\Gamma(U,\psi\cup \mu)$$ is the union of all these connected components.
\enl

We next show that the thinnings   $ T^{\ms{RCM,cl},2}_{U,\mu,\nu_n,\psi\cup \mu_i}$ converge when $n\to \ff$. That the thinning of the clusters not connected to $\mu\cup \nu_n$ stabilize when $n\to \ff$ can be shown as in Proposition \ref{prop:RCM_ff}. However, to show that the thinning of $\CC_\mu$ also stabilizes requires an additional proof, which makes use of Corollary \ref{cor:poiss_emb} and the perturbed PIs from Section \ref{ss:pert_RCM}. The resulting thinnings will have the property \eqref{eq:Tcl2_disag} below that exchanging the boundary conditions $\mu_1$ by $\mu_2$ only affects the clusters $\CC_\mu$ connected to $\mu$.

For the proof, we make some extra assumptions. Let $\Gamma$ denote an RCM $\Gamma(U,\es)$ with vertex set $\Pd_U$.
By dominated convergence, if the PI satisfies \eqref{cond:perc_RCM}, then, for every $n \ge 1$,
$$\lim_{n'\to \ff}\int_{U_n } \P(\xx\lrsa_{\Gamma_\xx} \Pd_{U\sm U_{n'}}) \la(\d \xx) =0. $$
Thus, after possibly replacing $W_n$ by a subsequence, we may assume in the proof that 
\begin{equation}\label{eq:RCM2-con1}
\int_{U_n } \P(\xx\lrsa_{\Gamma_\xx} \Pd_{U\sm U_{n+1}}) \la(\d \xx)<2^{-n}.
\end{equation}
Note that this means that Proposition \ref{prop:RCM2_ff} only holds for this subsequence, but this is sufficient for our purposes.
If \eqref{cond:local_int_v} holds, we may further assume
\begin{equation}\label{eq:RCM2-con2}
	\int_{\Xx \sm W_{n+1}} \int_{W_n} \pi(\xx,\yy) \la(\d \xx)\la (\d \yy) < 2^{-n}.
\end{equation}

\bepr\label{prop:RCM2_ff}
Let $U\in \BB$ with $U\cap W_1=\es$, $A\su U$ with $A\in \BB_0$ and $\psi\in \Nlf_{\Xx\sm U}$.
Let $\k$ be a PI with non-negative pair potential satisfying \eqref{cond:int_v}, \eqref{cond:local_int_v} and \eqref{cond:perc_RCM}.
 Assume $\mu_1\cup \mu_2\su\mu$ where $\mu$ is a Poisson process on $W_1\sm U$ with intensity measure $c\la$ for some $c>0$. Assume that the sequence $W_n$ satisfies \eqref{eq:RCM2-con1} and \eqref{eq:RCM2-con2}. Let $\CC_\mu$ denote the union of all clusters connected to $\mu$ in $\GG^{\ms{cl},2} (\Phi^*_{ U,\psi\cup \mu},\psi\cup \mu)$. Then, for $i=1,2$,
\begin{equation}\label{eq:Tcl2Cmu}
 \P\Big(\bigcap_{n\ge 1}\bigcup_{m\ge n} \Big\{ T^{\ms{RCM,cl},2}_{U_m,\mu,\nu_m ,\psi\cup\mu_i}(\Phi^*_{U_m,\psi\cup\mu \cup \nu_m})\cap \CC_\mu \neq  T^{\ms{RCM,cl},2}_{U_n,\mu,\nu_n ,\psi\cup \mu_i}(\Phi^*_{U_n,\psi\cup\mu \cup \nu_n})\cap \CC_\mu \Big\} \Big) =0.
 \end{equation}
If $A\su U_n$, then on the event $\{\Pd_A \cup \CC_\mu \not\lrsa_{\Gamma} \nu_n\}$, where $\Gamma=\GG^{\ms{cl},2} (\Phi^*_{ U,\psi \cup \mu },\psi\cup \mu)$, we have
\begin{equation}\label{eq:Tcl2A}
T^{\ms{RCM,cl},2}_{U_m,\mu,\nu_m,\psi\cup \mu_i}(\Phi^*_{U_m,\psi\cup \mu \cup \nu_m})\cap (A\sm \CC_\mu) =T^{\ms{RCM,cl},2}_{U_n,\mu,\nu_n,\psi\cup \mu_i}(\Phi^*_{U_n,\psi\cup \mu \cup \nu_n})\cap (A\sm \CC_\mu) 
\end{equation}
almost surely whenever $m\ge n$. Let 
$$T^{\ms{RCM},\ff,2}_{U,\mu,\psi\cup\mu_i}(\Phi^*_{U,\psi\cup \mu}) := \bigcup_{n_0\ge 1} \bigcap_{n\ge n_0} T^{\ms{RCM,cl},2}_{U_n,\mu,\nu_n,\psi\cup \mu_i }(\Phi^*_{U_n,\psi\cup \mu\cup \nu_n})$$
denote the limiting process. Then,
\begin{equation}\label{eq:Tcl2_disag}
\P\big(T^{\ms{RCM},\ff,2}_{U,\mu,\psi\cup\mu_1}(\Phi^*_{U,\psi\cup \mu})\cap A \ne T^{\ms{RCM},\ff,2}_{U,\mu,\psi\cup\mu_2}(\Phi^*_{U,\psi\cup \mu})\cap A\big) \le \P(\Pd_A \lrsa_{\Gamma} \mu).
\end{equation}
 If $\Psi\stackrel d=\X(\Xx,\es)\cap (\Xx\sm U)$ with $\Psi\cap W_1\su \mu$, then  
 \begin{equation}\label{eq:Tcl_gibbs}
\Psi\cup T^{\ms{RCM},\ff,2}_{U,\mu,\Psi}(\Phi^*_{U,\Psi\cup\mu})\stackrel d= \X(\Xx,\es).
\end{equation}
\enpr

A similar thinning $T^{\ms{cl},2}_{U_n, \mu, \nu_n,\psi\cup \mu_i}(\Pds_{U_n})$ can be given for the processes considered in Section \ref{sec:const}. Again, we first explore all clusters of $G(\Pd_{U_n})$ connected to $\mu$,  then all clusters connected to $ \nu_n$, and finally the remaining clusters. Under the assumption \eqref{cond:perc2}, we may assume that the sequence $W_n$ (or a subsequence) satisfies
\begin{equation}\label{eq:RCM2-con3}
	\int_{W_n} \P(\xx \lrsa_{\Pd_\xx} \Xx \sm W_{n+1}) \la(\d \xx)< 2^{-n},
\end{equation}
and under \eqref{cond:unif_sim} also
\begin{equation}\label{eq:RCM2-con4}
	\int_{\Xx\sm W_{n+1}}\int_{W_n} \one\{\xx\sim \yy\} \la(\dd \xx) \la(\dd \yy) < 2^{-n}.
\end{equation}
We have the following version of Proposition \ref{prop:RCM2_ff}, whose proof is similar to Proposition \ref{prop:RCM2_ff}, but simpler, so we omit it here.

\bepr\label{prop:2_ff}
Consider a PI satisfying \eqref{cond:PIbound},  \eqref{cond:PI_sim}, \eqref{cond:unif_sim} and \eqref{cond:perc2}.
Assume that the sequence $W_n$ satisfies \eqref{eq:RCM2-con3} and \eqref{eq:RCM2-con4}. Let $A\su U\su \Xx$ with $U\cap W_1=\es$,  $A\in \BB_0$ and $\psi\in \Nlf_{\Xx\sm (U\cup W_1)}$. Assume $\mu_1\cup \mu_2\su \mu$  where $\mu$ is a Poisson process on $W_1\sm U$ of intensity $c\la$ for some $c>0$. Let $\CC_\mu$ denote the union of all clusters connected to $\mu$ in $G(\Pd_U\cup \mu)$. Then, for $i=1,2$
$$\lim_{n\to \ff} \P\Big(\bigcup_{m\ge n} \{ T^{\ms{cl},2}_{U_m,\mu,\nu_m,\psi\cup\mu_i}(\Pds_{U_m})\cap \CC_\mu \neq  T^{\ms{cl},2}_{U_n,\mu,\nu_n,\psi\cup\mu_i}(\Pds_{U_n })\cap \CC_\mu \} \Big) =0.$$
On the event $\{\Pd_A \not\lrsa_{\Pd} \nu_n\}$,
$T^{\ms{cl},2}_{U_m,\mu,\nu_m,\psi\cup \mu_i}(\Pds_{U_m})\cap (A\sm \CC_\mu) =T^{\ms{cl},2}_{U_n,\mu,\nu_n,\psi\cup \mu_i}(\Pds_{U_n})\cap (A\sm \CC_\mu) $ almost surely whenever $m>n$. Let 
$$T^{\ff,2}_{U,\mu,\psi\cup \mu_i}(\Pds_U) : = \bigcup_{n_0\ge 1} \bigcap_{n\ge n_0} T^{\ms{cl},2}_{U_n,\mu,\nu_n,\psi\cup \mu_i}(\Pds_{U_n})$$
denote the limiting process. Then,
\begin{equation}\label{eq:Tcl2_disag}
	\P\big(T^{\ff,2}_{U,\mu,\psi\cup \mu_1}(\Pds_U)\cap A \ne T^{\ff,2}_{U,\mu,\psi\cup \mu_2}(\Pds_U)\cap A\big) \le \P(\Pd_A \lrsa_{\Gamma} \mu).
\end{equation}
If $\Psi\stackrel d=\X(\Xx,\es)\cap (\Xx\sm U)$ and $\Psi\cap W_i\su \mu$, then  $\Psi\cup T^{\ff,2}_{U,\mu,\Psi}(\Pds_U)\stackrel d= \X(\Xx,\es)$. 
\enpr

	\section{A qualitative CLT under weak stabilization}
\label{sec:weak}

The aim of this section is to prove Theorem \ref{thm:3}. Thus, we consider a Gibbs process on  marked Euclidean space $\Xx = \R^d\times \M$ with translation invariant PI.

%
%
\subsection{Gibbs processes with translation-invariant PI}\label{sec:ergodic}
We first show some general properties of Model (L) and (P) when the PI is translation-invariant. 

\bepr\label{prop:ergodic}
 Consider Model (L) or Model (P) with a  translation-invariant PI.  For $U\in \BB$ and $\psi\in \Nlf_{\Xx\sm U}$, $\X(U,\psi)+y\stackrel{d}{=}\X(U+y,\psi + y)$ for any $y\in \R^d$. In particular, $\X(\Xx,\es)$ is stationary. Moreover, $\X(\Xx,\es)$ is mixing and hence ergodic.
\enpr

\bep
We give the argument only for Model (P), since Model (L) is similar. 
If $\kappa$ is translation-invariant, then   $\X(U,\psi) $ satisfies the same GNZ equation \eqref{cond:GNZ_boundary} as $ \X(U+y,\psi + y)$ for $y \in \R^d$. The first claim follows since the Gibbs distribution satisfying the GNZ equation is unique by Theorem \ref{thm:unique1} and \ref{thm:unique2}, respectively.

To see that $\X:=\X(\Xx,\es)$ is mixing, we need to show
\begin{equation}\label{eq:mixing}
\lim_{\|x\|\to \infty}\E[f(\X +x)g(\X)]=\E[f(\X)]\E[g(\X)]
\end{equation}
for any $f,g:\Nlf_{\R^d\times \M} \to [0,1]$ of the form $f(\vp)=f(\vp\cap W_s)$, $g(\vp)=g(\vp\cap W_s)$ for some large $s$. We may assume that $\|x\|$ is so large that $W_s\cap (W_s-x) = \es$.

Construct realizations of  $\X(\Xx,\es)$ as follows: 
Set $\Psi=T_{\Xx, \es}^{\ms{RCM},\ff}(\Phi^*_{\Xx,\es})\cap W_s$ and $\tilde{\Psi}=T_{\Xx, \es}^{\ms{RCM},\ff}(\tilde{\Phi}^*_{\Xx,\es})\cap W_s$, where  $\tilde{\Phi}^*_{\Xx,\es}$ is an independent copy of $\Phi^*_{\Xx,\es}$. 
Then $\Psi\stackrel d= \tilde{\Psi} \stackrel d=\X(\Xx,\es)\cap W_s$ by Proposition \ref{prop:RCM_ff}, 
 and $\mu=\Pd_{W_s}\cup \tilde{\Pd}_{W_s}$ is a Poisson process of intensity $2\la$ with $\Psi\cup\tilde{\Psi}\su \mu$. 
Extend $\Psi$ and $\tilde{\Psi}$ to $\X = \Psi \cup \X_1$ and $\tilde{\X} = \tilde{\Psi} \cup \tilde{\X}_1$ where 
${\X}_1=T^{\ms{RCM},\ff,2}_{\Xx\sm W_s, \mu,\Psi}(\check{\Phi}^*_{\Xx\sm W_s,\mu})$ and $\tilde{\X}_1=T^{\ms{RCM},\ff,2}_{\Xx\sm W_s, \mu,\tilde{\Psi}}(\check{\Phi}^*_{\Xx\sm W_s,\mu})$ with $\check{\Phi}^*_{\Xx\sm W_s,\mu}$ independent  of $\Phi^*_{\Xx,\es}$ and $\tilde{\Phi}^*_{\Xx,\es}$. Then $\X \stackrel d= \tilde{\X} \stackrel{d}{=} \X(\Xx,\es)$ by Proposition \ref{prop:RCM2_ff}. Let $\check{\Gamma}:= \Gamma(\check{\Phi}^*_{\Xx\sm W_s,\mu},\mu)$.
If $\{\mu \not \lrsa_{\check{\Gamma}} \check{\PP}_{W_s-x}\}$, then $\X\cap (W_s-x) = \tilde{\X}\cap (W_s - x)$. Thus, 
\begin{align*}
&\big|	\E[f(\X +x)g(\X)] - \E[f(\X +x)]\E[ g(\X)]|\\   &=|	\E[f(\X\cap (W_s-x) +x)(g(\Psi) -\E[ g(\Psi)])]|\\
&\le \big|	\E[f(\tilde{\X}\cap (W_s-x) +x)(g(\Psi)-\E[ g(\Psi)])]\big| + 2\P(\mu \lrsa_{\check{\Gamma}} \check{\PP}_{W_s-x})\\
&= \big|	\E[\E[f(\tilde{\X}_1\cap (W_s-x) +x)\mid \Psi] (g(\Psi)-\E[ g(\Psi)])]\big| + 2\P(\mu  \lrsa_{\check{\Gamma}} \check{\PP}_{W_s-x})\\
&=\big|\E[f(\tilde{\X}_1 \cap (W_s-x) +x)]	\E[ g(\Psi)-\E[ g(\Psi)]]\big| + 2\P(\mu  \lrsa_{\check{\Gamma}} \check{\PP}_{W_s-x})\\
&=2\P(\mu  \lrsa_{\check{\Gamma}} \check{\PP}_{W_s-x}). 
\end{align*}
To obtain the second last equality, note that, even though the realization of $\mu$ affects the thinning algorithm and is correlated with $\Psi$, the conditional distribution of $\tilde{\X}_1$ given $\Psi $ is always $\X\big(\Xx\sm W_s, \tilde \Psi\big)$, which is independent of $\Psi$.  Since $\lim_{\|x\|\to \ff} \P(\mu  \lrsa_{\check{\Gamma}} \check{\PP}_{W_s-x}) =0$, we obtain \eqref{eq:mixing}. 
\enp

\subsection{Proof of Theorem \ref{thm:3}}

 Throughout this section we let $\X$ denote the infinite-volume Gibbs process on $\Xx=\R^d\times \M$ from Model (L) or Model (P).    
 Recall the notation from Section \ref{ss:thm3}. The key observation in the proof is that 
 $
	H_n(\X) - \E[H_n(\X)]= \sum_{i\le k_n} \De_{i,n}, 
	$
where 
$$\De_{i,n} := \De_{z_i,n}:= \E[H_n(\X) | \FF_{i,n}] - \E[H_n(\X) | \FF_{i-1,n}]. $$
For each $n$, the sequence $ \E[H_n(\X) | \FF_{i,n}]$ for $i=1,\ldots,k_n$ defines a martingale with respect to $\FF_{i,n} $. 
A CLT for $H_n(\X)$ can therefore be derived by a standard martingale method \cite{mcleish}. Conditions for this method to apply are given in \cite[Prop. 5.1]{hos}. Theorem \ref{thm:3} will be proved by verifying these conditions.

 We define sets
\begin{equation*}
	V_z = \bigcup_{w\preceq z} Q_{w,1},\quad U_z = \bigcup_{w \succ z} Q_{w,1},
\end{equation*}
where $\preceq$ refers to the lexicographic ordering on $\Z^d$. In Figure \ref{fig:lex_order}, $V_z$ corresponds to the dark grey set and $U_z$ corresponds to the union of the white and the light grey set. Define $V_{z,n}:=V_z \cap Q_n$ and $U_{z,n}:=U_z\cap Q_n$. In the following, for $z\in \Z^d$ we let $z_{-}$ denote the point in $\Z^d$ that comes just before $z$ in the lexicographic order.

Fix a $z\in \Z^d$. For the proof, we construct $\X(\Xx,\es)$ using the disagreement coupling construction for Gibbs processes in unbounded domains from Section \ref{sss:dcff}. In the case of Model (L), we first construct $\X=T^{\ff}_{\Xx,\es}(\Pds)$. Let ${\PP}^*_1$ and $\Pds_2$ be  independent copies of $\Pds$. Define $\mu_z=(\Pd\cup \Pd_2 ) \cap Q_{z,1}$ and
\begin{align*}
&\eta_z=\XX \cap V_z, \quad  \XX_z = \eta_z \cup T^{\ff,2}_{U_z,\mu_z,\eta_z}((\PP_1)^*_{U_z}),\\
 &\eta_{z_-}=\XX\cap V_{z_-}, \quad \zeta_{z_-}= T^\ff_{U_{z_-},\eta_{z_-}} ((\Pd_2)^*_{U_{z_-}})\cap Q_{z,1},\quad \XX_{z_-} = \eta_{z_-} \cup \zeta_{z_-} \cup T^{\ff,2}_{U_{z},\mu_{z},\eta_{z_-}\cup \zeta_{z_-}}((\PP_1)^*_{U_z}).
 \end{align*}
 For the thinnings, we use the sequence $W_n = Q_{z,n}$ and an ordering $\iota_z$ where  $\iota_z(W_{n+1}\sm W_{n}) \su (2^{-n-1},2^{-n})$. Then,  $\XX_z\sim \XX_{z_-} \sim \XX$ by Proposition \ref{pr:chimera} and \ref{prop:2_ff}. 

In the case of a PI of Model (P), we similarly  construct $\X=T^{\ms{RCM},\ff}_{\Xx,\es}(\Phi_{\Xx,\es}^*)$.  Let $(\Phi_2)^*_{U_{z_-},\Pd_{\Xx\sm U_{z_-}}}$ and  $(\Phi_1)^*_{U_{z},\Pd_{\Xx\sm U_{z}} \cup \mu_z}$ be independent of each other and of $\Phi_{\Xx,\es}^*$. 
Here, $\mu_z=(\Pd\cup \Pd_2 ) \cap Q_{z,1}$ where $\Pd$ and $\Pd_2$ are the Poisson processes underlying $\Phi_{\Xx,\es}^*$ and $(\Phi_2)^*_{U_{z_-},\Pd_{\Xx\sm U_{z_-}}}$, respectively.  Note that the second subscript $\psi$ in $\Phi_{U,\psi}^*$ is only used to indicate how coordinates are identified with points, and hence, this is not part of the independence requirement. Define
\begin{align*}
	&\eta_z=\XX \cap V_z, \quad  \XX_z = \eta_z \cup T^{\ms{RCM},\ff,2}_{U_z,\mu_z,\eta_z}((\Phi_1)^*_{U_{z},\Pd_{\Xx\sm U_{z}} \cup \mu_z}),\\
	&\eta_{z_-}=\XX\cap V_{z_-}, \quad \zeta_{z_-}= T^{\ms{RCM},\ff}_{U_{z_-},\eta_{z_-}} ((\Phi_2)^*_{U_{z_-},\Pd_{\Xx\sm U_{z_-}}})\cap Q_{z,1},\\
	&\XX_z = \eta_{z_-} \cup \zeta_{z_-} \cup T^{\ms{RCM},\ff,2}_{U_{z},\mu_{z},\eta_{z_-}\cup \zeta_{z_-}}((\Phi_1)^*_{U_{z},\Pd_{\Xx\sm U_{z}} \cup \mu_z}).
\end{align*}
We take $W_n$ and  $\iota_z$ as in Model (L).
Then, $\XX_z'\sim \XX_{z_-}' \sim \XX$ by Proposition~\ref{prop:RCM_ff} and \ref{prop:RCM2_ff}.

In both cases we get 
$ \E[H_n(\XX) | \eta_{z_-} ]= \E[H_n(\XX_{z_-}) | \eta_{z_-} ] = \E[H_n(\XX_{z_-}) | \eta_z ]$.
Thus, we can write the martingale differences as 
$$\Delta_{z,n} = \E[H_n(\XX_z) \mid \eta_z ]-\E[H_n(\XX_{z_{-}}) \mid \eta_{z_-}] = \E[H_n(\XX_z) - H_n(\XX_{z_{-}}) \mid \eta_z ].$$
By construction, $\XX_z$ and $\XX_{z_-}$ both equal $\eta_{z_-}$ on $V_{z_-}$. On $Q_{z,1}$ they may differ, but $(\XX_z\cup \XX_{z_-})\cap Q_{z,1} \su \mu_z$. Since they are extended over $U_z$ by the same cluster based thinning algorithm starting from $\mu_z$, they will agree on $\Xx\sm Q_{z,n}$ if $\{\mu_z \not\lrsa_{(\Pd_1)_{U_z}\cup \mu_z} (\Pd_1)_{U_z\sm Q_{z,n}}\}$ or $\{\mu_z \not\lrsa_{\Gamma} (\Pd_1)_{U_z\sm Q_{z,n} }\}$, respectively, where $\Gamma:= \Gamma((\Phi_1)^*_{U_{z}, \mu_z},\mu_z)$. In particular,  \eqref{cond:perc2} and \eqref{cond:perc_RCM}, respectively, ensure that there is almost surely an $N_z(\omega)<\ff$ such that $\XX_z\sm Q_{z,N_z} = \XX_{z_-} \sm Q_{z,N_z}$. 

With these constructions, the proof is almost identical to the proof of \cite[Thm. 5.3]{hos}. We give the details for completeness.

\bep[Proof of Theorem \ref{thm:3}]
The result follows from \cite[Prop 5.1]{hos}  with $a_n=n+\sqrt n$, $\YY=\XX$ and $\rho(a_n)=2\sqrt n$ if we can verify Condition (i) and (ii) of that proposition. Condition (i) is exactly assumption \eqref{cond:cond_moment}. It remains to show (ii). Since $\X$ is stationary and ergodic by Proposition \ref{prop:ergodic} this condition states:
\begin{itemize}
	\item[(ii)] 
For every $z\in \Z^d$, there is a random variable $\De_z = \De_z(\X)$ such that $\De_{z,n} \to \De_z$ in  $L^2$. 
The convergence must be uniform in the sense that
\begin{equation}
\label{eq:unif_conv}
\lim_{n\to \ff} \sup_{z\in Z_{n-\sqrt{n}}} \|\De_z - \De_{z,n}\|_{L^2} =0.
\end{equation}
The limit is shift invariant meaning that for any $z_0 \in \Z^d, \De_{z + z_0} (\X + z_0) = \De_z(\X)$.
\end{itemize}

Fix $z\in \Z^d$. As noted above, there is  an $N_z(\omega)$ such that $\XX_z\sm Q_{z,N_z}=\XX_{z_{-}}\sm Q_{z,N_z}$. Hence,
\begin{align*}
	 H(\XX_z) - H(\XX_{z_{-}}) &:=\lim_{n\to \infty} (H_n(\XX_z) - H_n(\XX_{z_{-}}))\\
	& =\lim_{n\to \infty} (H_n(\XX_z) - H_n(\XX_z\sm Q_{z,N_z}) - H_n(\XX_{z_-}) + H_n(\XX_{z_{-}}\sm Q_{z,N_z}))
\end{align*}
exists almost surely by the weak stabilization condition \eqref{cond:weak}. Moreover,
$$\E[\E[H(\XX_z) - H(\XX_{z_{-}})|\eta_z]^4] \le M:=\sup_{n\ge 1, z\in\Z^d\cap Q_n}\E[\Delta_{z,n}^4]$$
by the conditional Fatou's Lemma. 

 We now show that 
$$\E[H_n(\XX_z) - H_n(\XX_{z_{-}}) \mid \eta_z ] \to \E[H(\XX_z) - H(\XX_{z_{-}})\mid \eta_z] =:\Delta_z \text{ in } L^2 .$$
First, note that \eqref{cond:weak} implies that for any $z\in \Z^d$, $l\ge 1$, and $\eps>0$, there exists a random variable $R(z,l,\eps)$ such that whenever $Q_{z,R(z,l,\eps)} \su Q_{w,m}$,
$$|(H(\XX_z)-H(\XX_z\sm Q_{z,l})) -( H(\XX_z\cap Q_{w,m})-H(\XX_z\cap Q_{w,m} \sm Q_{z,l}))|<\eps,$$
see \cite[Proof of Thm. 5.3]{hos} for details. Similarly, there exists a random variable $R_-(z,l,\eps)$ such that whenever $Q_{z,R_-(z,l,\eps)} \su Q_{w,m}$,
$$|(H(\XX_{z-})-H(\XX_{z-}\sm Q_{z,l})) -( H(\XX_{z_-}\cap Q_{w,m})-H(\XX_{z_-}\cap Q_{w,m} \sm Q_{z,l}))|<\eps.$$
The distributions of $R(z,l,\eps)$ and $R_-(z,l,\eps)$ are independent of $z$.

Let $\eps>0$ be given and fix $z\in \Z^d$. For $l,n\in \N$, we define events
$$E_{z,l} :=\{\XX_z\sm Q_{z,l} =\XX_{z_{-}}\sm Q_{z,l} \} , \quad F_{z,l,n} :=\big\{\max\{R(z,l,\eps),R_-(z,l,\eps)\}\le \sqrt n\big\}.$$
Choose $l$ such that $\P(N_z>l)< \eps$. Then
$\P(E_{z,l})> 1-\eps$. For any $n \ge l$ so large that 
 $\P(F_{z,l,n}) \ge 1-\eps$, the conditional Jensen inequality and the Cauchy-Schwarz inequality yield
\begin{align*}
	&\big\| \E[H_n(\XX_z) - H_n(\XX_{z_{-}}) \mid \eta_z ]- \E[H(\XX_z) - H(\XX_{z_{-}}) \mid \eta_z ]\big\|_{L^2}^2\\ &\le 
	\big\| \one_{E_{z,l} \cap F_{z,l,n}}((H_n(\XX_z) - H_n(\XX_{z_{-}})) -( H(\XX_z) - H(\XX_{z_{-}}) ))\big\|^2_{L^2}\\ &\quad +4\E[\Delta_{z,n}^4 + \Delta_z^4]^{1/2} (\P(E_{z,l}^c) + \P(F_{z,l,n}^c))^{1/2}\\
	&\le 4\big\|\one_{F_{z,n,l}} (H_n(\XX_z) - H_n(\XX_z\sm Q_{z,l}) - H(\XX_z) + H(\XX_z\sm Q_{z,l}) )\big\|^2_{L^2} \\
	&\quad+4\big\|\one_{F_{z,n,l}}( H_n(\XX_{z_{-}}) -H_n(\XX_{z_{-}}\sm Q_{z,l}) - H(\XX_{z_{-}}) + H(\XX_{z_{-}}\sm Q_{z,l}) )\big\|^2_{L^2} +8 M^{1/2} \eps^{1/2}\\
	&\le 8\eps^2 +4\sqrt{6} M^{1/2} \eps^{1/2}
\end{align*}
Since the distributions of $N_z$, $R(z,\eps,l)$ and $R_-(z,\eps,l)$ can be assumed independent of $z$, the bound is uniform for all $z\in Q_{n-\sqrt n}$. Since $\eps$ was arbitrary, this shows \eqref{eq:unif_conv}. The fact that $\Delta_{z+z_0}(\XX+z_0) = \Delta_z(\XX)$ follows immediately by construction.
\enp

	\section{Palm couplings}
\label{sec:pac}
In this section, we discuss general results on Palm couplings, i.e., couplings involving some random measure and its Palm kernel. These results will play a central role in the proof of both the Poisson and the quantitative normal approximation results in Theorems \ref{thm:1} and \ref{thm:2}, respectively. Palm couplings were already used successfully for Poisson and normal approximation in \cite{BSY202,dp} and \cite{chen,hos}, respectively.
We first explain the overall idea behind this method. After that we show how to obtain Palm couplings using disagreement coupling. We will mainly focus on the case of pair-wise interactions, i.e., Model (P). Model (L) is similar and a bit easier. Hence, we only add some brief remarks on the steps that are substantially different. Throughout Section \ref{sec:pac}, \ref{sec:papx} and \ref{sec:malliavin_stein}, we leave out the subscript $n$ on $g_n$, $\La_n$ and $\Xi_n$.

%
%
Recall that we consider the random measure
	\begin{align}
	\Xi := \Xi[\XX]& :=\sum_{\xx \in \XX } g(\xx, \XX)\,\delta_{\xx} \label{eqn:poixidefXi}
\end{align}
 based on the score function $g$. This random measure perspective makes the task amenable to the use of Palm theory. We recall that the random measures $\{\Xi_{\xx}\}_{\xx \in \Xx}$ form a \emph{Palm kernel} of a random measure $\Xi$ if 
\begin{align}
	\E\Big[ \int_{\Xx} f(\xx,\Xi) \,\Xi(\d \xx)\Big]=\int_{\Xx} \E [f(\xx,\Xi_{\xx})] \,\La(\d \xx) \label{def:palm1}
\end{align}
for all non-negative measurable $f$, where $\La := \E[\Xi]$ denotes the intensity measure. We refer to \cite{randMeas} for details on Palm theory.

%
%
The Palm kernel is defined on a distributional level. Hence, the definition alone does not give any relation between the Palm kernel and the original process on the level of realizations. However, in the Poisson approximation setting where $\Xi$ is a point process, a key finding is that if one can provide a coupling of the original process $\Xi$ and its Palm kernel $\Xi_\xx$ in such a way that they are identical outside a small neighborhood of $\xx$, then $\Xi$ is close to a Poisson process. More precisely, from \cite[Theorem 3.1]{BSY202}, we obtain that
\begin{align}
	\label{eq:dkr}
	d_{\ms{KR}}\big(\Xi\cap {Q_n}, \mc M \cap {Q_n}\big) \le 2 \int_{Q_n} \E\big[(\Xi \Delta \Xi_{\xx}^!)(Q_n)\big] \La(\d \xx),
\end{align}
where $d_{\ms{KR}}$ is the Kantorovich-Rubinstein distance 
from \eqref{eq:ddkr}, $\mc M$ is a Poisson process of intensity $\La$, and $\Xi_{\xx}^!:=\Xi_{\xx}\sm\{\xx\}$ is a reduced Palm kernel. Likewise, in the normal approximation setting, a result of \cite{chen}, stated precisely in Theorem \ref{thm:rollin} below, provides a bound on $$\dk\bigg(\frac{\Xi(Q_n) -\E[\Xi(Q_n)]}{\var(\Xi(Q_n))^{1/2}} ,N(0,1)\bigg)$$ 
in terms of a coupling between $\Xi$ and $\Xi_{\xx}, \xx\in Q_n$. This motivates the construction of such a coupling. 

We write $\{\XX_{\xx}^\Xi\}_{\xx\in Q_n}$ for a Palm kernel of $\XX$ with respect to $\Xi$ in $Q_n$. That is, for all measurable $f:Q_n \ti \mathbf N \to [0,\ff)$,
\begin{align}
	\E \int_{Q_n} f(\xx,\XX) \,\Xi(\d \xx)= \int_{Q_n} \E\big[f(\xx,\XX_{\xx}^\Xi)\big] \, \La(\d \xx). \label{def:palm2}
\end{align}
We note that assumption \eqref{eq:m2} ensures the $\s$-finiteness of $\La$, and thereby the existence of a Palm kernel $\{\XX_{\xx}^\Xi\}_{\xx \in Q_n}$, see \cite{randMeas}. 
From \eqref{def:palm1} and \eqref{def:palm2} it follows that $\{\Xi[\XX_{\xx}^\Xi ] \}_{\xx\in Q_n}$ is a Palm kernel $\{\Xi_{\xx}\}_{\xx \in Q_n}$ of $\Xi$. Hence, a coupling between $\Xi$ and $\Xi_{\xx}$ will be provided via a coupling between $\X$ and $\X_{\xx}^{\Xi}$. 
The following lemma from \cite{hos} shows that the Palm process $\XX_{\xx}^{!,\Xi}:=\XX_{\xx}^\Xi\sm \{\xx\}$ is a Gibbs process and hence a coupling with $\X$ can be obtained by disagreement coupling.

%
%
\bel[Lemma 6.2 in \cite{hos}]
\label{lem:palmgibbs}
Let $\XX$ be a Gibbs process satisfying the conditions in Theorem \ref{thm:1}. Assume that $g$ satisfies the hereditary condition \eqref{eq:here}. Then, for almost all $\xx \in Q_n$, the reduced process $\XX_{\xx}^{!,\Xi}$ is a Gibbs process with PI $\k_{\xx}$ given by
$$\k_{\xx}(\yy,\vp):=\k(\yy,\vp \cup \{\xx\}) \f{g(\xx,\vp \cup \{\xx,\yy\})}{g(\xx,\vp \cup \{\xx\})}$$
where $0/0:=1$.
\enl

We now provide the precise construction of the Palm coupling in our setting. The specific coupling construction proceeds loosely as in \cite{dp} and \cite{hos}, but there are some differences. Most importantly, for Model (P) we use the RCM-adapted thinning algorithm from Section \ref{sec:RCM}.

\def\Tc{T^{\mathsf{RCM},\mathsf{dc}}}
\def\Tcx{T^{x,\mathsf{RCM},\mathsf{dc}}}
\def\Tcp{T^{\prime,\mathsf{RCM},\mathsf{dc}}}
\def\Tcpx{T^{\prime,x,\mathsf{RCM},\mathsf{dc}}}

%
%
Let $r:=6s := 6c_{\ms St} \log (\la(Q_n))$ with $c_{\ms St}>0$ a constant that will be specified later. While in the random measure $\Xi\cap Q_n$ we only sum over points in $Q_n$, the dependence through the score functions means that $\Xi_n$ also depends on the points outside $Q_n$. Therefore, we introduce the enlarged window $Q_n^+ = B_{2r}(Q_n)$. 
 
 We first construct $\X$. Let $\tilde{\X}$ be a Gibbs process on $\Xx$ with PI $\k$. Outside $Q_n^+$ we will use the boundary condition $\Psi:=\tilde \XX \sm Q_n^+$. Inside $Q_n^+$ we construct $\X$ by thinning an independent process $\Phi_{Q_n^+,\es}^*$. We put $\Psi$ as part of the PI, that is, we use the PI 
 \begin{equation}\label{eq:k'}
 \k'(\xx,\vp) = \k(\xx,(\vp\cap {Q_{n}^+})\cup \Psi)\one\{\xx\in Q_n^+\}
 \end{equation}
 and empty boundary conditions. This is a perturbed PI of the form \eqref{eq:tildekk'}, and hence we can construct it as a thinning of the form
 $\bar{T}^{\prime,\ms{RCM,cl}}_{Q_n^+,\es,\es}$ with dominating PI $\k$, where the "prime" refers to the PI $\k'$. Thus, we set 
\begin{equation}
\XX:=\bar{T}^{\prime,\ms{RCM,cl}}_{Q_n^+,\es,\es}(\Phi_{Q_n^+,\es}^*) \cup \Psi.\label{eqn:couplingX}
\end{equation}

Next, we construct $\XX_{\xx}^\Xi$. Let $\tilde \XX_{\xx}^\Xi$, $\xx \in Q_n$, be an independent family of Palm processes of $\XX$ with respect to $\Xi$ that are independent of $\tilde \X$ and $\Phi_{Q_n^+,\es}^*$. This is possible by the existence of uncountable product measures, see \cite[Corollary 6.18]{kallenberg}. Let $\Psi_{\xx}=\tilde \XX_{\xx}^\Xi \cap (B_s(\xx)\cup (\Xx \sm Q_n^+) )$ be the boundary conditions and define 
\begin{align*}
	{\k'_{\xx}}(\zz,\vp):=\k_{\xx}(\zz,\vp \cap (Q_n^+ \sm B_s(\xx)) \cup \Psi_{\xx} ) \one\{\zz \in Q_n^+ \sm B_s(\xx)\}.
\end{align*}
If $R(\xx, \tilde \XX_{\xx}^\Xi)\le s$, then 
 \begin{align}\label{eq:k'x}
 	{\k'_{\xx}}(\zz,\vp)=\k(\zz,\vp \cap (Q_n^+ \sm B_s(\xx)) \cup \Psi_{\xx} ) \one\{\zz \in Q_n^+ \sm B_s(\xx)\}.
 \end{align}
 We define 
\begin{align}
	\XX_{\xx}^{\Xi}:= \begin{cases}
		\bar{T}^{\prime,\xx,\ms{RCM,cl}}_{Q_n^+,\es,\es}(\Phi_{Q_n^+,\es}^{*})\cup \Psi_{\xx} & \text{ if $R(\xx, \tilde \XX_{\xx}^\Xi)\le s$,}\\
		{\tilde \XX_{\xx}^\Xi} & \text{ if $R(\xx,\tilde \XX_{\xx}^\Xi)> s$,}
	\end{cases} 
	\label{eqn:coupling}
\end{align} 
where $	\bar{T}^{\prime,\xx,\ms{RCM,cl}}_{Q_n^+,\es,\es}$ denotes thinning with respect to the PI $\k_{\xx}'$ with the same dominating RCM as $\bar{T}^{\prime,\ms{RCM,cl}}_{Q_n^+,\es,\es}$. 
It follows analogously to the argument given in \cite{hos} that $\XX_{\xx}^{\Xi}$ is a Palm process of $\XX$ with respect to $\Xi$ at $\xx$. 
In the thinnings, we always use a cube-wise ordering as introduced in Section \ref{ss:pert}

%
%
We note that $\k'$ and $\k'_{\xx}$ given by \eqref{eq:k'} and \eqref{eq:k'x} are in the form of Lemma \ref{lem:pa1} with $P_1=Q_n^+$, $P=P_{1,-} =Q_n^+\sm B_s(\xx)$, $P_3=P_{2,-}=\es$, $\psi=\Psi$ and $\psi_1=\Psi_{\xx}$. Hence, for $A:=Q_n\sm B_{4s}(\xx)$, Lemma \ref{lem:pa1} and Theorem~\ref{thm:Talpha}  yield
\begin{align}\nonumber
&\P(\X\cap (Q_n\sm B_{4s}(\xx)) \ne \X_{\xx}^{!,\Xi}\cap (Q_n\sm B_{4s}(\xx))) \le \P(R(\xx,\X_{\xx}^\Xi)>s) + \P(\Ed) \\ \nonumber
&+ \P(\EpRCM(Q_n\sm B_{4s}(\xx),s))+ 4\int_{B_{2s}(Q_n) \sm B_{2s}(\xx)}\P(\zz\lrsa_{\G_{\zz}(Q_n^+,\Psi\cup \Psi_{\xx})} \Pd_{B_s(\xx)}\cup (\Psi \De \Psi_{\xx}) )\la(\d \zz)\\
& + \int_{(Q_n^+ \sm B_{2s}(Q_n))\cup B_{2s}(\xx)} \int_{B_s(Q_n)\sm B_{3s}(\xx)}\pi(\zz,\yy) \la(\d \zz) \la(\d \yy),\label{eq:coupling_bound}
\end{align}
where the expectation on the right-hand side is with respect to $\Psi$ and $\Psi_{\xx}$. For Model (L), we get a similar bound using Theorem \ref{thm:disagree} instead of Theorem \ref{thm:Talpha}. The following lemma will prove helpful to control potential long-range effects both for Poisson and normal approximation.

\bel[Bound on disagreement between $\X$ and $\X_\xx^{!,\Xi}$]
\label{lem:1pei}
Let $m \ge 1$. Let $g: \Xx \to [0,\ff)$ be measurable and satisfying \eqref{eq:stop}, \eqref{eq:s2}, \eqref{eq:here} and \eqref{eq:m2}. Then, with $c_{\ms{St}} \ge 2(m+2)\max(c_{\ms{es},2}^{-1},c_{\ms{v},2}^{-1},c_{\ms{D},2}^{-1})$ in the definition of $s$,
$$\int_{Q_n} \P(\X\cap (Q_n\sm B_{4s}(\xx)) \ne \X_{\xx}^{!,\Xi_n}\cap (Q_n\sm B_{4s}(\xx))) \La_n (\d \xx)\in O(\la(Q_n)^{-m}).$$
\enl
\bep
We deal with the terms in \eqref{eq:coupling_bound} separately. 
\medskip

\ni{\bf First term.}
We use \eqref{def:palm2} and the GNZ equation \eqref{cond:GNZ} to see that 
 \begin{align*}
	\int_{Q_n}\P(R(\xx,\X_{\xx}^\Xi)>s) \La(\d \xx) &\le \int_{Q_n} \E\big[g(\xx, \XX\cup\{\xx\}) \one\{R(\xx,\X\cup \{\xx\})>s\}] \la(\d \xx)\\
	& \le \Big(\int_{Q_n} \E[g(\xx,\X\cup\{\xx\})^2] \la(\d \xx)\Big)^{1/2} \Big(\int_{Q_n}\P(R(\xx,\X\cup\{\xx\})>s) \la(\d \xx) \Big)^{1/2}.
\end{align*}
By the stabilization and moment conditions \eqref{eq:s2} and \eqref{eq:m2}, the right-hand side is of order $O(\la(Q_n)^{-m-1 })$.
\medskip

\ni{\bf Second term.}
By the Poisson concentration property, we have $\P(\Ed) \in O(\la(Q_n)^{-m-1})$, and therefore, since $\La(Q_n)\in O(\la(Q_n))$ by the GNZ-equation \eqref{cond:GNZ} and \eqref{eq:m2}, $\La(Q_n)\P(\Ed) \in O(\la(Q_n)^{-m})$.
\medskip

\ni{\bf Third term.}
Using the GNZ equation \eqref{cond:GNZ} as in the first step, we see 
\begin{align*}
	&\int_{Q_n}\P\big(\EpRCM(Q_n\sm B_{4s}(\xx),s)\big) \La(\d \xx)\\ 
	&\le \Big(\int_{Q_n} \E\big[g(\xx, \XX\cup\{\xx\})^2] \la(\d \xx)\Big)^{1/2} \Big( \int_{Q_n} \P\Big(\Pd_{Q_n\sm B_{4s}(\xx)}\lrsa_{\G(Q_n^+,\es)} \Pd_{(Q_n^+\sm B_s(Q_n))\cup B_{2s}(\xx)}\Big) \la (\d \xx) \Big)^{1/2}\\
	&\le \Big(\int_{Q_n} \E\big[g(\xx, \XX\cup\{\xx\})^2] \la(\d \xx)\Big)^{1/2} \Big( \int_{Q_n}\int_{Q_n\sm B_{4s}(\xx)} \P\Big(\zz \lrsa_{\G_{\zz}(Q_n^+,\es)} \Pd_{(Q_n^+\sm B_s(Q_n))\cup B_{2s}(\xx)}\Big)\la (\d \zz) \la (\d \xx) \Big)^{1/2}
\end{align*}
By the exponential cluster size property \eqref{cond:sharp_vv}, the second factor is of order $\la(Q_n)^{(-2m-2)/2}$. Thus, the total order is $O(\la(Q_n)^{-m-1/2})$. For Model (L) we use Condition \eqref{cond:sharp} instead.
\medskip

\ni{\bf Fourth term.} The fourth term is bounded by
\begin{align}\label{eq:fourth}
	4\int_{B_{2s}(Q_n) \sm B_{2s}(\xx)}\hspace{-.1cm}\P\big(\zz\lrsa_{\G_{\zz}(Q_n^+,\Psi )} \Pd_{B_s(\xx)}\cup\Psi \big)\la(\d \zz)+ 4\int_{B_{2s}(Q_n) \sm B_{2s}(\xx)}\hspace{-.1cm}\P\big(\zz\lrsa_{\G_{\zz}(Q_n^+, \Psi_{\xx})} \Psi_{\xx} \big)\la(\d \zz).
\end{align}
 To bound the first integral, we use that $\Psi$ can be assumed to be contained in a Poisson process $\tilde{\PP}$ of intensity $\la$ independent of $\Pd$ and hence, in order for $\zz$ to be connected to $\Pd_{B_s(\xx)}\cup\Psi $ in $\Gamma_{\zz} (Q_n^+,\Psi)$, it must be connected to a point in $\Pd_{B_s(\xx)}\cup \tilde{\PP}_{\Xx\sm (Q_n^+\sm B_s(\xx))}$ inside $\Gamma_{\zz} (Q_n^+ \sm B_s(\xx), \Pd_{B_s(\xx)}\cup \tilde{\PP}_{\Xx\sm (Q_n^+\sm B_s(\xx))} )$. The $\La$-integral of the first integral in \eqref{eq:fourth} is bounded by
 \begin{align*}
4\int_{Q_n} \int_{B_{2s}(Q_n) \sm B_{2s}(\xx)} 2\P\big(\zz \lrsa_{\Gamma_{\zz} (\Xx,\es )} \Pd_{\Xx \sm (Q_n^+\sm B_s(\xx))}\big)\la(\d\zz) \La(\d \xx).
 \end{align*}
 Using \eqref{cond:sharp_vv} for Model (P) and \eqref{cond:sharp} for Model (L), we get that this is of order $O(\la(Q_n)^{-2m-2})$.

For the second integral in \eqref{eq:fourth}, since $\Psi_\xx = \tilde{\X}_\xx^\Xi \sm (Q_n^+\sm B_s(\xx))$, we have from \eqref{def:palm2}, \eqref{def:palm1} and \eqref{cond:GNZ} that
\begin{align*}
 &\int_{Q_n}\int_{B_{2s}(Q_n) \sm B_{2s}(\xx)}\P\Big(\zz\lrsa_{\G_{\zz}(Q_n^+, \Psi_{\xx})} \Psi_{\xx} \Big)\la(\d \zz) \La(\d \xx ) \\
&\le\E\int_{Q_n}\int_{B_{2s}(Q_n) \sm B_{2s}(\xx)}\P\Big(\zz\lrsa_{\G_{\zz}(Q_n^+,\{\xx\}\cup \tilde{\XX}_{\Xx\sm (Q_n^+ \sm B_s(\xx))})} \{\xx\}\cup \tilde{\XX}_{\Xx\sm (Q_n^+ \sm B_s(\xx))} \Big)g(\xx,\tilde{\X}\cup \{\xx\}) \la(\d \zz) \la(\d \xx )\\
&\le \bigg( \int_{Q_n}\int_{B_{2s}(Q_n) \sm B_{2s}(\xx)} \E[ g(\xx,\tilde{\X}\cup \{\xx\})^2] \la(\d \zz) \la(\d \xx ) \bigg)^{1/2}\\
&\times \bigg( \int_{Q_n}\int_{B_{2s}(Q_n) \sm B_{2s}(\xx)} \P\Big(\zz\lrsa_{\G_{\zz}(Q_n^+,\{\xx\}\cup \tilde{\PP}_{\Xx\sm (Q_n^+ \sm B_s(\xx))})} \{\xx\}\cup \tilde{\PP}_{\Xx\sm (Q_n^+ \sm B_s(\xx))} \Big) \la(\d \zz) \la(\d \xx )\bigg)^{1/2},
\end{align*}
where we used that $\tilde{\XX}\su \tilde{\PP}$ for some Poisson process $\tilde{\PP}$ independent of $\G(Q_n^+,\es)$. The last factor may be bounded by
\begin{align*}
	&\bigg( \int_{Q_n}\int_{B_{2s}(Q_n) \sm B_{2s}(\xx)} (2\P(\zz\lrsa_{\G_{\zz}(\Xx, \es)} {\PP}_{\Xx\sm B_s(\zz)} ) + \P(\zz\lrsa_{\G_{\zz}(Q_n^+,\{\xx\})} \xx ) )\la(\d \zz) \la(\d \xx )\bigg)^{1/2}\\
	&\le \bigg( \int_{Q_n}\int_{B_{2s}(Q_n) \sm B_{2s}(\xx)} (3\P(\zz\lrsa_{\G_{\zz}(\Xx, \es)} {\PP}_{\Xx\sm B_s(\zz)} ) + \P(\xx\lrsa_{\G_{\xx}(\Xx,\es)} \Pd_{\Xx\sm B_s(\xx)} )  + \pi(\xx,\zz) ) \la(\d \zz) \la(\d \xx )\bigg)^{1/2}.
\end{align*}
To get the first bound, we used that in order for $\zz$ to be connected to $\{\xx\}\cup \tilde{\PP}_{\Xx\sm (Q_n^+ \sm B_s(\xx))}$ it must either be connected to $\xx$ or to a point in ${\PP}_{\Xx\sm (Q_n^+ \sm B_s(\xx))}\cup \tilde{\PP}_{\Xx\sm (Q_n^+ \sm B_s(\xx))}$ by a path in $Q_n^+ \sm B_s(\xx)$. The second bound used that since the distance between $\xx$ and $\zz$ is at least $2s$, either there is a direct edge between them, or at least one of them must be connected to a point of $\Pd$ at distance at least $s$ away. Using \eqref{cond:sharp_vv} and \eqref{cond:sharp_v} in the case of Model (P) and \eqref{cond:sharp} and \eqref{cond:exponential_sim} in case of Model (L), we deduce that the whole term is $O(\la(Q_n)^{-m})$. 
\medskip

\ni{\bf Fifth term.}
Since the distance between $\zz$ and $\yy$ is at least $s$, we may bound $\pi(\zz,\yy)$ using \eqref{cond:sharp_v} for Model (P) and \eqref{cond:exponential_sim} for Model (L).
This means that the integral is of order $O\big(\la(Q_n)^{-2m-2}\big)$ and hence the $\La$-integral of this is of order
$\La(Q_n)\la(Q_n)^{-2m-2} \in O(\la(Q_n)^{-2m-1}).$
\enp

%
%
\section{Poisson approximation}\label{sec:papx}
In this section, we prove the Poisson approximation property stated in Theorem \ref{thm:1} using the theory of Palm couplings set up in Section \ref{sec:pac}. 
We again focus on Model (P) as it is the more challenging one.

%
%
\bep[Proof of Theorem \ref{thm:1}]
 Take $c_{\ms{St}} :=14\max(c_{\ms{es},2}^{-1},c_{\ms{v},2}^{-1},c_{\ms{D},2}^{-1})$ such that Lemma \ref{lem:1pei} holds for with $m\le 5$.

Let $\X$ and $\X_{\xx}^\Xi$, $\xx\in Q_n$, be constructed as in Section \ref{sec:pac}. Let $\mc M$ be a Poisson process with intensity measure $\La$. For $\xx \in Q_n\su \R^d\times \M$, we let $S_{n,\xx}:=B_{3 s}(\xx)$. We bound the symmetric difference $(\Xi \Delta \Xi_{\xx}^!)(Q_n)$ in the integral on the right-hand side of \eqref{eq:dkr} by the following three terms:
\begin{align}\label{symmdiffbou}
	\int_{Q_n} \one\{S_{n,\xx} \cap S_{n,\yy} \ne \es\} \Xi(\d \yy)+\hspace{-.13cm}\int_{Q_n}\hspace{-.15cm} \one\{S_{n,\xx} \cap S_{n,\yy} \ne \es\} \Xi_{\xx}^!(\d \yy)+\hspace{-.1cm}\int_{Q_n}\hspace{-.1cm} \one\{S_{n,\xx} \cap S_{n,\yy} = \es\} (\Xi \De \Xi_{\xx}^!)(\d y).
\end{align}
We now analyze the integral with respect to $\La$ of the expectation of each of the three terms separately.
\medskip

\noindent{\bf First term.}
For the first term, we obtain
$$
\E \Big[\int_{Q_n}\int_{Q_n} \one \{S_{n,\xx} \cap S_{n,\yy} \ne \es\} g(\yy,\XX) \XX(\d \yy)\La(\d \xx) \Big] = \int_{Q_n}\int_{Q_n}\hspace{-.1cm} \one \{S_{n,\xx} \cap S_{n,\yy} \ne \es\} \La(\d \xx )\La(\d \yy).
$$

\noindent{\bf Second term.} Recall that $\XX_{\xx}^{!,\Xi}:=\XX_{\xx}^\Xi\sm \{\xx\}$, where $\XX_{\xx}^\Xi$ is a Palm kernel of $\XX$ with respect to $\Xi$. The $\La$-integral of the expectation of the second term in \eqref{symmdiffbou} is by the definition of $ \XX_{\xx}^\Xi$ \eqref{def:palm2} and by the GNZ equation \eqref{cond:GNZ} given by
\begin{align*}
&\E \Big[ \int_{Q_n}\int_{Q_n} \one \{S_{n,\xx} \cap S_{n,\yy} \ne \es\} g(\yy,\XX_{\xx}^{!,\Xi}) \XX_{\xx}^{!,\Xi}(\d \yy)\La(\d \xx) \Big]\\
	&\quad = \E \Big[ \int_{Q_n} \int_{Q_n} \one \{S_{n,\xx} \cap S_{n,\yy} \ne \es\} g(\xx,\XX) g(\yy,\XX\sm \{\xx\}) (\XX\sm \{\xx\})(\d \yy) \XX(\d \xx) \Big]\\
	&\quad =\int_{Q_n} \int_{Q_n} \one \{S_{n,\xx} \cap S_{n,\yy} \ne \es\} \E [g(\xx,\XX\cup\{\xx, \yy\}) g(\yy,\XX\cup\{\yy\}) \k(\xx,\XX\cup\{\yy\}) \k(\yy,\XX)] \la (\d \xx) \la(\d \yy).
\end{align*}

\noindent{\bf Third term.} We write the third term on the right-hand side of \eqref{symmdiffbou} as
$$
\int_{Q_n} \one\{S_{n,\xx} \cap S_{n,\yy} = \es\} (\Xi[\XX] \De (\Xi[\XX_{\xx}^\Xi] \setminus \{\xx\}))(\d \yy).
$$
Now, note that a point $\yy\in \Xi[\XX] \De (\Xi[\XX_{\xx}^\Xi] \setminus \{\xx\})$ with $S_{n,\xx} \cap S_{n,\yy} = \es$ can only contribute to the symmetric difference if $\yy\in \X$ and $R(\yy,\XX)>s$, if $\yy\in \XX_{\xx}^\Xi$ and $R(\yy,\XX_{\xx}^\Xi)>s$, or if $\max(R(\yy,\XX), R(\yy,\XX_{\xx}^\Xi)) \le s$ (with $R(\xx,\vp)=0$ if $\xx\notin \vp$) and $(\XX \De \XX_{\xx}^\Xi) \cap Q_n \not \su B_{4s}(\xx)$. Therefore, the third term in \eqref{symmdiffbou} is bounded by
%
%
\begin{align}\nonumber
	&\int_{Q_n} \one\{R(\yy,\XX)>s\} \XX(\d \yy) + \int_{Q_n} \one\{R(\yy,\XX_{\xx}^\Xi)>s\} \XX_{\xx}^{!,\Xi}(\d \yy)\\
	& + (\XX \cup \XX_{\xx}^\Xi)(Q_n\sm S_{n,\xx}) \one\{(\XX \De \XX_{\xx}^\Xi) \cap Q_n \not\su B_{4s}(\xx) \}.\label{papx:symmdiff}
\end{align}
The expectation of the first term in \eqref{papx:symmdiff} is given by
\begin{equation}\label{eq:Rint}
\E\Big[\int_{Q_n} \one\{R(\yy,\XX)>s\} \XX(\d \yy)\Big]=\int_{Q_n} \E[\one\{R(\yy,\XX\cup \{\yy\})>s\} \k(\yy,\XX)] \la(\d \yy).
\end{equation}
Hence, it follows from \eqref{cond:PIbound}, \eqref{eq:s2} and the fact that $s \ge 3c_{\ms{es},2}^{-1} \log (\la(Q_n))$ that the $\La$-integral of \eqref{eq:Rint} is $O(\la(Q_n)^{-1})$. Similarly, the $\La$-integral of the expectation of the second term in \eqref{papx:symmdiff} is given by
\begin{align*}
&\int_{Q_n}\E\Big[\int_{Q_n} \one\{R(\yy,\XX_{\xx}^\Xi)>s\} \XX_{\xx}^{!,\Xi}(\d \yy)\Big] \La(\d \xx)\\
&\quad =\int_{Q_n} \int_{Q_n} \E[\one\{R(\yy,\XX\cup \{\xx,\yy\})>s\} g_n(\xx,\XX\cup\{\xx, \yy\}) \k(\xx,\XX\cup\{\yy\}) \k(\yy,\XX)] \la(\d \yy)\la(\d \xx),
\end{align*}
which is bounded using Cauchy-Schwarz and a similar argument as \eqref{eq:Rint} using $s \ge 6c_{\ms{es},2}^{-1} \log (\la(Q_n))$. For the $\La$-integral of the expectation of the third term in \eqref{papx:symmdiff} we have by Cauchy-Schwarz and the fact that $\E[(\XX \cup \XX_{\xx}^{!, \Xi})(Q_n\sm S_{n,\xx})^2] \le 2\E[\PP(Q_n)^2] = 2\la(Q_n)^2+2\la(Q_n)$ by Poisson domination (since $\k_\xx'\le 1$), 
\begin{align}\label{papx:cauchy}
&\int_{Q_n} \E[(\XX \cup \XX_{\xx}^{!,\Xi})(Q_n\sm S_{n,\xx}) \one\{(\XX \De \XX_{\xx}^\Xi) \cap Q_n \not\su B_{4s}(\xx)\}] \La(\d \xx) \nonumber\\
&\le \Big(\int_{Q_n} \E[(\XX \cup \XX_{\xx}^{!,\Xi})(Q_n\sm S_{n,\xx})^2] \La(\d \xx)\Big)^{1/2} \Big(\int_{Q_n} \P\big((\XX \De \XX_{\xx}^\Xi) \cap Q_n \not \su B_{4s}(\xx)\big) \La (\d \xx)\Big)^{1/2}\nonumber\\
&\le\Big( 2\la(Q_n) ( \la(Q_n)^2+ \la(Q_n))\int_{Q_n} \P\big((\XX \De \XX_{\xx}^\Xi) \cap Q_n\not\su B_{4s}(\xx) \big)\La(\d \xx) \Big)^{1/2}.
\end{align}
By Lemma \ref{lem:1pei} with $m=5$, we have
\begin{align}\label{papx:symdiffX}
\int_{Q_n}\P\big((\XX \De \XX_{\xx}^\Xi) \cap Q_n \not\su B_{4s}(\xx) \big)\La(\d \xx) \in O(\la(Q_n)^{-5}).
\end{align}
Hence, we conclude that \eqref{papx:cauchy} is in $O(\la(Q_n)^{-1})$.
\enp

	\def\Ep{E_{\ms{perc}}}
\def\Tcpp{T^{'',\mathsf{RCM},\mathsf{dc}}}
\def\Tcppp{T^{''',\mathsf{RCM},\mathsf{dc}}}

%
%
\section{Quantitative normal approximation}
\label{sec:malliavin_stein}
In this section, we prove Theorem \ref{thm:2}, i.e., the quantitative normal approximation of the considered geometric functionals.  Similarly as in \cite[Section 6]{hos}, this will be achieved by applying the normal approximation method  from \cite{chen} using Palm couplings. However, one important difference to \cite{hos} is that now we need to work with possibly unbounded interaction ranges.	This requires us to invest substantial efforts into proving the novel Lemma \ref{lem:prprpr}, which is a disagreement result for a thinning based on a Poisson point process that is independently resampled in some region.
Besides this novel ingredient, we briefly recall the main ideas of the Palm coupling from \cite{chen} to make our presentation self-contained. As in Sections \ref{sec:pac} and \ref{sec:papx}, we concentrate  on Model (P), since Model (L) is very similar. 

We again work in  the setting of Section \ref{sec:pac}. That is, we consider a coupling of the random measure $\Xi$ from  \eqref{eqn:poixidefXi} and its Palm kernel $\{\Xi_{\xx}\}_{\xx \in Q}$ defined in  \eqref{def:palm1}.
To state the main result, we put
\begin{align}
\label{def:Ws}
	Y_{\xx} :=\Xi_{\xx}(Q)-\Xi(Q), \quad  \De_{\xx} :=\f{Y_{\xx}}\s=\f{\Xi_{\xx}(Q)-\Xi(Q)}\s, 
\end{align}
where $\s^2:=\Var(\Xi(Q))$ denotes the variance of the total mass of the random measure $\Xi$. To make our presentation self-contained, we  recall \cite[Theorem 3.1]{chen}, which will be the key tool in the proof.
\bet [Theorem 3.1 in \cite{chen}] \label{thm:rollin}
Let $\Xi$, $\{\Xi_{\xx}\}_{\xx \in Q}$, $	Y_{\xx} $, and $\De_{\xx}$ be as in \eqref{def:Ws}. Then,
\begin{align*}
	\dk\Big(\frac{\Xi(Q)-\E[\Xi(Q)]}{\s}, N(0,1)\Big) \le 2E_1+5.5E_2+5E_3+10E_4+7E_5.
\end{align*}
Here
\begin{align*}
	E_1&:=\f1{\s^2} \E\Big|\int_Q (Y_{\xx}\, \one \{|Y_{\xx}|\le \s\}-\E\big[Y_{\xx} \, \one \{|Y_{\xx}|\le \s\}\big] )\,\La(\d \xx)\Big|,\\
	E_2&:=\f1{\s^3}\int_Q \E \big[Y_{\xx}^2\,\one \{|Y_{\xx}| \le \s\}\big] \,\La(\d \xx),\\
	E_3&:=\f1{\s^2}\int_Q \E \big[|Y_{\xx}|\,\one\{|Y_{\xx}| > \s\}\big] \,\La(\d \xx),\\
	E_4&:=\f1{\s^2} \int_{-1}^1 \int_Q \int_Q \Cov\big(\phi_{\xx}(t), \phi_{\yy}(t)\big) \, \La(\d \xx)\,\La(\d \yy)\,\d t,\\
	E_5&:=\f1{\s}\Big(\int_{-1}^1 \int_Q \int_Q |t|\,\Cov\big(\phi_{\xx}(t), \phi_{\yy}(t)\big) \, \La(\d \xx)\,\La(\d \yy)\,\d t\Big)^{1/2},
\end{align*}
where to simplify the notation, we write
\begin{align*}
	\phi_{\xx}(t)=\begin{cases}
		\one\{1 \ge \De_{\xx}>t>0\},\quad &t>0,\\
		\one\{-1 \le \De_{\xx}<t<0\},\quad &t<0.
	\end{cases}
\end{align*}
\ent

\ber
A subtle difference between the Palm couplings needed for  \eqref{eq:dkr} and Theorem \ref{thm:rollin} is that in  \eqref{eq:dkr}, it is enough to have pairwise couplings of $\Xi$ and $\Xi_{\xx}$ for each $\xx$, while for Theorem \ref{thm:rollin} it is necessary  to couple $\Xi$  and $\Xi_{\xx}$ for all $\xx\in Q$ simultaneously. The couplings \eqref{eqn:couplingX} and \eqref{eqn:coupling} are targeted to both applications, however, the need for simultaneous couplings in Theorem \ref{thm:rollin} substantially complicates the disagreement coupling set-up.
\enr

As in Section \ref{sec:papx}, we always fix the scales $r :=6s := 6c_{\ms{St}}\log(\la(Q_n))$ where $c_{\ms{St}} >0$ is a constant that will be chosen according to  Lemma \ref{lem:1pei} as well as Lemma  \ref{lem:pbe1}, \ref{lem:Yintbou}, \ref{lem:thm2cov} and \ref{lem:prprpr} below.

%
%

\subsection{2-ball  coupling.}
In Section \ref{sec:papx}, we  provided the couplings between $\X$ and $\X_{\xx}^\Xi$, $\xx\in Q_n^+$, given in \eqref{eqn:couplingX} and \eqref{eqn:coupling}. In order to apply Theorem \ref{thm:rollin}, we need to compare not only the disagreement between $\X_{\xx}^\Xi$ with $\X$ but also between  $\X_{\xx}^\Xi$ and $\X_{\yy}^\Xi$ for $\yy \ne \xx$.
More precisely, we need to compare $\XX_{\xx}^\Xi$ and $\XX_{\yy}^\Xi$ for $\xx, \yy \in Q_n$ satisfying $\|x - y\|>  4r $.  Since this is challenging to achieve directly, we define an intermediate process, which can be compared to both $\XX_{\xx}^\Xi$ and $\XX_{\yy}^\Xi$. 
We first define the PI $\k''_{\xx,\yy}$ associated with the intermediate process. More precisely, 
\begin{align}
	\label{eq:k''}
	\k''_{\xx,\yy}(\zz,\vp):=\begin{cases}
		\k\big(\zz,(\vp \cap ( B_{r + 3s}(\xx)\sm B_s(\xx)) ) \cup  \Psi_{\xx}\big), & \zz \in B_{r + 3s}(\xx)\sm B_s(\xx),\\
		\k\big(\zz,(\vp \cap ( B_{r + 3s}(\yy)\sm B_s(\xx))) \cup  \Psi_{\yy}\big), & \zz \in B_{r + 3s}(\yy)\sm B_s(\yy),\\
			0, &\text{otherwise}.
	\end{cases} 
\end{align}
Then $\k_{\xx,\yy}''$ has the form of $\k''$ in Lemma \ref{lem:pa1}. Thus, we can construct it by thinning a dominating RCM corresponding to the PI $\k$. 
Denoting this thinning by $\bar{T}^{\prime\prime,\xx,\yy,\ms{RCM,cl}}$, we define
\begin{align*}
\XX''_{\xx,\yy} := \bar{T}^{\prime\prime,\xx,\yy,\ms{RCM,cl}}_{Q_n^+, \es,\es}(\Phi_{Q_n^+,\es}^*) \cup (\Psi_{\xx}\cap B_s(\xx)) \cup (\Psi_{\yy}\cap B_s(\yy)).
\end{align*}
Note that $\k'_{\xx}, \k_{\xx,\yy}''$ satisfy Lemma \ref{lem:pa1}  with $P_1=Q_n^+$, $P_{1,-}=B_{r+3s}(\xx)$, $P_3=B_s(\xx)$, $P_{2,-}=B_{r + 3s}(\yy)\sm B_s(\yy)$, $\psi_1=\Psi_{\xx}$, $\psi_2=\Psi_{\yy}$ and $\psi=\Psi_{\xx}$, see also Figure \ref{fig:pert_RCM}. Thus, by Theorem \ref{thm:Talpha} with  $A=B_r(\xx)$ and $P=B_{r+3s}(\xx)$,
\begin{align}\nonumber
	&\P\big(\X_{\xx}^{!,\Xi}\cap B_r(\xx) \ne \X''_{\xx,\yy}\cap B_r(\xx)\big)\\ \nonumber
	&\le \P\big(R(\xx,\X_{\xx}^\Xi)>s\big)+\P\Big(\bar{T}^{\prime,\xx,\ms{RCM},\ms{cl}}_{Q_n^+,\es,\es}\big(\Phi_{Q_n^+,\es}^{*}\big)\cap B_r(\xx)  \ne \bar{T}^{\prime\prime,\xx,\yy,\ms{RCM},\ms{cl}}_{Q_n^+,\es,\es}\big(\Phi_{Q_n^+,\es}^{*}\big)\cap B_r(\xx) \Big)\\ \nonumber
	& \le  \P\big(R(\xx,\X_{\xx}^\Xi)>s\big)+ \P(\Ed)+\P\big(\EpRCM(B_r(\xx),s)\big)+4 \int_{B_{r+2s}(\xx)}\P\big(\zz\lrsa_{\G_{\zz}(Q_n^+,\es)} \Pd_{Q_n^+\sm B_{r+3s}(\xx)}  \big) \la(\d \zz) \\ \label{eq:X''disagree}
	&\phantom{\le} + \int_{Q_n^+ \sm B_{r+2s}(\xx) }\int_{B_{r+s}(\xx)\cap Q_n^+} \pi(\zz,\zz') \la(\d \zz) \la(\d \zz') .
\end{align}
For Model (L), the only necessary change is to invoke Theorem \ref{thm:disagree} instead of Theorem \ref{thm:Talpha}.
The key task is to show that when integrated over $Q_n^2$, this disagreement expression is negligible as $n \tff$.  

%
%
\bel[Bound on 2-point error integrals]
\label{lem:pbe1}
Suppose the PI satisfies the conditions of Theorem \ref{thm:2}.
Let $m \ge 1$. If  $c_{\ms{St}} > 0$ is chosen sufficiently large (depending on $m$), then
$$\iint_{Q_n^2\cap \{\|x-y\|>4r\}} \P\big(\X_{\xx}^{!,\Xi}\cap B_r(\xx) \ne \X''_{\xx,\yy}\cap B_r(\xx)\big) \La^2(\d \xx, \d \yy)\in O(\la(Q_n)^{-m}).$$
\enl
\bep
The proof of Lemma \ref{lem:pbe1} is very similar to the proof of Lemma \ref{lem:1pei}. Therefore, we omit it.
\enp

\subsection{Two auxiliary lemmas}
Having specified the couplings, we now proceed with the application of the general Theorem \ref{thm:rollin} to our setting. While bounding the expressions in Theorem \ref{thm:rollin} for general Palm measures may be cumbersome, for score sums of Gibbs processes, the task simplifies substantially. In earlier work, we reduced the problem to two key auxiliary results, namely \cite[Lemmas 6.3, 6.4]{hos}. We restate them below to make the manuscript self-contained. The first of the results is a moment bound on $Y_{\xx}$.

%
%
\bel[Moment bound for $Y_{\xx}$]
\label{lem:Yintbou}
Assume that the PI $\k$ satisfies the conditions of Theorem \ref{thm:2}. For $m \le 4$ and $n \ge 1$, if $c_{\ms{St}}$ is sufficiently large, we have 
$$
\int_{Q_n} \E\big[|Y_{\xx}|^m\big]\, \La (\d \xx) \le c_{\ms Y} \la(Q_n) s^{dm},
$$
where $c_{\ms Y}=c_{\ms Y}(d,\a,c_{\ms{es},1}, c_{\ms{es},2},c_{\ms{D},1},c_{\ms{D},2}, c_{\ms{v},1},c_{\ms{v},2}, c_{\ms m})> 0$.
\enl 

%
%
The second of these results is a covariance bound. 
\bel[Covariance bound]
\label{lem:thm2cov}
Assume that the PI $\k$ satisfies the conditions of Theorem \ref{thm:2}. Let $F:\mathbf N \to [0,\ff)$ be measurable and $p>0$ and assume that there is some $q_0=q_0(d,\a,c_{\ms{es},1}, c_{\ms{es},2},c_{\ms{D},1},c_{\ms{D},2}, c_{\ms{v},1},c_{\ms{v},2}, c_{\ms m})>0$
such that $F(\vp)\le  \la(Q_n)^p \vp(\Xx)$  for $\vp \in \Nlf$ and $\la(Q_n)>q_0$. Then, for $c_{\ms{St}}$ sufficiently large and $\la(Q_n)>q_0$,  
$$
\iint_{Q_n^2 \cap \{\|x-y\|>4r\}} \Cov\big(F_\xx, F_\yy)\big)\, \La^2(\d \xx, \d \yy) \le  \la(Q_n)^{-1},
$$
where $F_{\zz}=F(\XX_{\zz}^\Xi \cap B_r(\zz))$. 
\enl
The analogue of Lemma \ref{lem:thm2cov} holds when one or more of $F_\xx$ and $F_\yy$ are replaced by $F(\XX_{\xx}^\Xi \cap B_r(\xx))$ or $F(\XX_{\yy}^\Xi \cap B_r(\yy))$, but we only give the proof in the  case stated in the lemma.
We now prove Lemma~\ref{lem:Yintbou} and \ref{lem:thm2cov}. The proofs rely on  disagreement coupling through Lemma \ref{lem:1pei} and~\ref{lem:pbe1}.

%
%
 %
 %
 \bep[Proof of Lemma \ref{lem:Yintbou}] 
 
 We decompose  $Y_{\xx}$ as in \cite{hos}. That is, set  
	\begin{align}
		\Xi_{\xx}^\le:= \sum_{\yy \in \XX_{\xx}^\Xi \cap Q_n}g(\yy,\XX_{\xx}^\Xi) \mathds 1\{R(\yy,\XX_{\xx}^\Xi)\le s\}\de_{\yy},\qquad \Xi^\le:= \sum_{\yy \in \XX \cap Q_n}g(\yy,\XX) \mathds 1\{R(\yy,\XX)\le s\}\de_{\yy}.\label{def:Xile}
	\end{align}
and let	$\Xi^>:=\Xi- \Xi^\le$ and $\Xi_{\xx}^>:=\Xi_{\xx} - \Xi_{\xx}^\le$.
 Then, we can decompose $Y_{\xx}$ as
	\begin{align*}
		\big(|\Xi_{\xx}^> |-|\Xi^>|\big) + \big(\Xi_{\xx}^\le(B_r(\xx))-\Xi^\le (B_r(\xx)) \big) +\big(\Xi_{\xx}^\le ( B_r(\xx)^c)-\Xi^\le (  B_r(\xx)^c)\big) \one \{\XX \Delta \XX_{\xx}^\Xi \not \su B_{4s}(\xx)\}.
	\end{align*}
The main idea of the proof is then to apply Lemma \ref{lem:1pei} to bound $\P(\XX \Delta \XX_{\xx}^\Xi \not \su B_{4s}(\xx))$.

From the Hölder inequality applied to both the expectation and the integral with respect to $\La$, we get
	\begin{align*}
		\int_{Q_n} \E[ | Y_{\xx}|^m\big] \, \La(\d \xx) \le \hspace{-0.4cm} \sum_{\substack{{I=(i_1,\dots,i_6)}, \\ i_1+\cdots+i_6=m}} d_I&\Bigg[\Big(\int_{Q_n} \E\big[\Xi_{\xx}^\le\big(B_r(\xx)\big)^m\big] \, \La(\d \xx)\Big)^{\f{i_1} m}  \Big(\int_{Q_n} \E\big[\Xi^\le\big(B_r(\xx)\big)^m\big] \, \La(\d \xx)\Big)^{\f{i_2} m}\nonumber\\
		&\times \Big(\int_{Q_n} \E\big[|\Xi_{\xx}^>|^m\big] \, \La(\d \xx)\Big)^{\f{i_3} m}  \Big(\int_{Q_n} \E\big[|\Xi^>|^m\big] \, \La(\d \xx)\Big)^{\f{i_4} m}\nonumber\\
		&\times \Big(\int_{Q_n} \E\big[|\Xi_{\xx}|^m  \one \{\XX \De\XX_{\xx}^\Xi \not\su B_{4s}(\xx)\}\big] \, \La(\d \xx)\Big)^{\f{i_5} m}\nonumber\\
		&\times \Big(\int_{Q_n} \E\big[|\Xi|^m  \one \{\XX \De\XX_{\xx}^\Xi \not\su B_{4s}(\xx)\}\big] \, \La(\d \xx)\Big)^{\f{i_6} m} \Bigg]
	\end{align*}
	for some $d_I > 0$. As explained in \cite[Lemma 6.3]{hos}, the most difficult task is the fifth integral. First, by the H\"older inequality, this integral is bounded above by 
$$\int_{Q_n} \E \big[|\Xi_{\xx}|^m  \one \{\XX \De\XX_{\xx}^\Xi \not\su B_{4s}(\xx)\}\big] \, \La(\d \xx) \le \Big(\int_{Q_n}\hspace{-.3cm} \E\big [|\Xi_{\xx}|^{m + 1}  \big] \, \La(\d \xx)\Big)^{\f m{m + 1}} \Big(\int_{Q_n}\hspace{-.3cm} \P\big(\XX \De\XX_{\xx}^\Xi \not\su B_{4s}(\xx)\big)\, \La(\d \xx)\Big)^{\f1{m + 1}}.$$
Now,  we can apply Lemma \ref{lem:1pei}  to control the integral
$\int_{Q_n} \P\big(\XX \De\XX_{\xx}^\Xi \not\su B_{4s}(\xx)\big)\, \La(\d \xx)$, thereby concluding the proof.
\enp

%
%
We next sketch the proof of Lemma \ref{lem:thm2cov} and explain the main differences from \cite{hos}.

%
%
\bep[Proof of Lemma \ref{lem:thm2cov}]
 Assume that $\|x - y\| > 4r$. For $\zz \in \{\xx,\yy\}$, we set $R_{\zz} := R(\zz,\XX_{\zz}^\Xi)$.
Then,
        \begin{align*}
        \Cov\big(F_{\xx}, F_{\yy}\big)
        &=\Cov\big(\tilde F_{\xx}, \tilde F_{\yy}\big) + \Cov\big(\tilde F_{\xx},  F_{\yy}\one\{R_{\yy}> s\}\big) +  \Cov\big(F_{\xx} \one\{R_{\xx}> s\}, F_{\yy} \big),
\end{align*}   
where 
      $\tilde F_{\zz}:= F\big(\XX_{\zz}^\Xi \cap B_r(\zz)\big) \one\big\{R_{\zz}\le s\big\}.$
We only discuss the first covariance since the second and third can be bounded  using \eqref{def:palm2} and the exponential stabilization \eqref{eq:s2} as in \cite[Lemma 6.4]{hos}. 

We introduce the agreement event $S_{\zz}:=\big\{(\XX_{\zz}^\Xi \De \XX_{\xx,\yy}'') \cap B_r(\zz)=\es \big\}$. Lemma \ref{lem:pbe1} provides bounds on the probability of the complementary event $\P(S_\xx^c)$.  This allows us to further reduce to the consideration of the covariance $\Cov\big(\tilde F_{\xx} \one_{S_{\xx}}, \tilde F_{\yy}\one_{S_{\yy}}\big)$. However, this covariance vanishes since $\Psi_{\xx}$ and $\Psi_{\yy}$ are independent by construction and given these, $\XX''_{\xx,\yy}\cap B_r(\xx)$ and $\XX''_{\xx,\yy}\cap B_r(\xx)$ are independent. To see this, note that one can construct $\XX_{\xx,\yy}''$  by a Poisson thinning that first considers $B_{r+3s}(\xx)$ and then $B_{r+3s}(\yy)$. Since the PI on $B_{r+3s}(\yy)$ does not depend on the result of thinning $B_{r+3s}(\xx)$, the distributions of $\XX''_{\xx,\yy}\cap B_{r+3s}(\xx)$ and $\XX''_{\xx,\yy}\cap B_{r+3s}(\yy)$ are independent.
\enp

\subsection{Proof of Theorem \ref{thm:2}}
In this section, $C>0$ denotes a constant which may change value from line to line. It may depend on $d,\a,c_{\ms{es},1}, c_{\ms{es},2},c_{\ms{D},1},c_{\ms{D},2}, c_{\ms{v},1},c_{\ms{v},2}, c_{\ms m}$, but not on $n$.


\bep[Proof of Theorem \ref{thm:2}] Given Lemma \ref{lem:Yintbou} and \ref{lem:thm2cov}, the proof of Theorem \ref{thm:2} follows exactly as in \cite{hos}. There is one error in that proof, which we now  correct. When bounding $E_4$ from Theorem \ref{thm:rollin}, the computation was reduced  to showing that for a given $l\ge 1$, $c_{\ms{St}}$ can be chosen so large that 
$$\int_\M\int_{\M}\Cov(\phi_{(x,m)}(t),\phi_{(y,m')}(t))\mu(\d m) \mu(\d m')\le C\la(Q_n)^{-l}$$
whenever $\|x-y\|>4r$. In the following, let $\xx = (x,m)$ and $\yy=(y,m')$.
Similar to \cite{hos}, we let 
$$U_\xx := \big\{\XX \De \XX_\xx^\Xi \su B_{4s}(\xx)\} \cap \big\{ \max_{\ww \in\XX \cap Q_n} R(\ww,\XX)\le s\big\} \cap \big\{\max_{\ww \in\XX_\xx^\Xi \cap Q_n}R(\ww,\XX_\xx^\Xi)\le s\big\}.$$ 
Then, since $0\le \phi_{\xx}(t)\le 1$,
$$|\Cov(\phi_{\xx}(t),\phi_{\yy}(t))-\Cov(\phi_\xx(t)\one\{U_\xx\},\phi_\yy(t)\one\{U_\yy\})| \le \P(U_\xx^c) + \P(U_\yy^c).$$
Under $U_\xx$,
\begin{align}\label{eq:phidef}
	\phi_{\xx}(t)=\phi_{\xx,r}(t):=\begin{cases}
		\one\{1 \ge \frac{\Xi [\XX_\xx^\Xi \cap B_r(\xx)] (B_{5s}(\xx)) - \Xi[\XX \cap B_{r}(\xx)](B_{5s}(\xx))}{\sigma} >t>0\},\quad &t>0,\\
		\one\{-1 \le \frac{\Xi [\XX_\xx^\Xi \cap B_{r}(\xx)] (B_{5s}(\xx)) - \Xi[\XX \cap B_{r}(\xx)](B_{5s}(\xx))}{\sigma}  <t<0\},\quad &t<0.
	\end{cases}
\end{align}
Thus,
$$|\Cov(\phi_{\xx,r}(t),\phi_{\yy,r}(t))-\Cov(\phi_\xx(t)\one\{U_\xx\},\phi_\yy(t)\one\{U_\yy\})| \le \P(U_\xx^c) + \P(U_\yy^c).$$
The integral of $\P(U_\xx^c) + \P(U_\yy^c)$ can be bounded using Lemma \ref{lem:1pei} and \eqref{eq:s2}. Hence, it is enough to bound the integral over the marks of $\Cov(\phi_{\xx,r}(t),\phi_{\yy,r}(t))$.
However, $\phi_{\xx,r}(t)$ is not of the form $F(\XX_{\xx}^\Xi\cap B_r(\xx))$, as claimed in \cite{hos}, since it  depends on the joint distribution of $\XX\cap B_{r}(\xx)$ and $\XX_{\xx}^\Xi \cap B_{r}(\xx)$.  Hence, Lemma \ref{lem:thm2cov} does not apply.

Instead, we introduce four interpolating PIs of the form 
\begin{align*}	&	\k'''_{\zz}({\bf w},\vp) = \k\big({\bf w},\vp \cap B_{r + 3s}(\zz)   \big)\one\{{\bf w}\in B_{r+3s}(\zz) \}\\
&\k^{\prime \prime \prime,\Xi}_{\zz}({\bf w},\vp) = \k\big({\bf w},(\vp \cap (B_{r + 3s}(\zz) \sm B_s(\zz )))\cup \Psi_{\zz} \big)\one\{ {\bf w}\in B_{r+3s}(\zz) \sm B_s(\zz )\},
\end{align*}	
where $\zz \in \{\xx,\yy\}$. Denote the corresponding processes by $\XX_{\zz}'''$ and $\XX_{\zz}^{\prime\prime\prime,\Xi}$.  
With high probability, $\X \cap B_{r}(\zz) = \X_{\zz}^{\prime\prime\prime} \cap B_{r}(\zz) $ and $\X^\Xi_{\zz} \cap B_{r}(\zz) = \X_{\zz}^{\prime\prime\prime,\Xi} \cap B_{r}(\zz) $ by the same argument as in Lemma \ref{lem:1pei} and \ref{lem:pbe1}. It follows that
   $$\big|\Cov({\phi}_{\xx,r}(t),{\phi}_{\yy,r}(t))- \Cov({\phi}_{\xx,r}'''(t),{\phi}_{\yy,r}'''(t))\big|\le C \la(Q_n)^{-l},$$
 where ${\phi}_{\zz,r}'''(t)$ are computed by inserting the processes $\X_{\zz}^{\prime\prime\prime}$ and $\X_{\zz}^{\prime\prime\prime,\Xi}$ in \eqref{eq:phidef}. 
 
 It remains to bound $\Cov({\phi}_{\xx,r}'''(t),{\phi}_{\yy,r}'''(t))$. While the processes $\XX_{\xx}'''$ and $ \XX_{\xx}^{\prime\prime\prime,\Xi}$ are pairwise independent of $\XX_{\yy}'''$ and $ \XX_{\yy}^{\prime\prime\prime,\Xi}$, the coupling $(\XX_{\xx}''', \XX_{\xx}^{\prime\prime\prime,\Xi})$ may be correlated with the coupling $(\XX_{\yy}''', \XX_{\yy}^{\prime\prime\prime,\Xi})$ because the underlying Poisson process is explored in the same order. 
 Loosely speaking, we will instead show that the joint distribution of $(\XX_{\xx}''', \XX_{\xx}^{\prime\prime\prime,\Xi})$ is close to independent of the joint distribution $(\XX_{\yy}''', \XX_{\yy}^{\prime\prime\prime,\Xi})$.

 To formalize this for Model (L), let  $\Check{\Check{\PP}}^*$ be an independent copy of $\Pds$. Define a new Poisson process ${\check{\PP}}^* = \PP_{ B_{2r}(\xx)}^* \cup \Check{\Check{\PP}}^*_{Q_n^+ \sm B_{2r}(\xx)} $. We write $\check{\XX}_{\zz}'''$, $\Check{\Check{\XX}}_{\zz}'''$ etc.\ depending on which Poisson process was used in the thinning and $\check{\phi}_{\zz,r}'''(t)$ or $\Check{\Check{\phi}}_{\zz,r}'''(t)$ when $\check{\XX}_{\zz}'''$ and $\check{\XX}_{\zz}^{\prime\prime\prime,\Xi}$ or $\Check{\Check{\XX}}_{\zz}'''$ and $\Check{\Check{\XX}}_{\zz}^{\prime\prime\prime,\Xi}$, respectively, was used to define \eqref{eq:phidef}.

  The joint distribution of the  four processes $\XX_{\zz}'''$ and $\XX_{\zz}^{\prime\prime\prime,\Xi}$, $\zz\in \{\xx,\yy\}$, will be the same as that of $\check{\XX}_{\zz}'''$ and $\check{\XX}_{\zz}^{\prime\prime\prime,\Xi}$. It follows that 
  $$\Cov({\phi}_{\xx,r}'''(t),{\phi}_{\yy,r}'''(t))= \Cov(\check{\phi}_{\xx,r}'''(t),\check{\phi}_{\yy,r}'''(t)).$$
On the other hand, Lemma \ref{lem:prprpr} below shows that 
$$\int_{\M}\int_{\M}| \Cov(\check{\phi}_{\xx,r}'''(t),\check{\phi}_{\yy,r}'''(t))- \Cov({\phi}_{\xx,r}'''(t),\Check{\Check{\phi}}_{\yy,r}'''(t))| \mu(\d m)\mu (\d m')\le 4 C\la(Q_n)^{-l}.$$
  But $\Cov({\phi}_{\xx,r}'''(t),\Check{\Check{\phi}}_{\yy,r}'''(t))=0$ since $\PP$ and $\Check{\Check{\PP}}$ are independent, which shows the claim for Model (L).
  
  For Model (P), we take independent copies $\Phi_{Q_n^+,\es}^*$ and $\Check{\Check{\Phi}}_{Q_n^+,\es}^*$. From the underlying Poisson processes $\PP$ and $\Check{\Check{\PP}}$, we construct $\check \PP$ as under Model (L). Given $\check \PP$, we add i.i.d.\ marks to form $\check{\Phi}_{Q_n^+,\es}^*$ in such a way that marks $r_{\zz,\ww}$ with $\zz,\ww\in B_{2r}(\xx)$ or $\zz,\ww\in Q_n^+ \sm B_{2r}(\xx)$ agree with those in  $\Phi_{Q_n^+,\es}^*$ or $\Check{\Check{\Phi}}_{Q_n^+,\es}^*$, respectively. From there, the argument proceeds the same.
  \enp
  
The following lemma was used in the proof of Theorem \ref{thm:2}. Since the proof of the lemma involves some technical disagreement coupling arguments, we postpone it to Section \ref{sec:proofthm2}.

   \bel \label{lem:prprpr}
 Let $l\ge 1$.  Under the conditions of Theorem \ref{thm:2}, there is a $C$ (depending on $l$) such that if $c_{\ms{St}}$ is chosen sufficiently large,
 \begin{align*}
\max\bigg\{{}&\int_{\M}\int_{\M}\P(\XX_{\xx}''' \ne \check \XX_{\xx}''') \mu(\d m) \mu(\d m'), \int_{\M}\int_{\M} \P(\XX_{\xx}^{\prime\prime\prime,\Xi}\ne \check \XX_{\xx}^{\prime\prime\prime,\Xi})\mu(\d m) \mu(\d m'),\\ 
&\int_{\M}\int_{\M}\P(\Check{\Check{\XX}}_{\yy}''' \ne \check \XX_{\yy}''') \mu(\d m) \mu(\d m'), \int_{\M}\int_{\M} \P(\Check{\Check{\XX}}_{\yy}^{\prime\prime\prime,\Xi}\ne \check \XX_{\yy}^{\prime\prime\prime,\Xi})\mu(\d m) \mu(\d m')\bigg\}\le C\la(Q_n)^{-l}.
\end{align*}
  \enl

	\section{Proofs for Section \ref{sec:const}}\label{sec:proofs_const}

\bep[Proof of Corollary \ref{cor:poiss_emb}]
We divide $A$ into disjoint sets $A_1,\ldots,A_N$ such that $\iota(\xx)<\iota(\yy)$ for $\xx\in A_i$, $\yy\in A_{i+1}$ and we assume $\la(A_i) ={\la(A)}/{N}$  for all $i$. Such a partition can be found by using $\iota$ to map $A$ into $\R$ and dividing $\R$ into intervals of equal push-forward measure. Let $T_i=T_{Q,\es}(\Pds_Q)\cap A_i$ and $T_i'=T_{Q,\es}'(\Pds_Q)\cap A_i$. Then, using Lemma \ref{lem:distv},
\begin{align*}
	&\P\big(T_{Q,\es}(\Pds_Q)\cap A\ne T'_{Q,\es}(\Pds_Q) \cap A\big) = \sum_{i=1}^N  \P( T_i\ne  T_i', T_j=T_j',j< i)\\
	&\le\sum_{i=1}^N  \P(\Pd(A_i)>1) + \sum_{i=1}^N  (2+\la(A_i))\\
	&\quad\times \E \Big[\one\{ T_j=T_j',j< i\}  \dtv\Big(\X\Big(Q\sm \bigcup_{j<i}A_i,  \bigcup_{j<i} T_j\Big)\cap A_i ,\X'\Big(Q\sm \bigcup_{j<i}A_i,  \bigcup_{j<i} T_j'\Big)\cap A_i \mid T_j,T_j',j<i\Big)\Big]\\
	&\le\sum_{i=1}^N  \P(\Pd(A_i)>1) + \sum_{i=1}^N  (2+\la(A_i) )\P(\Pd_{A_i}\lrsa_{\Pd \cup \bar\eta}   \Pd_{Q\sm P_1}\cup \bar\eta)\\
	&\le\sum_{i=1}^N  \P(\Pd(A_i)>1) + \sum_{i=1}^N  2\P(\PP_{A_i} \lrsa_{\Pd \cup \bar\eta}   \Pd_{Q\sm P_1}\cup \bar\eta ) + \sum_{i=1}^N\la(A_i)\P(\Pd(A_i)>0)
\end{align*}
Here, the total variation distance is in the conditional distribution given $T_j,T_j',j<i$. This was bounded in the second inequality by  constructing the two Gibbs processes on $Q\sm \cup_{j<i} A_j $ as follows: First construct both  processes on $Q\sm P $ by some thinning of $\Pds_{Q\sm  P}$. Then, we use a cluster-based thinning that starts by exploring all clusters connected to $ B= \Pd_{Q\sm P}\cup \bar \eta $. After all such clusters are explored, the PIs will look the same on $P\sm \cup_{j<i} A_j$ and the boundary conditions will agree since $T_j=T_j'$ for $j<i$. Thus, the result follows from Proposition \ref{pr:rad_prop}. 
The first and the last sum are of order $N^{-1}$ when $N\to\ff$, while the middle sum is bounded by
$$\sum_{i=1}^N  2\int_{A_i}\P(\xx\lrsa_{\Pd_{\xx} \cup \bar\eta}  \bar\eta\cup \Pd_{Q\sm P_1})\la(\d\xx) = 2\int_{A}\P(\xx\lrsa_{\Pd_{\xx} \cup \bar\eta}  \bar\eta\cup \Pd_{Q\sm P_1})\la(\d\xx),$$
as asserted.
\enp

Next, we move to the proof of Theorem \ref{thm:disagree}. Recall from Section \ref{ss:pert} that we divide $\R^d \times \M$ into (marked)  cubes of spatial diameter at most 1 and let $C_l=C_l^0\times \M$ denote  those cubes intersecting $Q$, $l=1,2,\ldots,L$.

\bep[Proof of Theorem \ref{thm:disagree}] 
Both $\X:=T^{\ms{cl}}_{Q,\es,\es}(\Pds_Q)$ and $\X':=T^{\prime,\ms{cl}}_{Q,\es,\es}(\Pds_Q)$ visit the points of $\Pd_Q$ in the same order. 
The strategy is to bound the error probabilities for each step in the algorithmic construction. For that purpose, let $V_i=S_i\sm S_{i-1}$ be the set explored in the $i$th step and $\X_i=\X\cap V_i$, $\X_i'=\X'\cap V_i$ be the parts of the Gibbs processes generated in the $i$th step.  Moreover, let $\GGc_i$ be the part of the graph $\GGc$ obtained in the $i$th step.

Note that under the event $\Ed^c$, there are at most $M:=4\la(Q)$ steps. Under the event $\Ep(A,s)$, if $\X$ and $\X'$ agree on all clusters completely contained in $B_s(A)$, then they agree on $A$. Hence, 
\begin{align*}
\P\big(\XX\cap A\ne \XX'\cap A, \Ep(A,s)^c\cap \Ed^c\big) \le \E\Big[ \sum_{i=1}^M \P(\X_i \ne  \X_i',  E_i^1\cup E_i^2\mid \Pds_{S_{i-1}})\Big],
\end{align*}
where $E_i^1$ denotes the event that   for all $j\le i-1$ with $\Pd _{V_{j}}\lrsa_{\GGc_{i-1}} \Pd_{V_{i-1}} $  it holds that $\Pd_{V_j}\su B_s(A) $ and  $\X_{j} = \X_{j}'$ (that is, we continue exploring a cluster so far contained in $B_s(A)$ and the two Gibbs processes so far agree on the cluster) and $E_i^2$ is the event that $\Pd_{V_{i-1}}=\es$ and $\Pd_{V_i}\su B_s(A)$ (so we start a new cluster from a point in $B_s(A)$). 

Under $E_i^1$, we have the bound from Corollary \ref{cor:poiss_emb} 
\begin{align*}
	&\P\big(\X_i \ne  \X_i', E_i^1\mid \Pds_{S_{i-1}}\big) \le 2\int_{V_i} \P(\xx\lrsa_{\Pd_\xx\cup \bar\eta} \Pd_{Q\sm P} \cup \bar{\eta})\la(\d\xx)\one\big\{\Pd_{V_{i-1}}\su B_{s}(A)\big\}\\
	&\le 2\int_{V_i\cap B_{2s}(A)} \P(\xx\lrsa_{\Pd_\xx\cup \bar\eta} \Pd_{Q\sm P} \cup \bar{\eta})\la(\d\xx)+ 2\la(N(\Pd_{V_{i-1}})\cap (Q \sm B_{2s}(A))) \one\big\{\Pd_{V_{i-1}}\su B_{s}(A)\big\}
\end{align*}
since $V_i \su N(V_{i-1})$.
Summing over $i$, and remembering that the $V_i$ are disjoint, we get
\begin{align*}
	&\sum_{i=1}^M \P(\X_i \ne  \X_i', E_i^1\mid \Pds_{S_{i-1}}) \le 2\int_{ B_{2s}(A)} \P(\xx\lrsa_{\Pd_\xx\cup \bar\eta} \Pd_{Q\sm P} \cup \bar{\eta})\la(\d\xx)+ 2\sum_{\xx\in \Pd_{B_s(A)}}\la\big(N(\xx)\cap (Q \sm B_{2s}(A))\big). 
\end{align*}
Taking expectations, we get
\begin{align}
	\E \Big[\sum_{i=1}^M \P(\X_i \ne  \X_i', E_i^1\mid \Pds_{S_{i-1}}) \Big]
	\hspace{-.1cm}\le \hspace{-.1cm}2\int_{ B_{2s}(A)}\hspace{-.8cm} \P(\xx\lrsa_{\Pd_\xx\cup \bar\eta} \Pd_{Q\sm P} \cup \bar{\eta})\la(\d\xx)\hspace{-.1cm}+\hspace{-.1cm} 2 \int_{B_s(A)}\int_{Q\sm B_{2s}(A)}\hspace{-.95cm}\one\{\xx\sim \yy\}\la (\d \yy)\la(\d\xx). \label{eq:u2}
\end{align}
Under $E_i^2$, we construct the Gibbs process on $Q\sm S_{i-1}$ by a standard Poisson embedding and take the thinning of the first Poisson point as $\X_i$ and $\X'_i$, respectively. By the choice of ordering, the Poisson embedding on $Q\sm S_{i-1}$ is first made inside $(Q\sm S_{i-1})\cap C_L$. Conditionally on the result, we make a Poisson embedding on $(Q\sm S_{i-1})\cap C_{L-1}$ and so on. Thus, Corollary \ref{cor:poiss_emb} gives
\begin{align}\nonumber
	&\P\big(\X_i \ne  \X_i', E_i^2 \mid \Pds_{S_{i-1}}\big) \\ \nonumber
	&\le \sum_{l=1}^L \one\big\{C_l \cap B_{s}(A) \ne \es\big\} \P\Big( \inf{}_\iota (\X_i\cap C_l) \ne \inf{}_\iota(\X_i'\cap C_l) \mid \Pd_{C_{l'}\cap(Q\sm S_{i-1})}=\es, l'>l,\Pds_{S_{i-1}}\Big)\\\nonumber
	&\le \sum_{l=1}^L 2\int_{B_{s+1}(A)\cap C_l}  \P\big(\xx\lrsa_{\Pd_\xx\cup \bar\eta} \Pd_{Q\sm P} \cup \bar{\eta}\big) \la(\d\xx) \\
	&\le 2\int_{B_{s+1}(A)}  \P\big(\xx\lrsa_{\Pd_\xx\cup \bar\eta} \Pd_{Q\sm P} \cup \bar{\eta}\big) \la(\d\xx) . \label{eq:u1}
\end{align}
Adding  \eqref{eq:u2}, \eqref{eq:u1}, and the probability that we are not under $\Ed^c\cup \Ep(A,s)^c$ yields the result.
\enp

\bep[Proof of Proposition \ref{pr:xxff}]
The cluster-based thinning algorithm $T^{\ms{cl}}_{U_n,\xi_n, \psi}$ will first move through all the clusters having points connected to $\xi_n$. These are contained in the clusters of $\GG(\Pd_U)$ intersecting $U\sm W_n$. Afterwards, it explores the remaining clusters in the order of their $\iota$-smallest point. Let $\xx\in \Pd$ be the $\iota$-smallest point connected to $\Pd_A$ and suppose  $\xx\in W_n$. Then, by the assumption on $\iota$, all clusters with $\iota$-smallest point larger than $\xx$ must also be contained in $W_n$. Thus, the order in which these clusters are visited will be the same for all $m\ge n$. When the thinning algorithm reaches $\xx$, the unexplored part of $U$ will be the points in  ${U_{(\xx,\infty)}}$ that are not neighbors to a cluster that contains a point in $\Pd_{U_{(-\infty,\xx)}}$. Again, this set does not depend on $m\ge n$. Hence, on the event $\big\{\Pd_A \not \lrsa_{\Pd_U} \Pd_{U \sm W_n}\big\}$, $T^{\ms{cl}}_{U_m , \xi_m,\psi}(\Pds_{U_m })$ will agree on all clusters connected to $A$ whenever $m\ge n$.

For the last statement, we need to show that $T^{\ff}_{U,\psi}(\Pds_U)$ satisfies the GNZ equation \eqref{cond:GNZ} for any $f:\Xx \times \Nlf \to [0,1]$ with support in $W_s\times \Nlf$ and such that $f(\xx,\vp) =f(\xx,\vp\cap W_s)$ for some $s>0$. Let $n>s$. Conditionally on $\xi_n$, $T^{\ms{cl}}_{U_n ,\xi_n, \psi }(\Pds_{U_n})$ is distributed as a $ \X(U_n ,\psi)$-process with PI $\k$, so this also holds unconditionally. Hence it satisfies the GNZ equation \eqref{cond:GNZ_boundary} for $f$. Taking limits and applying dominated convergence, we see that also $T^\ff_{U,\psi}(\Pds_U)$ satisfies the GNZ equation \eqref{cond:GNZ_boundary} for $f$ and $\XX(U,\psi)$ and hence is an infinite-volume Gibbs point process  with PI $\k$, which is unique by Theorem \ref{thm:unique1}.  
\enp

\bep[Proof of Proposition \ref{pr:chimera}]
In the following, we write $\X_U:=\X\cap U $, $B=\Xx\sm U$, and $\X_B:= \X\cap B$. From the GNZ equations applied to a function of the form $f\cdot h$ where $f:U\times \Nlf_U \to [0,1]$ and $h: \Nlf_B \to [0,1]$, we get
\begin{align*}
	&\int \E \Big[\int_U f(\xx,\X_U) \X(\d\xx)|\X_B=\psi\Big] h(\psi) \P_{\X_B}(\d \psi)\\
	&= \int \E\Big[\int_U f(\xx,\X_U\cup \{\xx\}) \k(\xx,\X_U\cup\psi) \lambda(\d\xx)|\X_B=\psi\Big]h(\psi)\P_{\X_B}(\d \psi) .
\end{align*}
Thus, for $\P_{\X_B}$-almost all $\psi$, $\X_U \mid \X_B =\psi$ satisfies the GNZ equation of $\X(U,\psi)$. Since the distribution of $\X(U,\psi) $ is unique by Theorem \ref{thm:unique1}, we deduce $\X_U \mid \X_B=\psi \sim T^{\ff}_{U,\psi}({\PP}^{*}_U)$ and hence the claim follows.
\enp

\section{Proofs for Section \ref{sec:RCM}}\label{sec:proofsRCM}

\subsection{Set-up}

The model with pairwise interaction does not fit the framework described in Section \ref{sec:const} directly because the RCM does not come from a symmetric relation on $\Xx$. Therefore, we start by discussing the construction from \cite{betsch} of an approximating Gibbs process on some enlarged space. There is a symmetric relation on this enlarged space such that \eqref{cond:PI_sim} holds and the associated graph approximates the RCM. 
This construction will be the main ingredient in the proofs. The approximation of the RCM is given in Section \ref{sss:approx} and the approximating Gibbs process is defined in Section \ref{sss:gibbsapprox}.

%
%
\subsubsection{Approximation of the RCM}
\label{sss:approx}

We first review the construction from \cite{betsch}. The idea is to approximate Gibbs processes with a non-negative pair potential on $Q$ by constructing an approximating Gibbs process on a marked space
$Q \times [0,1]^{\N}$
equipped with the product measure $\lambda\otimes \Q^\N$.
More precisely, for each $\delta>0$ we choose a partition of $\Xx $ into disjoint Borel sets $D_1^\delta, D_2^{\delta},\dots$ such that for  any $\xx,\yy\in\Xx$ there is a $\delta_0(\xx,\yy) $ such that $\xx$ and $\yy$ are separated by the partition for $\delta<\delta_0$. The sequence should be nested in the sense that for $\de'\le \de$ and $i\in \N$ there is a $j\in \N$ such that $D_i^{\de'}\su D_j^\de$. The index set $\N$ of $[0,1]^\N$ corresponds to the indices of the sets $D_i^\delta$. We write $\Uu^{\delta}=[0,1]^{\N}$ to emphasize this correspondence and $Q^\de = Q\times \Uu^\de$ for $Q\su \Xx$ etc. For each $\de>0$, we have the maps $h^\de:{\Xx^\de}\to \Xx$
and $h^{*,\de}:{(\Xx^*)^\de}\to \Xx$ given by 
$$h^\de(\xx, \rr) =\xx, \qquad h^{*,\de}(\xx, u, \rr) =\xx.$$
We assume throughout that $\vp^\de\in \Nlf_{\Xx^\de}$ denotes a configuration of points that satisfies $h^\de(\vp^\de)=\vp$, and $\vp^{\de,0}\in \Nlf_{\Xx^\de}$ denotes the special case where all the points in $\vp $ are equipped with the zero-sequence as mark. The latter notation is used in contexts where we want to lift a given $\vp\in \Nlf_{\Xx}$ to an element of $\Nlf_{\Xx^\de}$ and the choice of marks plays no role.
In practice, it will be enough to consider a countable sequence of $\de$'s tending to 0, so that properties that hold almost surely for each $\de$ also hold almost surely for all $\de$ simultaneous. However, we will not make this sequence explicit. 

We lift the ordering $\iota$ to a map $\iota^\de:\Xx \times \Uu^\de \to [0,1]$ by $\iota^\delta(\xx,\rr)=\iota(\xx)$, which defines a partial order on $\Xx^\delta$ by $(\xx,\rr)\le_{\iota^\de}(\yy,\ss)$ if $\xx<_\iota \yy$ or $(\xx,\rr)=(\yy,\ss)$. 
Define
\begin{equation}\label{eq:def_Rde}
R^\delta((\xx,\rr),(\yy,\ss)) = \sum_{i,j\in\N} \one_{D^\delta_i}(\xx) \one_{D^\delta_j}(\yy) (\one\{(\yy,\ss)<_{\iota^\de} (\xx,\rr)\} r_{j} + \one\{(\xx,\rr)<_{\iota^\de} (\yy,\ss)\}s_{i}) 
\end{equation}
where $\rr=(r_i)_{i\in \N}$ and $\ss=(s_i)_{i\in \N}$.
Then, a symmetric relation $\sim_\delta$ on $\Xx^\de$ is given by 
\begin{equation}\label{eq:def_simde}
(\xx,\rr)\sim_\delta (\yy,\ss) \text{ iff } R^\delta((\xx,\rr),(\yy,\ss))\leq \pi(\xx,\yy).
\end{equation}
For $\vp^\de,\psi^\de \in \Nlf_{\Xx^\de}$, we can define the graph $G^\de(\vp^\de,\psi^\de)$ on $\Xx^\de$ with vertex set $\vp^\de\cup\psi^\de$ and an edge between $(\xx,\rr)\in \vp^\de$ and $(\yy,\ss)\in \vp^\de\cup\psi^\de$ if $(\xx,\rr)\sim_\de (\yy,\ss)$.

To see how we can approximate the RCM on $Q$ by graphs on $Q^\de$, let $\Phi_{Q,\psi}^*$ be a Poisson process on $\bar Q^*$ and let $\psi\in \Nlf_{\Xx\sm Q}$ be finite. Define Poisson processes $\Phi^{*,\de}_{Q,\psi}$ on $Q^{*, \delta}:=Q^*\times \Uu^\de$ coupled with $\Phi_{Q,\psi}^*$ as follows.
A point $(\xx,u,\rr)\in \Phi^*_{Q,\psi}$ corresponds to  $(\xx,u,\rr^\de)\in \Phi^{*,\delta}_{Q,\psi}$, where 
$r_i^\de=r_{\xx,\yy_i}$
for $\yy_i=\inf_{\iota}\big\{\yy \in D_i^\delta\cap ( \Pd_Q\cup \psi ) \big \}$
when $  D_i^\delta \cap ( \Pd_Q\cup \psi ) \neq \es$, and otherwise $r_i^\delta$ is chosen uniformly at random. Then, $\Phi^{*,\delta}_{Q,\psi}$ projected to $ Q^*$ is $\Pds_Q$ and hence a Poisson process. Given $\Pds_Q$, the marks in $\Uu^\delta$ are independent uniform i.i.d.\ sequences.  Thus, $\Phi^{*,\delta}_{Q,\psi} $ is a Poisson process on $Q^*\times \Uu^\delta $. 
We can form the graph $G^\de\big(\Phi_{Q,\psi}^{*,\de},\psi^{\de,0}\big)$ on $\Xx^\de$.   If $\psi$ is finite, there is a $\de_0(\omega) $  for almost all $\omega$ such that $\Pd_Q\cup \psi$ is separated by the partition whenever $\de < \de_0$. For such $\de$, the projected graph $h^\de\big(G^\de\big(\Phi_{Q,\psi}^{*,\de}, \psi^{\de,0}\big)\big)$ on $\Xx$ will agree with $\Gamma(\Phi_{Q,\psi}^{*}, \psi)$. Therefore, we think of $h^\de\big(G^\de\big(\Phi_{Q,\psi}^{*,\de}, \psi^{\de,0}\big)\big)$ as approximations of the RCM $\Gamma\big(\Phi_{Q,\psi}^{*}, \psi\big)$ that converge almost surely when $\de\to 0$.

\subsubsection{Approximation of the Gibbs process}\label{sss:gibbsapprox}

Define a PI on $\Xx^\de$ defined for $(\xx,\rr)\in\Xx^\de$ and $\vp^\de\in \Nlf_{\Xx^\de}$ by
\begin{equation}\label{eq:def_kde}
\kappa^\delta((\xx,\rr),\vp^\de) = \one\{(\xx,\rr)\not \sim_\delta \vp^\de \}.
\end{equation}
Let $\X^\de(Q^{\delta},\psi^{\de,0})$ denote a Gibbs process on $Q^\de$ with PI~$\kappa^\delta$ with boundary condition  $\psi^{\de,0}\in\Nlf_{\Xx^\de\sm Q^\de}$ where $\psi\in \Nlf_{\Xx \sm Q}$. This $\kappa^\de$ satisfies \eqref{cond:PIbound} and \eqref{cond:PI_sim} with the relation $\sim_\de$. 
It is shown in \cite[Lem. 3.1]{betsch} that for finite $\psi$, 
$$\P(\X(Q,\psi) \in E) = \lim_{\delta \to 0} \P( h^\de(\X^\de(Q^\de,\psi^{\de,0})) \in E)$$
for any $E\in \NN$. In fact, the proof shows that 
\begin{equation}\label{eq:delta_dtv}
	\lim_{\delta\to 0}\dtv(h^\de(\X^\de(Q^\de,\psi^{\de,0})) ,\X(Q,\psi)) = 0.
\end{equation}
However, no explicit couplings of these processes were given in \cite{betsch}. Such couplings can be provided via the coupled Poisson processes  $\Phi_{Q,\psi}^{*,\de}$. The processes $T_{Q^\de,\psi^{\de,0},\iota^\de}(\Phi_{Q,\psi}^{*,\de})$
 provide a coupling between the processes $\X^\de(Q^\de,\psi^{\de,0})$ and we will see in the next section how to couple these approximating processes with $\X(Q,\psi)$.

 Note here that the map $\iota^\de$ is not injective, but $(\la \otimes \Leb_{[0,1]}^\N)((\iota^\de)^{-1}(t))=0$ for all $t\in \R$. Hence $\iota^\de$ defines a total ordering on $\Phi_{Q,\psi}^{*,\de}$ a.s., and Proposition \ref{pr:iota} shows that the Poisson embedding $T_{Q^\de,\psi^{\de,0},\iota^\de}(\Phi_{Q,\psi}^{*,\de})$ is well-defined a.s.~and that the results of Section \ref{ss:fw} carry over. This was shown in the Euclidean case in \cite[Prop. 4.2]{hos}. The case of general Borel spaces is similar, thus we only sketch the  generalization below.  

\bepr \label{pr:iota}
Assume that $Q\in \BB_0$ and $\iota: Q \to \R$ is a  measurable map such that $
	\lambda(\iota^{-1}(t) )= 0
	$ for all $t\in \R$.
Then, $T_{Q, \psi,\iota}(\Pds_Q)$ is well-defined almost surely and is distributed as $\XX(Q, \psi)$.  Moreover, Proposition \ref{th:dc_prop}, Lemma \ref{lem:distv}, Corollary \ref{cor:poiss_emb}, and Proposition \ref{pr:xxff} also hold with this $\iota$. 
\enpr

\bep
First, we note that  $\iota$ defines a total ordering on $\Pds_Q$ almost surely and the thinning probabilities $p(\xx,Q,\psi)$ are measurable functions of $\xx$ and $\psi$. Hence, the standard Poisson embedding $T_{Q,\psi,\iota}(\Pds_Q)$ is well-defined almost surely. 
To see that this is distributed as $\X(Q,\psi)$, consider for each $k\ge 1$ a partition  of $Q$ of the form $D_i^k = \iota^{-1}((i2^{-k},(i+1)2^{-k}])$, $i\in \Z$. Choose an injective ordering  $\iota_i^k:D_i^k\to (i2^{-k},(i+1)2^{-k}] $ on each $D_i^k$ (e.g.\ choose an injective ordering  $\tilde{\iota}: Q  \to  (0,1]$ and let $\iota_i^k(\xx) = 2^{-k}(i+\tilde{\iota}(\xx))$). These $\iota^k_i$ are then pieced together to form an injective map $\iota^k:Q\to \R^d$.  Using these orderings, we obtain a Poisson embedding map $T_{Q,\psi,\iota^k}$ for each $k$ such that $T_{Q,\psi,\iota^k}(\Pds_Q)\stackrel{d}{=}\X(Q,\psi)$. One can now show, using the same reasoning as in the proof of \cite[Proposition 4.3]{hos},  that almost surely, there is a $K(\omega)$ such that for all $k>K$
$$T_{Q,\psi,\iota^k}(\Pds_Q){=}T_{Q,\psi,\iota}(\Pds_Q)$$
and that this implies $T_{Q,\psi,\iota}(\Pds_Q)\stackrel{d}{=} \XX(Q,\psi)$. 
For Lemma \ref{lem:distv}, we may assume that $\iota(A)\su (-\ff,0)$ and $\iota(Q\sm A)\su (1,\ff)$. Then, the lemma holds for each $\iota^k$ with the same $A$, i.e.\
$$
\P\big(\inf{}_{\iota^k} (T_{Q,\psi,\iota^k}(\Pds_Q) \cap A) \ne \inf{}_{\iota^k} (T'_{Q,\psi',\iota^k}(\Pds_Q)\cap A)\big)\le \Big(2+\int_{A} \alpha(\xx) \lambda(\d \xx )\Big)\,\dtv\big(\XX(Q,\psi)\cap A,\XX'(Q,\psi')\cap A\big).
$$
Letting $k\to \ff$, the right hand side is constant, while the left hand side converges to $\P\big(\inf{}_{\iota} (T_{Q,\psi,\iota}(\Pds_Q) \cap A) \ne \inf{}_{\iota} (T'_{Q,\psi',\iota}(\Pds_Q)\cap A)\big)$. Corollary \ref{cor:poiss_emb}  follows from Lemma \ref{lem:distv} (since the pushforward measure of $\iota$ has no atoms).
The proofs of Proposition \ref{th:dc_prop} and \ref{pr:xxff} carry over immediately.
\enp

\subsection{Proofs for Section \ref{ss:setup}}\label{ss:proofs_se_RCM}

Henceforth, we let $E^\de$ denote the set of finite configurations $\vp\in \Nlf_{\Xx}$ for which $\sup_{i\ge1}\vp(D_i^\de)\le 1$, that is, the $D_i^\de$ separate the points of $\vp$.

\begin{proof}[Proof of Proposition \ref{prop:pois_emb_RCM}]
First, consider the case  where $\psi$ is finite.   We  show that the Poisson embeddings $T_{Q^\delta,\psi^{\de,0}} (\Phi^{*,\delta}_{Q,\psi})$ project to coupled processes $h^\de\big(T_{Q^\delta,\psi^{\de,0}} (\Phi^{*,\delta}_{Q,\psi})\big)$ on $Q$ such that there is almost surely a $\delta_1(\omega)$ with  $\Trr_{Q,\psi} (\Phi^*_{Q,\psi}) = h^\de(T_{Q^\delta,\psi^{\de,0}} (\Phi^{*,\delta}_{Q,\psi}))$  for  $\de<\de_1$.

If $\{\xx_1,\dots, \xx_m\}\cup \psi \in E^\delta$, we have \cite[Lem. 3.1]{betsch} that
\begin{align*}
	\int_{\Uu^m}\tilde{\kappa}^\delta((\xx_1,\rr_1),\dots,(\xx_m,\rr_m), \psi^{\de,0})(\Q^{\N})^m(\d (\rr_1,\dots,\rr_m))& =\prod_i \prod_{\yy\in \psi} (1-\pi(\xx_i,\yy))\prod_{i<j} (1-\pi(\xx_i,\xx_j)) \\
	&= \kappa_m(\xx_1,\dots,\xx_m,\psi).
\end{align*}
It follows (see computation in \cite[Lem. 3.1]{betsch}) for any event $F$,
$$\lim_{\delta\to 0}\P(h^\de(\X^\de(Q^\de,\psi^{\de,0}))\in F) = \P(\X(Q;\psi)\in F),$$
and for $B\su Q$,
$$|Z_{B^\delta}(\psi^{\de,0}) - Z_{B}(\psi)|\leq \P(\psi \cup \Pd_{B} \notin E^\delta).$$
In particular, if $(\xx,\rr)\in Q^\delta$ and $\phi^\de\in \Nlf_{ \Xx_{(-\ff,(\xx,\rr))}^\de  }$ is finite, then the retention probabilities for  $T_{Q^\de,\psi^{\de,0}}$ are
\begin{align*}
	&p^\delta((\xx,\rr),Q^\delta,\phi^\de) = \kappa^\delta((\xx,\rr),\phi^\de)\frac{Z_{Q_{((\xx,\rr),\ff)}^\de}(\phi^\de \cup \{(\xx,\rr)\})}{Z_{Q_{((\xx,\rr),\ff)}^\de }(\phi^\de)} \\
	&= \kappa(\xx,\phi)\frac{\prod_{\yy\in\phi}\sum_i\one\{\yy\in D_i^\de,r_i> \pi(\xx,\yy)\}}{\prod_{\yy\in\phi}(1-\pi(\xx,\yy))} \frac{Z_{Q_{((\xx,\rr),\ff)}^\de}(\phi^\de \cup \{(\xx,\rr)\})}{Z_{Q_{((\xx,\rr),\ff)}^\de}(\phi^\de)}\\
	&= p(\xx,Q,\phi) \frac{\prod_{\yy\in\phi}\sum_i\one\{\yy\in D_i^\de,r_i> \pi(\xx,\yy)\}}{\prod_{\yy\in\phi}(1-\pi(\xx,\yy))}\frac{Z_{Q_{((\xx,\rr),\ff)}^\de}(\phi^\de \cup \{ (\xx,\rr)\})}{Z_{Q_{(\xx,\ff)}}(\phi\cup \{\xx\})} \frac{Z_{Q_{(\xx,\ff)}}(\phi )}{Z_{Q_{((\xx,\rr),\ff)}^\de}(\phi^\de)}
\end{align*}
(where $0/0=0$). Assume that $(\xx,u,\rr)\in \Phi^{*,\de}_{Q,\psi }$ and $\phi\su \Pd_{Q_{(-\ff,\xx)}} \cup \psi$. Then, for $\delta$ small enough, there is at most one $\yy\in  D_i^\de\cap \phi^\de$ for each $i$, and if such a $\yy$ exists, we have $r_i=r_{\xx,\yy}$. Letting $\delta\to 0$, we thus get almost surely
$$ \lim_{\de \to 0} p^\delta((\xx,\rr),Q^\de,\phi^\delta) = p(\xx,Q,\phi) \frac{\one\{\xx\nsim_{\Gamma(\Phi_{Q,\psi}^*,\psi)} \phi\}}{\prod_{\yy\in\phi}(1-\pi(\xx,\yy))}=p^{\ms{RCM}}(\xx,Q,\phi,\Gamma(\Phi_{Q,\psi}^*,\psi)).$$

This shows that the Poisson embeddings $T_{Q^\delta,\psi^{\delta,0}} (\Phi^{*,\delta}_{Q,\psi})$ project to coupled processes $h^\de(T_{Q^\delta,\psi^{\delta,0}} (\Phi^{*,\delta}_{Q,\psi}))$ which, whenever $\delta<\de_1(\omega)$,  equals  $T_{Q,\psi}^{\ms{RCM}} (\Phi^*_{Q,\psi})$ almost surely. By \eqref{eq:delta_dtv}, the  distribution of $T_{Q,\psi}^{\ms{RCM}} (\Phi^*_{Q,\psi})$ must be that of $\X(Q,\psi)$.

For the case where $\psi$ is infinite, let $\phi_n= \phi \cap W_n$. Then $p(\xx,Q,\phi_n,\Gamma_n) \to p(\xx,Q,\phi, \Gamma)$ by dominated convergence applied to the partition functions. Thus, there is almost surely an $N(\omega)$ such that for $n\ge N(\omega)$,  $T_{Q,\psi_n}^{\ms{RCM}} (\Phi^*_{Q,\psi})$ will agree with $T_{Q,\psi}^{\ms{RCM}} (\Phi^*_{Q,\psi})$. Taking limits in the GNZ-equation for $ T_{Q,\psi_n}^{\ms{RCM}} (\Phi^*_{Q,\psi}) \stackrel d= \X(Q,\psi_n)$ we see that $ T_{Q,\psi}^{\ms{RCM}} (\Phi^*_{Q,\psi}) $ satisfy the GNZ-equation for $ \X(Q,\psi)$ and hence has the claimed distribution.
\end{proof}

We note for later that the Poisson embedding has the  following property: If $\xx$ is the $\iota$-smallest point in $\Pds_Q$ and $\xi\in \{\{\xx\},\es\}$ denotes the thinning of $\xx$, then  
\begin{equation}\label{eq:RCM_se_first}
	\Trr_{Q,\psi}(\Phi_{Q,\psi}^*)=	\xi \cup \Trr_{Q_{(\xx,\ff)},\psi\cup \xi}(\Phi_{Q_{(\xx,\ff)},\psi\cup \xi }^*),
\end{equation}
where $\Phi_{Q_{(\xx,\ff)},\psi\cup \xi}^*$ is constructed from $\Phi_{Q,\psi}^*$ in the obvious way.

\subsection{Proofs for Section \ref{ss:dc_RCM}}\label{ss:proofs_dc_RCM}
\bep[Proof of Lemma \ref{lem:RCM_it}]
Clearly,  $\Pd_Q \sm V$ is an independent thinning of $\Pd_Q$ with retention probabilities given by $\rho_\nu$, so the first claim follows from e.g.\ \cite[Prop. 3.7]{rasmus}.  Thus, conditionally on $\GGG_\nu$, $\Pd_Q\sm V$ is a Poisson process of intensity $\rho_\nu \la$, and the marks used for defining the edges in $\GGB(\Phi^*_{Q\sm V,\psi\cup\nu} ,V\cup (\psi\sm \nu))$ are uniform i.i.d.\ and independent of $\GGG_\nu$. Hence $\GGB(\Phi^*_{Q\sm V,\psi\cup\nu} ,V\cup (\psi\sm \nu))$ is an RCM $\Gamma^{\rho_\nu\la}(Q,V\cup (\psi\sm \nu) )$ and $\Gamma_\nu(\Phi^*_{Q,\psi\cup\nu},\psi)\stackrel d= \Gamma(\Phi^*_{Q,\psi\cup\nu},\psi) \stackrel d=\Gamma(Q,\psi)$ for any $\GGG_\nu$. Thus, the same holds unconditionally.
\enp

\begin{proof}[Proof of Lemma \ref{lem:TRCM_V}]
	Consider first the situation where $\psi$ and $\nu$ are finite. Lift $\psi,\nu$ to $\psi^{\delta,0},\nu^{\delta,0} \in \Nlf_{\Xx^\de } $ by adding the zero mark. This is needed for formal reasons only and  we will never use these marks in what follows. 
	Consider a standard Poisson embedding $T_{Q^\de,\psi^{\de,0},\iota_\de}\big({\Phi}^{*,\delta}_{Q,\psi \cup \nu}\big)$ where we define the deterministic set $V_\de= \{(\xx,\rr)\in Q^\de \mid \xx\sim_\de \nu^{\de,0} \}$ and the ordering 
	\begin{equation}\label{eq:idelta}
		\iota_\de(\xx,\rr) =\tfrac12\iota(\xx)\one\{ (\xx,\rr) \in V_\de\} + \big(\tfrac12 + \tfrac12\iota(\xx)\big)\one\{ (\xx,\rr) \notin V_\de \}. 
	\end{equation}
	This ordering is used instead of $\iota^\de$ in \eqref{eq:def_Rde} leading to a new symmetric relation $\sim^\de$ in \eqref{eq:def_simde} which is then used for defining the PI $\k_\de$ as in \eqref{eq:def_kde}. This $\k_\de$ is used for defining $T_{Q^\de,\psi^{\de,0},\iota_\de}\big({\Phi}^{*,\delta}_{Q,\psi \cup \nu}\big)$.
 	We will show that almost surely there is a $\de_1(\omega)$ such that for $\de<\de_1$, $h^\de(T_{Q^\de,\psi^{\de,0},\iota_\de}({\Phi}^{*,\delta}_{Q,\psi\cup \nu})) = T_{Q,\nu,\psi}^{\ms{RCM}}(\Phi^*_{Q,\psi\cup \nu})$. 
	
	Let $E^\de$ be as in the proof of Proposition \ref{prop:pois_emb_RCM}.
	To determine whether $(\xx,u,\rr)\in V_\delta^* \cap {\Phi}^{*,\delta}_{Q,\psi\cup \nu}$, we only need to know the coordinates  $r_i$ of $\rr$ for which $D_i\cap \nu \ne \es$. Thus, if $\nu\cup\psi \cup \Pd_Q\in E^\de$, then $r_i=r_{\xx,\yy}$ whenever $\yy\in D_i\cap \nu$ and hence $h^{*,\de}(V_\delta^* \cap {\Phi}^{*,\delta}_{Q,\psi\cup \nu})=V$. There is almost surely a $\de_0(\omega)$, such that $\nu\cup \psi \cup\Pd_Q\in E^\de$ for $\de < \de_0$ and hence $h^{*,\de}(V_\delta^* \cap {\Phi}^{*,\delta}_{Q,\psi\cup \nu})=V$  for $\delta<\delta_0(\omega)$. Moreover, under $h^{\de}$ the graph $G^\de({\Phi}^{*,\delta}_{Q,\psi\cup \nu}, \psi\cup \de)$ associated with $\sim^\de$ is mapped to  $h^{\de}(G^\de({\Phi}^{*,\delta}_{Q,\psi\cup \nu}\cup \psi\cup \nu)) =\G_\nu({\Phi}^{*}_{Q,\psi\cup \nu}, \psi\cup \nu)$ whenever $\delta<\delta_0(\omega)$ and the ordering of the vertices associated with $\iota_\de$ agrees with the thinning order for $T_{Q,\nu,\psi}^{\ms{RCM}}(\Phi^*_{Q,\psi\cup \nu})$.
	
	It remains to identify the retention probabilities. Recall that for a given $\delta$ and $\phi^\de\su  Q_{(-\ff,(\xx,\rr))}^\de \cup \psi^{\de,0}$, the retention probabilities on $Q^\de$ are
	$$p^\delta((\xx,\rr),Q^\de, \phi^\de) = \kappa_\delta((\xx,\rr),\phi^\de) \frac{Z_{Q_{((\xx,\rr),\ff)}^\de}(\phi^\de \cup \{(\xx,\rr)\} )}{Z_{Q_{((\xx,\rr),\ff)}^\de}(\phi^\de)}.$$
	If $\iota_\de(\xx,\rr)<\iota_\de(\xx_1,\rr_1)<\iota_\de(\xx_2,\rr_2) < \dotsm < \iota_\de(\xx_m,\rr_m)$  and $\phi\cup \{\xx_1,\dots,\xx_m\}  \in E^\delta$, then
	$$\tilde{\k}^\delta((\xx_1,\rr_1),\dots,(\xx_m,\rr_m);\phi^\de) = \prod_i\prod_{\yy\in\phi} \bigg(\sum_{\yy\in D_j}\one\{(\rr_i)_j >\pi(\xx_i,\yy)\} \bigg)\prod_{i<j} \bigg(\sum_{\xx_i\in D_l}\one\{(\rr_j)_l>\pi(\xx_i,\xx_j)\}\bigg)$$
	where $(\rr_i)_j$ is the $j$th coordinate of $\rr_i$.
	In particular, if $\delta <\de_0$, $\phi^\de\su  (\Phi^{*,\de}_{Q,\psi\cup\nu}\cap Q_{(-\ff,(\xx,\rr))}^\de ) \cup \psi^{\de,0}$ 
	and $(\xx,u,\rr)\in \Phi^{*,\de}_{Q,\psi\cup\nu}$, it holds almost surely that
	$$\kappa^\delta((\xx,\rr),\phi^\de) = \prod_{\yy\in\phi^\de}\one\{r_{\xx,\yy} >\pi(\xx,\yy)\}  =\one\{\xx\nsim_{\Gamma_\nu(\Phi^{*}_{Q,\psi\cup\nu},\psi\cup\nu)} \phi \}.$$
	
	Next, we determine the limit of $Z_{Q_{((\xx,\rr),\ff)}^\de}(\phi^\de)$ when $\delta \to 0$. Suppose $(\xx,\rr)\in V_\de$. Then, the set $Q_{((\xx,\rr),\ff)}^\de$ is necessarily of the form 
	\begin{align*}
		A_\delta &= (V_\de \cap (h^\de)^{-1}(Q_{(\xx,\ff)}) )\cup (Q^\de \sm V_\de)  = (h^\de)^{-1}(Q_{(\xx,\ff)}) \cup \{(\yy,\rr) \nsim_{\delta}  \nu^{\delta,0},  \yy \le_\iota \xx \}.
	\end{align*}
	We have the Poisson expansion
	\begin{align*}
		&	e^{(\lambda\otimes \Q^{\delta})(A_\de)}Z_{A_\delta}(\phi^\de)\\
		& = 1 + \sum_{m=1}^\ff\frac1{m!} \int_{(A_\delta)^m}  \tilde{\kappa}^\delta((\xx_1,\rr_1),\dots,(\xx_m,\rr_m);\phi^\de) (\lambda\otimes \Q^{\delta})^m (\d(\xx_1,\rr_1),\dots,\d (\xx_m,\rr_m))\\
	\end{align*}
	Assume $\phi\cup\{\xx_1,\dots,\xx_m\} \cup \nu \in E^\delta$ and  write $({\rr}_i)_{\yy}$ for the coordinate of ${\rr}_i$ corresponding to $\yy\in  \phi \cup \{\xx_1,\dots,\xx_m\}\cup \nu$.
	 If   $\phi^\de \su Q_{(-\ff,(\xx,\rr)]}^\de$ and $\iota_\de(\xx,\rr) <\iota_\de(\xx_1,\rr_1)<\dotsm < \iota_\de(\xx_m,\rr_m)$,
	\begin{align}\nonumber
		&\prod_l \one\{(\xx_l,{\rr}_l)\in A_\delta\}\kappa_m^\delta((\xx_1,{\rr}_1),\dots,(\xx_m,{\rr}_m);\phi^\de) \\ \nonumber
		&=\prod_{l: \xx_l<_\iota \xx} \prod_{\yy\in\nu} \one\{({\rr}_l)_\yy>\pi(\xx_l,\yy)\} \prod_{k \le m}\prod_{\yy\in \phi}\one\{({\rr}_k)_\yy > \pi(\xx_k,\yy)\} \prod_{i<j}\one\{({\rr}_{j})_{\xx_i}> \pi(\xx_i,\xx_j)\} \\
		&=\prod_{l: \xx_l <_\iota x} \one\{\xx_l \nsim_{\Gamma_m} \nu\}\prod_{k \le m}\one\{\xx_k \nsim_{\Gamma_m} \phi\} \prod_{i< j}\one\{\xx_i \nsim_{\Gamma_m} \xx_j\}, \label{eq:Adelta}
	\end{align}	
	where $\Gamma_m=\G(\{(\xx_1,{\rr}_1),\dots,(\xx_m,{\rr}_m)\}, \phi \cup \nu )$.
	Thus, 
	\begin{align*}
		&\bigg|  Z_{A_\delta}(\phi^\de) -e^{-(\lambda\otimes \Q^{\delta})(A_\de)}\bigg( 1 + \sum_{m=1}^\ff\frac1{m!} \int_{Q^m} \prod_{l: \xx_l <_\iota \xx} \one\{\xx_l \nsim_{\Gamma_m} \nu\}\prod_l\one\{\xx_l \nsim_{\Gamma_m} \phi\} \prod_{i< j}\one\{\xx_i \nsim_{\Gamma_m} \xx_j\} \\
		&\qquad \qquad \big(\la\otimes \Q^{\de}\big)^m (\d (\xx_1,{\rr}_1),\dots,\d(\xx_m,{\rr}_m) )\bigg) \bigg|\\
		& \leq \P(\phi \cup \nu\cup  \Pd_Q \in E^\delta).
	\end{align*}
	In particular, it follows from \eqref{eq:Adelta} with $m=1$ and $\phi=\es$ that $$\lim_{\de \to 0}\big(\lambda\otimes \Q^{\delta}\big)(A_\de) = \int_Q \rho_\nu^\xx (\yy)\lambda(\d \yy) =: c(\xx,\nu),$$ so 
	\begin{align}\nonumber 
		&\lim_{\delta \to 0} Z_{A_\delta}(\phi^\de) = e^{-c(\xx,\nu)}\bigg( 1 + \sum_{m=1}^\ff\frac1{m!} \int_{Q^m} \prod_{\xx_l<_\iota \xx} \one\{\xx_l \nsim_{\Gamma_m}  \nu\cup\phi\}\prod_{\xx_k >_\iota \xx}\one\{\xx_k \nsim_{\Gamma_m} \phi\} \prod_{i< j}\one\{\xx_i \nsim_{\Gamma_m} \xx_j\}\\ \nonumber
		&\qquad \qquad \big(\la\otimes \Q^{\de }\big)^m  \big(\d (\xx_1,{\rr}_1),\dots,\d (\xx_m,{\rr}_m)\big) \bigg)\\ \nonumber
		&=  e^{-c(\xx,\nu)}\bigg( 1 + \sum_{m=1}^\ff\frac1{m!} \int_{Q^m}\prod_l\rho_\nu^\xx(\xx_l) \E\bigg[ \prod_{\xx_l<_\iota \xx}\one\{\xx_l \nsim_{\Gamma_m} \phi \sm \nu\} \prod_{k:\xx_k>_\iota \xx}\one\{\xx_k \nsim_{\Gamma_m} \phi\}\prod_{i< j}\one\{\xx_i \nsim_{\Gamma_m} \xx_j\}\bigg]\\ \nonumber
		&\qquad \qquad  \la(\d \xx_1)\ldots \la(\d \xx_m)\bigg)\\
		&=  {\E}\bigg[ \prod_{\yy\in \PP_{Q_{(-\ff,\xx)}}^{\rho_\nu^\xx\lambda}} \one\{ \yy \nsim_{\Gamma^{\rho_\nu^\xx\lambda}(Q,\phi\cup \nu)} \phi\sm \nu\}\prod_{\yy\in \PP_{Q_{(\xx,\ff)}}^{\rho_\nu^\xx\lambda}}\one\{ \yy \nsim_{\Gamma^{\rho_\nu^\xx\lambda}(Q,\phi\cup \nu)}\phi\} \prod_{\yy_1,\yy_2\in \PP_Q^{\rho_\nu^\xx\lambda }} \one\{\yy_i \nsim_{\Gamma^{\rho_\nu^\xx\lambda}(Q,\phi\cup \nu)} \yy_j\} \bigg].\label{eq:lim_Zdelta}
	\end{align}
	In the second equality, we first integrated over the $\nu$-coordinates of ${\rr}_l$ for those $l$ with $\xx_l<_\iota \xx$ and then with respect to the remaining $\rr$-coordinates. In the last equality, we recognized the expression as an expectation with respect to a Poisson process $\PP_Q^{\rho_\nu^\xx\la}$ of intensity $\rho_\nu^\xx\la$.
	
	After thinning $V$, we have that $\Pd_Q\sm V$ is an independent Poisson process of intensity $\rho_\nu\lambda$ and for $\de<\de_0$, $h^{*,\de}(\Phi^{*,\de}_{Q,\psi\cup \nu}\sm V_\de)=\Pd_Q\sm V$ with the ordering given by $\iota$. Computing as above, we find when $\psi \cup \nu \cup \Pd_Q\in E^\de$ and $(\xx,u,\rr)\in \Phi^{*,\de}_{Q,\psi\cup\nu}\sm V_\de^*$,
	$$\kappa^\delta((\xx,\rr),\phi^\de) =\one\{\xx\nsim_{\Gamma_\nu(\Phi^*_{Q,\psi\cup\nu},\psi\sm \nu)} \phi\sm \nu\}$$
	and for $(\xx,\rr)\notin V_\de$ and $\phi \su \psi \cup V \cup \Pd_{Q_{(-\ff,\xx)}}$
	\begin{align*}
		&\lim_{\de \to 0}{Z}_{Q^\de_{((\xx,\rr),\ff)}}(\phi^\de)={\E}\bigg(	\prod_{\yy\in \PP_{Q_{(\xx,\ff)}}^{\rho_\nu\la}} \one \{\yy \nsim_{\Gamma^{\rho_\nu \la}(\PP_{Q_{(\xx,\ff)}}^{\rho_\nu\la},\phi\sm \nu)} \phi\sm \nu \} \prod_{\yy_1,\yy_2\in \PP_{Q_{(\xx,\ff)}}^{\rho_\nu\la}}\one \{\yy_1 \nsim_{\Gamma^{\rho_\nu \la}(\PP_{Q_{(\xx,\ff)}}^{\rho_\nu\la},\phi\sm \nu)} \yy_2\} \bigg)\\
		&=Z_{Q_{(\xx,\ff)}}^{\rho_\nu\lambda}(\phi\sm \nu ).
	\end{align*}
	Thus, we recognize $\lim_{\de\to 0}p^\de(\xx,Q^\de,\phi^\de)$ as the retention probabilities \eqref{eq:ret_RCM} for the embedding of $\X^{\rho_\nu\la}(Q, \xi\cup (\psi \sm \nu))$ in an RCM $\Gamma^{\rho_\nu\la}(Q,\xi \cup (\psi \sm \nu))$.
	
	In total, we have that for almost all realizations of $\Phi^{*}_{Q,\psi\cup\nu}$, there is a $\delta_1(\omega) > 0$ such that  $$h^\de\big(T_{Q^\de ,\nu^\de,\psi^{\de,0}}(\Phi^{*,\de}_{Q,\psi\cup\nu})\big) = \Trr_{Q,\nu,\psi}(\Phi^*_{Q,\psi\cup\nu})$$
	for all $\de <\de_1$. To see that $\Trr_{Q,\nu,\psi}(\Phi^*_{Q,\psi\cup\nu}) \stackrel d= \X(Q,\psi)$, note that
	\begin{align*}
	&e^{\la(Q)}\P(h^\de\big(T_{Q^\de ,\nu^\de,\psi^{\de,0}}(\Phi^{*,\de}_{Q,\psi\cup\nu})\big)\in F)= 1 +\\	&  \sum_{m=1}^\ff\frac1{m!} \int_{(Q^\de)^m} \one\{\{\xx_1,\ldots,\xx_m\}\in F\} \tilde{\kappa}^\delta((\xx_1,\rr_1),\dots,(\xx_m,\rr_m);\psi^{\de,0}) (\lambda\otimes \Q^{\delta})^m (\d(\xx_1,\rr_1),\dots,\d (\xx_m,\rr_m)).
		\end{align*}
		By the considerations above, the right hand side equals
		$$
	 1 + \sum_{m=1}^\ff\frac1{m!} \int_{Q^\de} \one\{\{\xx_1,\ldots,\xx_m\}\in F\} \tilde{\kappa}(\xx_1,\dots,\xx_m;\psi) \lambda( \d\xx_1) ,\dots,\la(\d \xx_m)
$$
plus an error term of at most $\P(\Pd_Q\cup\psi\cup\nu \notin E^\de)$. Letting $\de\to \infty$ and using dominated convergence, we get that $\P(\Trr_{Q,\nu,\psi}(\Phi^*_{Q,\psi\cup\nu})\in F) = \lim_{\de \to 0}\P(h^\de\big(T_{Q^\de ,\nu^\de,\psi^{\de,0}}(\Phi^{*,\de}_{Q,\psi\cup\nu})\big)\in F) = \P(\X(Q,\psi)\in F)$.
By construction, we see that conditionally on $\xi$, $\Trr_{Q,\nu,\psi}(\Phi^*_{Q,\psi\cup\nu})\sm V \stackrel d= \X^{\rho_\nu\la}(Q,\xi\cup (\psi\sm \nu))$.
	
	The case where $\psi$ is infinite follows by taking limits of $\psi_n = \psi \cap W_n$ for $n\to \ff$ as in the proof of Proposition~\ref{prop:pois_emb_RCM}. When $\nu$ is infinite, the theorem holds for $\nu_n=\nu\cap W_n$. As $n\to \ff$, the set $V\su \vp$ increases, and since $\vp$ is finite, $V$ stabilizes as the set of all points connected to $\nu_n$ for $n\ge M$ for some $M(\omega)$.
	Thus, it is again enough to see that the retention probabilities in both stages of the thinning procedure converge. For the thinning of $V$, note that $\rho_{\nu_n}$ and $\rho_{\nu_n}^\xx$ both decrease towards $\rho_{\nu}$ and $\rho_{\nu}^\xx$, respectively, when $n\to \ff$. To take the limit in \eqref{eq:Z_bar_def}, note that we may couple the RCMs $\Gamma^{\rho_{\nu_n}^\xx}(Q,\psi\cup \nu_n)$ to a common RCM $\Gamma(Q\times[0,1],\psi\cup \nu)$. Indeed, we take the vertex set of $\Gamma^{\rho_{\nu_n}^\xx}(Q,\psi\cup \nu_n)$ to be $\psi\cup \nu_n$ together with those $\xx\in Q$ for which $(\xx,u) $ a vertex of $\Gamma(Q\times[0,1],\psi\cup \nu)$ and $u\leq \rho_{\nu_n}^\xx(\xx)$. As edges of $\Gamma^{\rho_{\nu_n}^\xx}(Q,\psi\cup \nu_n)$ we take all edges induced from $\Gamma(Q\times[0,1],\psi\cup \nu)$. Then, almost surely, $\Gamma^{\rho_{\nu_n}^\xx}(Q,\psi\cup \nu_n)$ will equal $\Gamma^{\rho_{\nu}^\xx}(Q,\psi\cup \nu)$ for $n$ sufficiently large (except for the extra points $\nu\sm \nu_n$). Hence, dominated convergence shows that $Z_Q(\phi,\nu_n,\xx) \to Z_Q(\phi,\nu,\xx)$ when $n\to \ff$. A similar argument applies for the retention probabilities in \eqref{eq:ret_RCM}.

\end{proof}

\begin{remark}
	In the proof, we redefine the relation $\sim_\de$ to $\sim^\de$ when we change the ordering to $\iota_\de$, which is different from what we do in the first step of a disagreement coupling for Model (L) where the ordering is fixed throughout. This ensures that the thinning probabilities depend only on the marks of the unexplored points. Without this assumption, the limit of $Z_{A^\de}(\phi^\de)$ would depend on the marks of $\phi^\de$ and we would not be able to obtain formulas for the retention probabilities in the limit $\de\to 0$.  This fact is going to complicate the formulation of disagreement coupling for Model (P) in the following.
\end{remark}

\bep[Proof of Proposition \ref{prop:dc_RCM}] 
We first show by induction that conditionally on $(\GGd_k)^*$, $\Pd_Q\sm S_k$ is a Poisson process of intensity $\rho_k\la$ on $Q_{(\xx_k,\ff)}$ where $\rho_0= 1$ and
$$ \rho_k(\xx) := \rho_{k-1}(\xx) \prod_{\yy\in \nu_k}(1-\pi(\xx,\yy)).$$
Assume that this holds for some $k$ and condition on $(\GGd_{k+1})^*$. Then, $\Pd_Q\sm (S_k\cup V_{k+1})$ is a Poisson process of intensity $\rho_{k+1}\la$ on $Q_{(\xx_k,\ff)}$  by Lemma \ref{lem:RCM_it}. When $V_{k+1}\ne \es$, $S_k\cup V_{k+1}=S_{k+1}$ and $Q_{(\xx_k,\ff)}= Q_{(\xx_{k+1},\ff)}$ showing the claim in this case. When $V_{k+1}=\es$, $$\Pd_Q\sm S_{k+1} = \Pd_Q\sm (S_k\cup V_{k+1} \cup \{\xx_{k+1}\}) = (\Pd_Q\sm (S_k\cup V_{k+1})) \cap Q_{(\xx_{k+1},\ff)}. $$ Since $Q_{(\xx_{k+1},\ff)}$ is a stopping set and we have conditioned on $\xx_{k+1}$, $\Pd_Q\sm S_k \cap Q_{(\xx_{k+1},\ff)}$ is a Poisson process of intensity $\rho_{k+1}\la$ on $Q_{(\xx_{k+1},\ff)}$ showing the induction step in this case.

Define the graph 
\begin{equation}\label{eq:Gammak}
\Gamma_k:=\GGd_k \cup \Gamma_{\nu_{k+1}}(\Phi^*_{Q\sm S_k,\psi\cup\nu},S_k\sm\cup_{l\le k}\nu_l).
\end{equation}
It follows from Lemma \ref{lem:RCM_it} and the Poisson distribution of $\Pd_Q\sm S_k$ that conditionally on $(\GGd_k)^*$ we have 
\begin{equation}\label{eq:gam_nu_dist}
	\Gamma_{\nu_{k+1}}(\Phi^*_{Q,\psi\cup\nu}\sm \bar S_k^*,S_k\sm\cup_{l\le k}\nu_l) \stackrel d=  \Gamma^{\rho_k\la}(Q_{(\xx_k,\ff)},S_k \sm(\cup_{l\le k} \nu_l)).
\end{equation}

We will show by induction that $\Gamma_k\stackrel d=\Gamma(Q,\psi)$.  Note that $\Gamma_0=\Gamma_{\nu_1}(\Phi^*_{Q,\psi\cup \nu},\psi)\stackrel d=\Gamma(Q,\psi)$ by Lemma \ref{lem:RCM_it}. Suppose it holds for some $k\ge 1$.
By Definition \ref{def:rcm_nu},
$$\Gamma_k= \GGd_k \cup  \GGG_{\nu_{k+1}}(\Phi^*_{Q,\psi\cup \nu}\sm \bar S_k^*, S_k \sm(\cup_{l\le k} \nu_l))\cup \GGB_{\nu_{k+1}}(\Phi^*_{Q,\psi\cup \nu}\sm \bar S_k^*,V_{k+1} \cup S_k \sm(\cup_{l\le k+1} \nu_l)).$$
Now, condition on  $(\GGd_{k+1})^*$. If $V_{k+1}\neq \es$, then  $\GGd_k \cup \GGG_{\nu_{k+1}} = \GGd_{k+1}$. By Lemma \ref{lem:RCM_it} and the Poisson distribution of $\Pd_Q\sm S_k$,
$$\GGB_{\nu_{k+1}}(\Phi^*_{Q,\psi\cup\nu}\sm \bar S_k^*,V_{k+1} \cup S_k \sm(\cup_{l\le k+1} \nu_l))\stackrel d= \Gamma^{\rho_{k+1}\la}(Q_{(\xx_k,\ff)},S_{k+1}   \sm \cup_{l\le k+1} \nu_l).$$ 
If $V_{k+1}=\es$, then $\GGd_k \cup \GG_{\{\xx_{k+1}\}} = \GGd_{k+1}$. In this case, $S_{k+1}=S_k\cup \{\xx_{k+1}\}$ and $$\GGB_{\nu_{k+1}}(\Phi^*_{Q,\psi\cup\nu}\sm \bar S_k^*,V_{k+1} \cup S_k \sm(\cup_{l\le k+1} \nu_l))= \GG_{\{\xx_{k+1}\}}\cup \Gamma(\Phi^*_{Q,\psi\cup \nu}\sm \bar S_{k+1}^*, S_{k+1}  \sm (\cup_{l\le k+1} \nu_l)) $$
is the RCM constructed as in Definition \ref{def:qrcm}, which can be constructed by first finding the $\iota$-smallest point $\xx_{k+1}$ and all its connections to $S_k\sm (\cup_{l\le k+1}\nu_l)$ forming the graph $\GG_{\{\xx_{k+1}\}}$ and then connecting all points in $\Pd_Q\sm S_{k+1}$ to $\{\xx_{k+1}\} \cup (S_k\sm (\cup_{l\le k+1}\nu_l))=S_{k+1}\sm \cup_{l\le k+1}\nu_l$  forming the graph $\Gamma(\Phi^*_{Q,\psi\cup \nu}\sm \bar S_{k+1}^*, S_{k+1}\sm \cup_{l\le k+1}\nu_l)$. Thus,
$$\Gamma_k =\GGd_k \cup \GG_{\{\xx_{k+1}\}}\cup \Gamma(\Phi^*_{Q,\psi\cup \nu}\sm \bar S_{k+1}^*, S_{k+1}\sm \cup_{l\le k+1}\nu_l) ) = \GGd_{k+1} \cup  \Gamma(\Phi^*_{Q,\psi\cup \nu}\sm \bar S_{k+1}^*, S_{k+1}\sm \cup_{l\le k+1}\nu_l)  ). $$
Since  $\Pd_Q\sm S_{k+1}$ is a Poisson process on $Q_{(\xx_{k+1},\ff)}$ of intensity $\rho_{k+1}\la$, we have
$$\Gamma(\Phi^*_{Q,\psi\cup \nu} \sm \bar S_{k+1}^*, S_{k+1}\sm \cup_{l\le k+1}\nu_l)) \stackrel d= \Gamma^{\rho_{k+1}\la}(Q_{(\xx_{k+1},\ff)}, S_{k+1}  \sm (\cup_{l\le k+1} \nu_l)).$$
Thus, in both cases ($V_{k+1}=\es$ and $V_{k+1}\ne \es$), we get that conditionally on $(\GGd_{k+1})^*$,
$$\Gamma_k\stackrel d= \GGd_{k+1}  \cup \Gamma^{\rho_{k+1}\la}({Q_{(\xx_{k+1},\ff)}},S_{k+1}  \sm (\cup_{l\le k+1} \nu_l))\stackrel d= \Gamma_{k+1},$$
where the second equality in distribution comes from \eqref{eq:gam_nu_dist} with $k$ replaced by $k+1$.
Since $\Gamma_k \stackrel d= \Gamma(Q,\psi)$ unconditionally by the induction hypothesis, also $\Gamma_{k+1} \stackrel d= \Gamma(Q,\psi)$ unconditionally. Since $\GGd(\Phi^*_{Q,\psi\cup \nu},\psi)$ equals $\Gamma_k$ for $k\ge K(\omega)$, we obtain the result.
\enp

\begin{remark}\label{rem:directed}
	For arguing about the clusters of $\GGc_\nu$ in the proofs, it is sometimes convenient to view the cluster-based graph  $\GGc(\bar\vp^*,\psi)$ as a subgraph of a directed graph $\overrightarrow{\GG}(\bar\vp^*,\psi\cup\nu)$.  The vertex set of $\overrightarrow{\GG}$ is $\vp\cup \psi\cup \nu$ and for $\xx\in \vp$ and $\yy\in \vp\cup \psi\cup \nu$, there is a directed edge $\xx\to \yy$ if $r_{\xx,\yy} \leq \pi(\xx,\yy)$.   If we only keep edges $\xx\to \yy$ in $\overrightarrow{\GG}$ with $\xx >_\iota \yy$ (we refer to these as order-reversing edges) and ignore edge directions, we obtain  $\Gamma(\bar\vp^*,\psi\cup\nu)$. We obtain $V_\nu(\bar\vp^*,\psi)$ as the set of vertices in $\vp$ having a directed edge towards a point in $\nu$. If we keep the edges of $\overrightarrow{\GG}$ pointing from $\vp$ to $\psi $, edges pointing from $\vp \sm V_\nu$ to $V_\nu$ and order-reversing edges  $\xx\to \yy$ with  $\xx, \yy\in V_\nu$ or $\xx,\yy\in\vp \sm V_\nu$, then we get $\Gamma_\nu(\bar\vp^*,\psi)$. 
	
	The cluster-based graph $\GGc$ is built from $\overrightarrow{\GG}$ by first building the clusters connected to $\nu$. These clusters consist of all vertices in $\overrightarrow{\GG}$ having a directed path in $\overrightarrow{\GG}$ to a point in $\nu$. When all such vertices have been added, we start a new cluster at the $\iota$-smallest point not yet added to $\GGc$, say $\xx_k$. This $\xx_k$ has no directed path to a $\iota$-smaller point in $\vp$.  The cluster started at $\xx_k$ consists of all vertices in $\overrightarrow{\GG}$ having a directed path towards $\xx_k$, but not to any $\iota$-smaller point in $\vp$. The order in which the points of the cluster is considered for thinning is given by the graph distance to $\xx_k$ (i.e.\ the minimal number edges in a directed path to $\xx_k$). Two points with the same graph distance are ordered by $\iota$.
	When all points of this cluster have been explored, we start a new cluster from the $\iota$-smallest point in $\vp$ not yet considered and so on. 	
\end{remark}

\bep[Proof of Theorem \ref{thm:dc_RCM}]
Let 
$T^k:=\xi_1\cup \dots \cup \xi_{k-1} \cup \Trr_{Q_{(\xx_k,\ff)},\nu_k,B_k}(\Phi^*_{Q_{(\xx_k,\ff)},\nu\cup\psi}\sm \bar S_{k-1}^*).$
We show by induction that  $T^k\stackrel d= \X(Q,\psi)$. Note that  $T^1=\Trr_{Q_{(\xx_1,\ff)},\nu_1,\psi}(\Phi^*_{Q,\nu\cup\psi})\stackrel d= \X(Q,\psi)$ by Lemma \ref{lem:TRCM_V}.

Assume that the statement holds for $T^k$. Recall from the proof of Proposition \ref{prop:dc_RCM} that $\Pd_{Q_{(\xx_k,\ff)}}\sm  S_{k-1}^*$ is a Poisson process of intensity $\rho_{k-1}\la$ on $Q_{(\xx_k,\ff)}$. Given $(\GGd_k)^*$ with $V_k\ne \es$, we have by Lemma \ref{lem:TRCM_V} that  
$$T^k \stackrel d=\xi_1\cup\dots\cup  \xi_k \cup \X^{\rho_k\la} (Q_{(\xx_k,\ff)},\xi_k\cup B_k\sm \nu_k ) =\xi_1\cup\dots\cup  \xi_k \cup \X^{\rho_k\la} (Q_{(\xx_{k+1},\ff)}, B_{k+1}).$$
If $V_k=\es$, $\Trr_{Q_{(\xx_k,\ff)},\nu_k,B_k}$ is a standard Poisson embedding in an RCM, and hence the same holds by  \eqref{eq:RCM_se_first}.  By Lemma \ref{lem:TRCM_V}, 
the thinning 
$ \Trr_{Q_{(\xx_{k+1},\ff)},\nu_{k+1},B_{k+1}}(\Phi_{Q_{(\xx_{k+1},\ff),\psi\cup \nu}}^*\sm \bar S_k^*)$ can be used to  construct the process $\X^{\rho_k\la} (Q_{(\xx_{k+1},\ff)},B_{k+1})$. Thus, still conditionally on $(\GGd_k)^*$, 
$$T^{k} \stackrel d= \xi_1\cup\dots\cup  \xi_k \cup \Trr_{Q_{(\xx_{k+1},\ff)},\nu_{k+1},B_{k+1}}(\Phi_{Q_{(\xx_{k+1},\ff)}\sm S_k}^*)= T^{k+1}. $$
Since this holds for any $(\GGd_k)^*$, it also holds unconditionally,
so $T^k\stackrel d= \X(Q,\psi)$ implies $T^{k+1}\stackrel d= \X(Q,\psi)$. Letting $k\to \ff$, we find $\Trd_{Q,\nu,\psi}(\Phi^*)\stackrel d= \X(Q,\psi)$.
\enp

\subsection{Proofs for Section \ref{ss:pert_RCM}}\label{ss:proofs_perturb}

To construct $\bar T^{\ms{RCM,cl}}_{Q,\nu,\psi}$, we again use an approximation on $Q^\de$.
We introduce the PI on $\Xx^\de$ defined for $(\xx,\rr)\in\Xx^\de$ and $\vp^\de\in \Nlf_{\Xx^\de}$ by
\begin{equation*}\label{eq:simsim}
\bar{\kappa}^\delta((\xx,\rr),\vp^\de) = \bar\a(\xx)\one\{(\xx,\rr)\not \approx_\delta \vp^\de \},
\end{equation*}
where we used the symmetric relation $(\xx,\rr) \approx_\de (\yy,\ss)$ if and only if $R^\de((\xx,\rr),(\yy,\ss))\le \bar{\pi}(\xx,\yy)$, with $\bar{\pi}(\xx,\yy)=1-e^{-\bar{v}(\xx,\yy)}$. In particular, $(\xx,\rr) \approx_\de (\yy,\ss)$ implies $(\xx,\rr) \sim_\de (\yy,\ss)$, where $\sim_\de$ is the relation given by $R^\de((\xx,\rr),(\yy,\ss))\le {\pi}(\xx,\yy)$. This means that $\bar{\k}^\de$ satisfies \eqref{cond:PI_sim} with the symmetric relation $\sim_\de$. One can then define a thinning from fixed points by following the same steps as in the proof of Lemma \ref{lem:TRCM_V}  by setting  $V_\de =\{(\xx,\rr)\mid (\xx,\rr)\sim_\de \nu^{\de,0}\}$. We again use $V_\de$ to define a new ordering $\iota_\de$ (which agrees with the old one on $V_\de$). This new ordering $\iota_\de$ defines a relation $\approx^\de$ which can be used in \eqref{eq:simsim} to get a new PI $\bar \k_\de$.  We now use $\bar\k_\de$ and $\iota_\de$ in a Poisson embedding $T_{Q^\de,\nu^{\de,0},\psi^{\de,0}}(\Phi_{Q,\psi\cup\nu}^{*,\de})$ on $Q^\de$. One can show as in Lemma \ref{lem:TRCM_V} that when projecting to $Q$ and letting  $\de \to 0$ one obtains a copy of $\X(Q,\psi)$. We keep the part projected from $V_\de$ (i.e.\ the thinning of the points in $\G_\nu(\Phi_{Q,\nu\cup\psi}^*,\nu\cup\psi,\pi)$ connected to $\nu$) and iterate the construction as in Definition \ref{def:trd_RCM} following  a  cluster-based exploration of $\GGc(\Phi_{Q,\nu\cup\psi}^*,\nu\cup\psi,\pi)$. We denote the resulting point process by $\bar T^{\ms{RCM,cl}}_{Q,\nu,\psi}\big(\Phi_{Q,\psi\cup \nu}^{*}\big)$. One may show as in Theorem \ref{thm:dc_RCM}, that the result is a  point process with PI $\bar{\k}$.  We will not give the explicit formulas for retention probabilities, but they are very similar to \eqref{eq:ret_RCM} and \eqref{eq:p_nu}, only one has to be careful about which of the relations $\sim$ and $\approx$ are used where and how the activity shows up.

\bep[Proof of Theorem  \ref{thm:Talpha}]
  
 Let $\X_i $ and $\X_i'$ denote the part of the processes generated in the $i$th step of the algorithm  and let $\X_i^\de$ and $\X_i^{\prime,\de}$ denote the corresponding part of the processes generated on $Q^\de$. Then
 \begin{align*}
 &\P(\X_i\cap A \ne \X_i'\cap A, \X_j\cap A =\X_j'\cap A,j<i, \EpRCM(A,s)^c \cup \Ed^c\mid \GGc_{i-1}(\Phi_{Q,\es}^{*})^* ) \\
 &\le \lim_{\de \to 0} \P(\X_i^\de\cap A^\de \ne \X_i^{\prime,\de}\cap A^\de, \X_j \cap A=\X_j^{\prime}\cap A,j<i, \EpRCM(A^\de,s)^c \cup \Ed^c\mid \GGc_{i-1}(\Phi_{Q,\es}^{*})^* ).
 \end{align*} 
 One may argue as in the proof of Theorem \ref{thm:disagree} to see that 
 \begin{align*}
 	& \P(\X_i^\de\cap A^\de \ne \X_i^{\prime,\de}\cap A^\de, \X_j \cap A=\X_j^{\prime}\cap A ,j<i, \EpRCM(A^\de,s)^c \cup \Ed^c\mid \GGc_{i-1}(\Phi_{Q,\es}^{*})^*) \\
 	& \le 4\E\int_{B_{2s}(A^\de)}\P\Big((\xx,\rr)\lrsa_{\{(\xx,\rr) \}\cup \Phi_{Q\sm S_{i-1},\bar \eta}^{\de,*}\cup\bar{\eta}^\de} (\Phi_{Q\sm S_{i-1},\bar \eta}^{*,\de}\cap (Q^\de\sm P^\de))\cup \bar{\eta}^\de \Big) (\la\otimes \Leb_{[0,1]}^\de)(\d(\xx,\rr)) \\
 	&+ 2\E \int_{Q^\de \sm B_{2s}(A^\de) }\int_{B_{s}(A^\de)} \P((\xx,\rr)\sim^\de (\yy,\rr)) (\la\otimes \Leb_{[0,1]}^\de)(\d(\xx,\rr)) (\la\otimes \Leb_{[0,1]}^\de)(\d (\yy,\ss)). 
 \end{align*}
 Here $\Phi_{Q\sm S_{i-1},\bar \eta}^{*,\de}$  is the unexplored process before step $i$.  Letting $\de \to 0$, we get
 \begin{align*}
 	& \P(\X_i\cap A \ne \X_i'\cap A, \X_j\cap A =\X_j'\cap A,j<i, \EpRCM(A,s)^c \cup \Ed^c\mid \GGc_{i-1}(\Phi_{Q,\es}^{*})^* ) \\
 	& \le 4\E\int_{B_{2s}(A)}\P\Big(\xx\lrsa_{\G(Q,\bar \eta) } \Pd_{Q\sm P} \cup \bar{\eta} \Big) \la(\d \xx) + 2\E \int_{Q \sm B_{2s}(A) }\int_{B_{s}(A)} \pi( \xx, \yy) \la(\d\xx) \la(\d \yy). 
 \end{align*} 
Summing the above and using that $i\le M=4\la(Q)$  under $\Ed$, we obtain the claim. 
\enp

\bep[Proof of Lemma \ref{lem:pa1}]

We may write
\begin{align*}
	\k'(\xx,\vp) &= \one\{\xx\in P_1\sm P_3\} e^{-\sum_{\yy\in \psi_1}v(\xx,\yy)} e^{-\sum_{\yy\in \vp} v(\xx,\yy)\one\{\yy\in P_1\sm P_3\}}
\end{align*}
which has the form \eqref{eq:ktilde} with $\bar{\a}(\xx)=e^{-\sum_{\yy\in \psi_1}v(\xx,\yy)}\one\{\xx\in P_1\}$ and $\bar{v}(\xx,\yy)=v(\xx,\yy)\one\{\xx,\yy\in P_1\}$. 
Moreover,
\begin{align*}
	&\k''(\xx,\vp) =
	e^{-\sum_{\yy\in \psi_1}v(\xx,\yy)}\one\{\xx\in P_{1,-}\sm P_3\} e^{-\sum_{\yy\in \vp} v(\xx,\yy)\one\{\xx,\yy\in P_{1,-}\sm P_3\}}\\
	& + e^{-\sum_{\yy\in \psi_2}v(\xx,\yy)}\one\{\xx\in P_{2,-}\} e^{-\sum_{\yy\in \vp} v(\xx,\yy)\one\{\xx,\yy\in P_{2,-}\}}\\
	&=\Big(e^{-\sum_{\yy\in \psi_1}v(\xx,\yy)}\one\{\xx\in P_{1,-}\sm P_3\} + e^{-\sum_{\yy\in \psi_2}v(\xx,\yy)}\one\{\xx\in P_{2,-}\}\Big) e^{-\sum_{\yy\in \vp}  v(\xx,\yy)\one\{\xx,\yy\in P_{1,-}\sm P_3\text{ or }\xx,\yy\in P_{2,-}\} },
\end{align*}
which again has form \eqref{eq:ktilde}. By construction $\k'$ and $\k''$ satisfy \eqref{eq:tildekk'} with $\a(\xx)=\one\{\xx\notin P_3\}$.
\enp

\subsection{Proofs for Section \ref{ss:ff_RCM}}\label{ss:proffs_ff}

\bep[Proof of Lemma \ref{lem:tilde_iotan}]
In the following, the intersection of a  graph $\GG$ with a set $A$ means intersecting the vertex set with $A$ and taking all induced edges. For the proof, we consider the directed graph $\overrightarrow{\GG}:=\overrightarrow{\GG}(\Phi^*_{U,\psi},\psi)$ with vertex set $\Pd_U\cup \psi$ and directed edges $\xx\to \yy$ if $r_{\xx,\yy}\le \pi(\xx,\yy)$ and $\xx\in \Pd_U$.

Recall the construction of the random graphs $\GG^n:= \GGc (\Phi^*_{ U_n,\psi\cup \nu_n},\psi) $ and $\overrightarrow{\GG}^n:= \overrightarrow{\GG} (\Phi^*_{ U_n,\psi\cup\nu_n},\psi\cup \nu_n) $ from Definition \ref{def:GGd} and Remark \ref{rem:directed}. Let $n'>n$. Then $\overrightarrow{\GG}^n \su \overrightarrow{\GG}^{n'}$ and $\overrightarrow{\GG}^n\cap U_n = \overrightarrow{\GG}^{n'}\cap U_n$. A connected component $\CC^n$ in $\GG^n$ that is not connected to $\nu_n$ is started from a point $\xx_0\in\Pd_{U_n}$ that has no directed path in $\overrightarrow{\GG}^n$ to a $\iota$-smaller point (we ignore edges to $\psi$ in the following because these coincide with edges to $\psi$ in $\overrightarrow{\GG}$ for all $\GG^n$). Then, there is also no directed path from $\xx_0$ to an $\iota$-smaller point  in $\overrightarrow{\GG}^{n'}$ since $\overrightarrow{\GG}^n\cap U_n = \overrightarrow{\GG}^{n'}\cap U_n$ implies that such a path would have to exit $U_n$ at some point, meaning that  $\xx_0$ would be connected to a point in $\nu_n$ in $\GG^n$, which is a contradiction. Thus, $\xx_0$ is also a starting point of a cluster $\CC^{n'}$ in $\GG^{n'}$.  The cluster $\CC^n $ consists of all $\xx>_{\iota} \xx_0$ having a directed path in $\overrightarrow{\GG}^n$ to $\xx_0$ and not to any $\iota$-smaller point. This path would also be a  directed path from $\xx$ to $\xx_0$ in $\overrightarrow{\GG}^{n'}$. If $\xx$ had a path to an $\iota$-smaller point in $\overrightarrow{\GG}^{n'}$, the part of this path until its first exit of $U_n$ would again form a path from $\xx$ to $\nu_n $ in $\GG^n$, which is a contradiction. Hence $\xx\in \CC^{n'}$. On the other hand, if $\xx\in \CC^{n'}$, then it has a directed path to $\xx_0$ and not to any smaller point. This path would lie in $\GG^n$ since otherwise a subpath would connect $\xx$ to a smaller point than $\xx_0$. This shows $\xx\in \CC^n$. In total, we have shown $\CC^n=\CC^{n'}$.  This implies that $\GG^n$ and $\GG^{n'}$ agree on all connected components of $\GG^n$ not connected to $\nu_n$ and thus, we may define $\GGc (\Phi^*_{ U,\psi},\psi) $ as in Lemma \ref{lem:tilde_iotan}.

We need to see that $\GGc (\Phi^*_{ U,\psi},\psi)\stackrel d=\Gamma(U,\psi) $. By Theorem \ref{thm:dc_RCM}, each $\GG^n$ is an RCM $\Gamma(U_n,\psi)$. This can be extended to an RCM $\Gamma^n(U,\psi)\stackrel d=\Gamma(U,\psi) $ on $\Xx$ by adding edges independently between points in $\nu_n$ and from $\nu_n$ to $\psi$. Then, for any $ W_n$ and $m>n$,  
$$ \P\big( \Gamma^m(U,\psi) \cap W_n \ne \GGc (\Phi^*_{ U,\psi},\psi) \cap W_n\big) \le \P\big(\PP_{U_n} \lrsa_{\Gamma^m(U,\psi)} \nu_m \big).$$
Each $\Gamma^m(U,\psi)$ is an RCM so
$$ \P\big(\PP_{U_n} \lrsa_{\Gamma^m(U,\psi)} \nu_m \big) \le \E \sum_{\xx\in \Pd_{U_n}} \one\{\xx \lrsa_{\Gamma_{\xx}(U,\psi)} \PP_{U \sm U_m} \} =\int_{U_n} \P(\xx \lrsa_{\Gamma_{\xx}(U,\psi)} \PP_{U \sm U_m}) \la(\dd\xx) $$
tends to 0 as $m\to \infty$ by assumption \eqref{cond:perc_RCM}. Hence, $\GGc (\Phi^*_{ U,\psi},\psi)$ has the same distribution as each $\Gamma^m(U,\psi)$.
\enp

\bep[Proof of Proposition \ref{prop:RCM_ff}]
Suppose we are on the event $\{\Pd_{A\cap U} \not\lrsa_{\Gamma} \nu_n\}$. 
Let $m\ge n$ and $\xx_0$ be the $\iota$-smallest point in the clusters connected to $A$. By Remark \ref{rem:directed}, $\Tf_{U_m,\nu_m,\psi}$ will first consider all clusters of $\GG^m= \GGc (\Phi^*_{ U_m,\psi\cup \nu_m},\psi) $ having a point $\iota$-smaller than $\xx_0$. From there, the points will be visited in the same order for any $m$. 
Moreover, the retention probabilities  $p^{\ms{RCM}}(\xx,Q_k,\nu_k,B_k,\GGc_k)$ will only  depend on the clusters with no $\iota$-smaller point than $\xx_0$, which are the same for all $m$.  

The limiting Gibbs distribution is first verified for $\Tf_{\Xx,\es}(\Phi^*_{\Xx,\es})$ by checking the GNZ-equation for bounded $f$ of the form $f(\xx,\vp) = f(\xx,\vp\cap W_s)\one\{\xx\in W_s\}$ and then for $\Psi\cup \Tf_{U,\Psi}(\Phi^*_{U,\Psi})$ as in Proposition~\ref{pr:chimera}.
\enp

\bep[Proof of Proposition \ref{prop:RCM2_ff}]
To show \eqref{eq:Tcl2Cmu}, let $\XX_m=T^{\ms{RCM,cl},2}_{U_m,\mu,\nu_m ,\psi\cup \mu_i }(\Phi^*_{U_m,\psi \cup \mu\cup \nu_m})$ and let $n+1<m<m'$. Then,
\begin{align*}
	&\P\Big(\XX_m\cap \CC_\mu \neq \XX_{m'} \cap \CC_\mu,\, \CC_\mu\su W_n \Big)\\
	& \le
	\sum_{k \ge  1} \P(\XX_m\cap V_k \neq\XX_{m'} \cap V_k , \X_m\cap S_{k-1}=\XX_{m'} \cap S_{k-1},S_k\su W_n ,k\le k_0)\\
	& \le
	\sum_{k \ge  1}\E\bigg[\one\{\X_m\cap S_{k-1}=\XX_{m'} \cap S_{k-1},k\le k_0\} \P(\XX_m\cap V_k \neq\XX_{m'} \cap V_k ,S_k\su W_n \mid   (\GG^{\ms{cl},2}_{k-1})^*)\bigg],
\end{align*}
where $k_0$ is the (random) number of steps it takes to explore $\CC_\mu$. 
As long as we are only thinning $\CC_\mu$, it does not matter in which order we explore the later clusters. Hence,
$$\XX_m\cap \CC_\mu=T^{\ms{RCM,cl},2}_{U_m,\mu,\nu_m ,\psi\cup \mu_i }(\Phi^*_{U_m,\psi \cup \mu\cup \nu_m})\cap \CC_\mu =T^{\ms{RCM,cl}}_{U_m,\mu, \psi\cup \mu_i }(\Phi^*_{U_m,\psi \cup \mu\cup \nu_m})\cap \CC_\mu ,$$
and a similar statement holds for $m'$. 
Note that $\X(U_m,\psi\cup \mu_i)$ 
can be considered as a Gibbs process  $\bar{\X}(U_{m'},\psi\cup\mu_i)$ on $U_{m'}$ with a perturbed PI of the form 
$$\bar{\k}(\xx,\vp) = \k\big(\xx, \vp\sm (U_{m'}\sm U_m)  \big)\one\{\xx\in U_m\}.$$
Construct $\bar{\X}(U_{m'},\psi\cup\mu_i)$ by a thinning $\bar T^{\ms{RCM,cl}}_{U_{m'},\mu ,\psi\cup \mu_i }$ as explained in Section \ref{ss:pert_RCM}. When $\CC_\mu\su W_n$,
$$\XX_m\cap \CC_\mu=T^{\ms{RCM,cl}}_{U_m,\mu ,\psi\cup \mu_i }(\Phi^*_{U_m,\psi \cup \mu\cup \nu_m})\cap \CC_\mu =\bar T^{\ms{RCM,cl}}_{U_{m'},\mu ,\psi\cup \mu_i }(\Phi^*_{U_{m'},\psi \cup \mu\cup \nu_{m'}})\cap \CC_\mu $$
since the points of $\CC_\mu$ are visited in the same order by $T^{\ms{RCM,cl}}_{U_m,\mu ,\psi\cup \mu_i }$  and $\bar T^{\ms{RCM,cl}}_{U_{m'},\mu,\psi\cup \mu_i}$ and the retention  probabilities using $\k'$ are the same (indeed, in the formulas for the retention probabilities \eqref{eq:ret_RCM} and \eqref{eq:p_nu}, the partition functions in the numerator and denominator are multiplied by the same constant, namely the probability that the relevant Poisson process has no points in $U_{m'}\sm U_m$).  Set $\bar{\X}_m = \bar T^{\ms{RCM,cl}}_{U_{m'},\mu,\psi\cup \mu_i }(\Phi^*_{U_{m'},\psi \cup \mu\cup \nu_{m'}})$.
Then, 
\begin{align}\nonumber
	&	\sum_{k \ge  1}\E\big[\one\{\X_m\cap S_{k-1}=\XX_{m'} \cap S_{k-1},k\le k_0\} \P(\XX_m\cap V_k \neq\XX_{m'} \cap V_k ,S_k\su W_n ,\mid   (\GG^{\ms{cl},2}_{k-1})^*)\big]\\ \nonumber
	& \le	\sum_{k \ge  1}\E\big[\one\{\bar{\X}_m\cap S_{k-1}=\XX_{m'} \cap S_{k-1},k\le k_0,S_{k-1}\su W_n\} \P(\bar{\XX}_m\cap V_k \neq {\XX}_{m'} \cap V_k \mid   (\GG_{k-1}^{\ms{cl}})^*)\big]\\\nonumber
	&\le \sum_{k \ge  1}	\E\bigg[\one\{k\le k_0, S_{k-1}\su W_n\}\\ 
	&\quad \times\E\bigg[ 2\int_{U_{m'}}  \one\{ \xx \sim_{\Gamma_{k-1,\xx}} V_{k-1}\cap W_n \}\one \{\xx\lrsa_{\Gamma_{k-1,\xx}}  \nu_m\} \la(\d\xx)\mid   (\GG_{k-1}^{\ms{cl}})^*\bigg] \bigg] ,\label{eq:compu}
\end{align}
where $\Gamma_{k-1,\xx}$ denotes the graph $\Gamma_{k-1}$ from \eqref{eq:Gammak} under the Palm distribution at $\xx$ of the underlying Poisson process. 
To see where the second inequality came from, we consider the $\de$-approximations as in the proof of Theorem \ref{thm:Talpha}. Let $V_1^\de = \{(\xx,\rr)\in U_{m'}^\de: (\xx,\rr)\sim_\de \mu \}$ and $V_i^\de = \{(\xx,\rr)\in U_{m'}^\de: (\xx,\rr)\sim_\de V_{i-1}^\de\}$ and let $\X^\de_{m,i}$ and $\bar \X_{m,i}^\de$ be the part of the Gibbs processes constructed in step $i$. Conditionally on the points explored up to step $k-1$, i.e.\ $\Phi_{U_{m'},\psi\cup \mu}^{*} \cap S_{k-1}$, if the thinnings agreed so far, we have from Corollary \ref{cor:poiss_emb} that
\begin{align*}
	& \one\{ \X_{m,i} \ne \bar \X_{m,i} ,i\le k-1 \}\P\big(\X_k^\de \ne \bar \X_k^\de  \mid \Phi_{U_{m'},\psi\cup \mu}^{*} \cap S_{k-1} \big)\\
	& \le  2\int_{V_k^\de} \P\Big(\xx\lrsa_{\{\xx\}\cup \Phi_{U_{m'}, \psi\cup \mu}^{*,\de}} \Phi_{U_{m'},\psi\cup\mu}^{*,\de}\cap (U^\de_{m'}\sm U_m^\de)\Big) \big(\la\otimes \Leb_{[0,1]}^\de\big)(\d\xx).
\end{align*}
Letting $\de\to 0$, we obtain 
\begin{align*}
	&\one\{\bar{\X}_m\cap S_{k-1}=\XX_{m'} \cap S_{k-1}\} \P(\bar{\XX}_m\cap V_k \neq {\XX}_{m'} \cap V_k \mid   (\GG_{k-1}^{\ms{cl}})^*)
	\le  2\int_{U_{m'}} \P(\xx\sim V_{k-1} , \xx\lrsa_{\G_{k-1,\xx}} \Pd_{U\sm U_m} )\la(\d\xx),
\end{align*}
which is what was used in \eqref{eq:compu}.

Note that only the conditional distribution of $\Gamma_{k-1}=\GGd_{k-1}\cup \GGG_{V_k}\cup \GGB_{V_k}$ given $(\GGc_{k-1})^*$ matters in \eqref{eq:compu}, and this is the same as the distribution of $\GG^{\ms{cl},2}$. 
Splitting the integral after whether $\xx \in U_{m-1}$ or $\xx\in U_{m'}\sm U_{m-1} $, we obtain the bound on \eqref{eq:compu}
\begin{align*}
	&	2\sum_{k \ge 1}	\E\bigg[ \one\{k\le k_0\}  \E\bigg[\int_{U_{m'}\sm U_{m-1}}  \one\{\xx \notin S_{k-1,\xx}, \xx \sim_{\GG_{\xx}^{\ms{cl},2}} V_{k-1,\xx}\cap W_n \} \la(\d\xx)\mid (\GG^{\ms{cl},2}_{k-1})^*  \bigg]\bigg]  \\ 
	& \quad + 2  \sum_{k \ge 1}	\E\bigg[ \one\{k\le k_0\}  \E\bigg[ \int_{U_{m-1}} \one\{\xx \notin S_{k-1,\xx}, \xx \sim_{\GG_{\xx}^{\ms{cl},2}} V_{k-1,\xx} \}\one \{\xx\lrsa_{\GG^{\ms{cl},2}_{\xx}}  \nu_m\} \la(\d\xx)\mid (\GG^{\ms{cl},2}_{k-1})^* \bigg]\bigg]\\
	&=	2\E\bigg[ \sum_{k = 1}^{k_0}	\int_{U_{m'}\sm U_{m-1}}  \one\{\xx \notin S_{k-1,\xx}, \xx \sim_{\GG_{\xx}^{\ms{cl},2}} V_{k-1,\xx}\cap W_n \} \la(\d\xx)\bigg]  \\ 
	& \quad + 2 \E\bigg[ \sum_{k = 1}^{k_0}	  \int_{U_{m-1}} \one\{\xx \notin S_{k-1,\xx}, \xx \sim_{\GG_{\xx}^{\ms{cl},2}} V_{k-1,\xx} \}\one \{\xx\lrsa_{\GG^{\ms{cl},2}_{\xx}}  \nu_m\} \la(\d\xx)\bigg]\\
	& \le 2	\E\Big[ \int_{U_{m'}\sm U_{m-1}} \one\{\xx\sim_{\GG^{\ms{cl},2}_{\xx}} \mu\cup \Pd_{W_n} \} \la(\d\xx)\Big] + 2	\E\Big[ \int_{U_{m-1}}  \one \{\xx\lrsa_{\GG_{\xx}^{\ms{cl},2}}  \nu_m\} \la(\d\xx)\Big]. 
\end{align*}
In the inequality, we used that the events  $\{\xx \notin S_{k-1,\xx}, \xx \sim_{\GG^{\ms{cl},2}_{\xx}} V_{k-1,\xx}\cap W_n \}\su \{\xx\in V_{k,\xx}\}$ are disjoint.   We further get the bound
\begin{align*}
	&	 2	\E  \int_{U\sm U_{m-1}} \sum_{\yy\in \mu}\pi(\xx,\yy)\la(\d\xx) 
	+ 2 \E\int_{U \sm U_{m-1}}\sum_{\yy\in \Pd_{W_n}}\pi(\xx,\yy)   \la(\d \yy) \la(\d\xx) 
	+ 2	\E\Big[ \int_{U_{m-1}}  \one \{\xx\lrsa_{\Gamma_{\xx}}  \nu_m\} \la(\d\xx)\Big]\\
	&	\le 2	c \int_{U\sm U_{m-1}} \int_{W_1}\pi(\xx,\yy)\la(\d \yy)\la(\d\xx) 
	+ 2 \int_{U \sm U_{m-1}}\int_{W_n} \pi(\xx,\yy)   \la(\d \yy) \la(\d\xx) 
	+ 2	\int_{U_{m-1}}  \P(\xx\lrsa_{\Gamma_{\xx}}  \nu_m) \la(\d\xx).
\end{align*}
Taking $m=n+2$ and $m'=n+3$, we get
\begin{align*}
	&\P\big(\XX_{n+2}\cap \CC_\mu \neq \XX_{n+3} \cap \CC_\mu\big)\\ 
	&\le 2	c \int_{\Xx\sm W_{n+1}} \int_{W_1}\pi(\xx,\yy)\la(\d \yy)\la(\d\xx) 
	+ 2 \int_{\Xx \sm W_{n+1}}\int_{W_n} \pi(\xx,\yy)   \la(\d \yy) \la(\d\xx) 
	\\&\quad + 2	\int_{W_{n+1}}  \P(\xx\lrsa_{\Gamma_{\xx}}  \Pd_{\Xx\sm W_{n+2}}) \la(\d\xx) + \P(\CC_\mu \cap (\Xx\sm W_n)\ne \es)\\
	& \le 2c\cdot 2^{-n } + 2\cdot 2^{-n} + 2\cdot 2^{-n-1} + c\cdot 2^{-n+1}.
\end{align*}
These probabilities are summable, hence by the Borel-Cantelli lemma, there is an $N(\omega)$ s.t. $\X_n\cap \CC_\mu=\X_N\cap \CC_\mu$ for all $n\ge N$. This shows \eqref{eq:Tcl2Cmu}.

To show \eqref{eq:Tcl2A}, note that we first construct $\CC_\mu$ as the set of points in $\Pd_U$ having a directed path to $\mu$ in $\overrightarrow{\GG}(\Psi_{U_m,\mu\cup \psi\cup \nu_m})$. The clusters of the remaining points are constructed as in the proof of Lemma \ref{lem:tilde_iotan}. Hence, the argument that the thinnings agree outside $\CC_\mu$ from the proof of Proposition \ref{prop:RCM_ff} carries over to show the second statement. Since the boundary conditions $\mu_i$ only affect the thinning of $\CC_\mu$, \eqref{eq:Tcl2_disag} follows immediately.

The limiting Gibbs distribution in \eqref{eq:Tcl_gibbs} is verified by checking the GNZ equations  for $\Psi \cup T^{\ms{RCM},\ff,2}_{U,\mu,\Psi}(\Phi^*_{U,\Psi\cup \mu})$, noting that $T^{\ms{RCM},\ff,2}_{U,\mu,\Psi}(\Phi^*_{U,\Psi\cup \mu}) $ equals $T^{\ms{RCM,cl},2}_{U_m,\mu,\nu_m,\Psi}(\Phi^*_{U_m,\Psi\cup \mu \cup \nu_m})$ for  $m\ge M$ for some $M(\omega)$.
\enp


 \section{Disagreement coupling proofs for Section \ref{sec:malliavin_stein}}\label{sec:proofthm2}

The main challenge in proving Lemma \ref{lem:prprpr} is that when we modify the Poisson process we are thinning, then the unexplored part of the domain in a given step of the cluster-based thinning algorithm may look different. Thus, a central ingredient in the proof is the following  lemma that bounds the effect of changing the domain.
 
 \bel\label{lem:Delta}
 Consider a  PI $\k$ from Model (L) on a set $P$. Define two PIs   $\k_i(\xx,\vp)=\k(\xx,\vp)\one\{\xx\in P_i\}$  where $P_i\su P$, $i=1,2$. Denote the corresponding Poisson embeddings by $T^i_{P,\psi}(\Pds_P)$. Let $M>0$.
 Then, 
 $$\P(T^1_{P,\psi}(\Pds_P)\ne T^2_{P,\psi}(\Pds_P )) \le M\big(2 + \la(P)\big)\la(P_1\De P_2)+\P( \PP(P) > M).$$
 \enl

 \bep
 From Lemma \ref{lem:distv}, we get that
 $$ \P(\inf_\iota T^1_{P,\psi}(\Pds_P)\cap P \ne \inf_\iota T^2_{P,\psi}(\Pds_P)\cap P)\le \big(2 + \la(P)\big)\dtv(T^1_{P,\psi}(\Pds_P)\cap P , T^2_{P,\psi}(\Pds_P)\cap P) .$$ 
 We can construct Gibbs processes with the same distribution as $T^1_{P,\psi}(\Pds_P)$ and $T^2_{P,\psi}(\Pds_P)$ by a Poisson embedding using an ordering where $P_1\De P_2$ comes before $P\sm (P_1\De P_2)$. If $\PP_{P_1\De P_2}=\es$, the PIs look the same on $P\sm (P_1\De P_2)$ and have the same boundary conditions. Hence, the Poisson embeddings will agree. Thus,
 $$\P\big(\inf_\iota T^1_{P,\psi}(\Pds_P)\cap P \ne \inf_\iota T^2_{P,\psi}(\Pds_P)\cap P\big) \le \big(2 + \la(P)\big)\P(\PP_{P_1\De P_2}\ne \es) \le \big(2 + \la(P)\big)\la(P_1\De P_2).$$
 Using this argument repeatedly when $\PP(P) \le M$, we get
 $$ \P\big( T^1_{P,\psi}(\Pds_P) \ne  T^2_{P,\psi}(\Pds_P)\big) \le M \big(2 + \la(P)\big)\la(P_1\De P_2) + \P\big(\PP(P) > M\big).$$ 
 \enp

Now, we have all the ingredients to prove Lemma \ref{lem:prprpr}. Again, $C>0$ denotes a constant which may change value from line to line.  To shorten notation, we write $Q$ for $Q_n^+$.
 
   \bep[Proof of Lemma \ref{lem:prprpr}] 
 	We show the inequality for $\P(\XX_{\xx}''' \ne \check \XX_{\xx}''')$. The remaining inequalities are similar.  We first consider Model (L), which is conceptually simpler, and then sketch the generalization to Model (P).
	\medskip

\begin{figure}
\begin{tikzpicture}[scale=1.0, every node/.style={font=\small}]
\tikzset{
  ptP/.style={circle, fill=black, inner sep=1.6pt},
  ptC/.style={circle, fill=green, inner sep=1.6pt},
  Yring/.style={circle, draw=black, line width=0.9pt, inner sep=2.8pt},
  dashshared/.style={dashed, line width=1.0pt},
  dashY/.style={dashed, red, line width=0.9pt},
  box/.style={draw=black, line width=1pt},
}

\def\W{10.2}
\def\H{6.2}

\def\xb{5.1}
\def\yb{3.1}
\def\r{2.2}

\draw[box] (0,0) rectangle (\W,\H);
\node[anchor=north west] at (0.3,\H-0.3) {$Q$};

\draw[dashed,orange] (\xb,\yb) circle (\r);
\node[anchor=west,orange] at (6,3.1) {$B_{2r}(\xx)$};

\draw[ blue] (6.0,4.1) -- (5.3,4.2);
\draw[ blue] (5.3,4.2) -- (4.8,3.6);
\draw[ blue] (5.1,1.8) -- (4.3,1.7)-- (4.6,2.2) -- (5.1,1.8) --(5.7,2.4);

\draw (1.0,5.2) -- (1.6,4.4);
\draw (2.6,5.6)--(3.2,4.8);
\draw (1.3,1.0) -- (2.1,1.6) -- (2.9,0.9);
\draw (8.8,4.6) -- (9.2,3.7);
\draw (8.5,1.6) -- (9.1,2.1);

\draw[] (4.2,4.9) -- (3.2,4.8); 

\draw[ green] (2.9,2.4) -- (3.3,3.1); 
\draw[ green] (2.4,3.3) -- (3.3,3.1); 
\draw[green] (6.6,1.9) -- (7.4,2.4); 
\draw[ black] (6.6,1.9) -- (7.4,1.3);
\draw[green] (3.3,3.1) -- (3.9,3.1);
\draw[black,dashed] (3.3,3.1) -- (3.9,3.1);


\draw[green]  (0.8,3.4)--(1.4,2.7); 
\draw[green] (2.4,3.3) -- (2.9,2.4);
\draw[green] (8.4,3.5) -- (8.1,2.8);
\draw[green] (7.4,2.4) -- (8.1,2.8);

\foreach \p in {
  (4.8,3.6),(5.3,4.2),
  (4.6,2.2),(5.1,1.8),(5.7,2.4),(6.0,4.1),
  (4.3,1.7),
   	 (3.3,3.1),
 	 (3.9,3.1),
 (4.2,4.9),
 (6.6,1.9)
}{
  \node[ptC] at \p {};
  \node[ptP] at \p {};
}

\foreach \p in {
	 (3.3,3.1),
	(4.2,4.9),
	(6.6,1.9),
		 (3.9,3.1)
}{
  \node[Yring] at \p {};
}

\foreach \p in {
  (1.0,5.2),(1.6,4.4),(2.6,5.6),(3.2,4.8),
  (7.9,5.4),(8.8,4.6),(9.2,3.7),
  (7.4,1.3),(8.5,1.6),(9.1,2.1),
  (1.3,1.0),(2.1,1.6),(2.9,0.9)
}{
  \node[ptP] at \p {};
  \node[Yring] at \p {};
}

\foreach \p in {
  (0.8,3.4),(1.4,2.7),(2.4,3.3),(2.9,2.4),
  (6.8,5.3),(8.4,3.5),(9.4,4.4),
  (7.4,2.4),(8.1,2.8),(9.1,1.2)
}{
  \node[ptC] at \p {};
  \node[Yring] at \p {};
}


\draw[green]  (0.8,3.4)--(1.4,2.7); 
\draw[green] (2.4,3.3) -- (2.9,2.4);
\draw[green] (8.4,3.5) -- (8.1,2.8);
\draw[green] (7.4,2.4) -- (8.1,2.8);

\node[circle, fill=orange, inner sep=1.8pt] at (\xb,\yb) {};

\node[orange, anchor=north] at (\xb,\yb) {$\xx$};

\end{tikzpicture}
\caption{
Coupled configurations $\PP$ and $\check\PP$ in the domain $Q$.
Black  points belong to $\PP$. Inside $B_{2r}(\xx)$, the two processes agree, while  green points represent $\check\PP \sm B_{2r}(\xx)$.   Circles indicate members of the set $\YY$, that is, points connected to $Q\setminus B_{2r}(\xx)$ by a path in either $\PP$ or $\check\PP$.
The clusters in $\PP$ intersecting $Q\sm B_{2r}(\xx)$ are shown with black edges, while clusters in $\check \PP$ intersecting $Q\sm B_{2r}(\xx)$. The dashed  black and green edge represents two points that form their own cluster in $\PP$ but are part of a larger cluster in $\check \PP$. The blue clusters inside $B_{2r}(\xx)$ are common for the two processes. These are thinned in the same order, but different black or green clusters may be visited in between. 
}
\label{fig:qpp}
\end{figure}
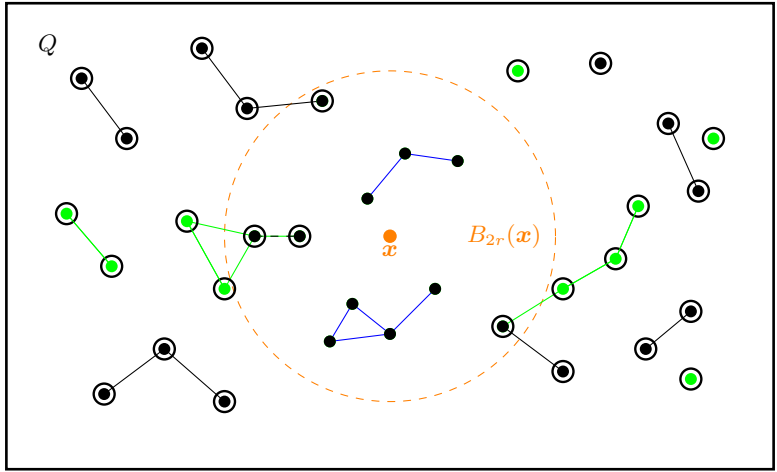

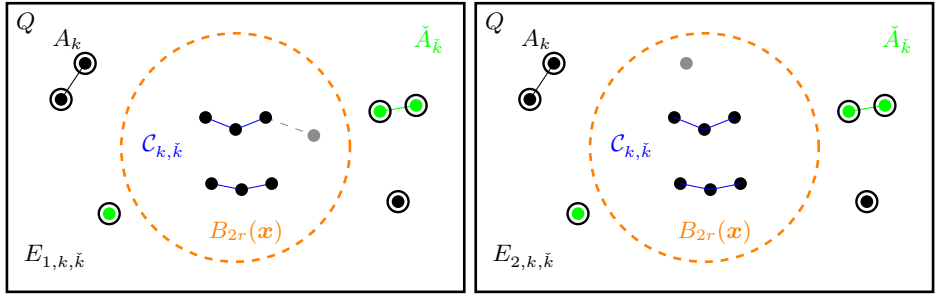
\begin{figure}
\begin{tikzpicture}[scale=0.8, every node/.style={font=\small}]
\tikzset{
  ptP/.style={circle, fill=black, inner sep=1.7pt},
  ptC/.style={circle, fill=black!45, inner sep=1.7pt},
  edgeP/.style={dashed,color=black!45},
  edgeC/.style={ color=blue},
  extraP/.style={circle, fill=black, line width=0.9pt, inner sep=1.7pt},
  extraC/.style={circle, fill=green, line width=0.9pt, inner sep=1.7pt},
  extra/.style={circle, draw=red, line width=0.9pt, inner sep=2.3pt},
  reject/.style={circle, draw=black!45, inner sep=2.2pt},
  ball/.style={dashed, line width=1pt},
  box/.style={draw=black, line width=1pt},
    Yring/.style={circle, draw=black, line width=0.9pt, inner sep=2.8pt},
}

\draw[box] (0,0) rectangle (7.6,4.8);
\node[anchor=north west] at (0,4.8) {$Q$};

\draw[ball, orange] (3.8,2.4) circle (1.9);
\node[anchor=west,orange] at (3.2,1) {$B_{2r}(\xx)$};

\coordinate (a) at (3.3,2.9);
\coordinate (b) at (3.8,2.7);
\coordinate (c) at (4.3,2.9);
\coordinate (d) at (3.4,1.8);
\coordinate (e) at (3.9,1.7);
\coordinate (f) at (4.4,1.8);
\coordinate (v) at (5.1,2.6); 

\draw[edgeC] (a)--(b)--(c);
\draw[edgeC] (d)--(e)--(f);
\draw[edgeP] (c)--(v);

\foreach \p in {a,b,c,d,e,f} {
  \node[ptP] at (\p) {};
}

\node[blue] at (2.6,2.4) {$ \CC_{k,\check k}$};

\node[extraP] at (6.5,1.5) {};
\node[extraP] at (1.3,3.8) {};
\node[extraP] at (0.9,3.2) {};
\node[Yring] at (0.9,3.2) {};
\node[Yring] at (1.3,3.8) {};
\node[Yring] at (6.5,1.5) {};
\draw (0.9,3.2) -- (1.3,3.8);

\node[anchor=west] at (0.6,4.2) {$A_k$};

\draw[green](6.8,3.1) -- (6.2,3.0);
\node[extraC] at (6.2,3.0) {};
\node[extraC] at (6.8,3.1) {};
\node[extraC] at (1.7,1.3) {};
\node[Yring] at (6.2,3.0) {};
\node[Yring] at (1.7,1.3) {};
\node[Yring] at (6.8,3.1) {};

\node[anchor=west, green] at (6.6,4.2) {$\check A_{\check k}$};

\coordinate (v) at (5.1,2.6);

\node[ptC] at (v) {};


\node[anchor=south] at (0.8,0.2)
{\(E_{1,k,\check k}\)};
\end{tikzpicture}
\begin{tikzpicture}[scale=0.8, every node/.style={font=\small}]
	\tikzset{
		ptP/.style={circle, fill=black, inner sep=1.7pt},
		ptC/.style={circle, fill=black!45, inner sep=1.7pt},
		edgeP/.style={dashed,color=blue},
		edgeC/.style={ color=blue},
		extraP/.style={circle, fill=black, line width=0.9pt, inner sep=1.7pt},
		extraC/.style={circle, fill=green, line width=0.9pt, inner sep=1.7pt},
		extra/.style={circle, draw=red, line width=0.9pt, inner sep=2.3pt},
		reject/.style={circle, draw=black!45, inner sep=2.2pt},
		ball/.style={dashed, line width=1pt},
		box/.style={draw=black, line width=1pt},
		Yring/.style={circle, draw=black, line width=0.9pt, inner sep=2.8pt},
	}
	
	\draw[box] (0,0) rectangle (7.6,4.8);
	\node[anchor=north west] at (0,4.8) {$Q$};
	
	\draw[ball, orange] (3.8,2.4) circle (1.9);
	\node[anchor=west,orange] at (3.2,1) {$B_{2r}(\xx)$};
	
	\coordinate (a) at (3.3,2.9);
	\coordinate (b) at (3.8,2.7);
	\coordinate (c) at (4.3,2.9);
	\coordinate (d) at (3.4,1.8);
	\coordinate (e) at (3.9,1.7);
	\coordinate (f) at (4.4,1.8);
	
	\foreach \p in {a,b,c,d,e,f} {
		\node[ptP] at (\p) {};
	}

	\draw[edgeC] (a)--(b)--(c);
	\draw[edgeC] (d)--(e)--(f);
	
	\node[blue] at (2.6,2.4) {$ \CC_{k,\check k}$};

	\node[extraP] at (6.5,1.5) {};
	\node[extraP] at (1.3,3.8) {};
	\node[extraP] at (0.9,3.2) {};
	\node[Yring] at (0.9,3.2) {};
	\node[Yring] at (1.3,3.8) {};
	\node[Yring] at (6.5,1.5) {};
	\draw (0.9,3.2) -- (1.3,3.8);
	
	\node[anchor=west] at (0.6,4.2) {$A_k$};
	
	\draw[green](6.8,3.1) -- (6.2,3.0);
	\node[extraC] at (6.2,3.0) {};
	\node[extraC] at (6.8,3.1) {};
	\node[extraC] at (1.7,1.3) {};
	\node[Yring] at (6.2,3.0) {};
	\node[Yring] at (1.7,1.3) {};
	\node[Yring] at (6.8,3.1) {};
	
	\node[anchor=west, green] at (6.6,4.2) {$\check A_{\check k}$};

	\coordinate (v) at (3.5,3.8);
	
	\node[ptC] at (v) {};
	
	
	\node[anchor=south] at (0.8,0.2)
	{\(E_{2,k,\check k}\)};
\end{tikzpicture}

\caption{Illustration of events $E_{1,k,\check k}$ (left) and $E_{2,k,\check k}$ (right).
In both panels, black points belong to $\PP$, while green points belong to $\check\PP \sm \PP$.
The orange circle indicates the ball $B_{2r}(\xx)$. Up to step $k$ and $\check k$, respectively, we have explored the common blue clusters ($\CC_{k,\check k})$. Additionally, we have explored the black clusters ($A_k$) in $\PP $ and the green clusters ($\check A_{\check k}$) in $\check \PP$. Under $E_{1,k,\check k}$, we continue exploring the same blue cluster via the dashed path. Under $E_{2,k,\check k}$, the blue cluster has been completed and we start a new cluster from the grey point. 
}
\label{fig:E1E2}
\end{figure}

{\bf Model (L).} 
 	Let $$\YY=\{ \ww \in \PP\cup \check{\PP} \mid \ww \lrsa_{\PP } \PP_{Q\sm B_{2r}(\xx)} \text{ or } \ww \lrsa_{ \check \PP } \check\PP_{Q\sm B_{2r}(\xx)}\}.$$ 
 	Recall that $\PP$ and $\check \PP$ agree on $B_{2r}(\xx)$, while $\PP\cup \check \PP$ forms a Poisson process of double intensity  outside this ball. We refer to Figure \ref{fig:qpp} for an illustration.
 	The clusters of $\PP\sm \YY \su B_{2r}(\xx)$ correspond  to clusters in both $\PP$ and $\check \PP$ that do not contain a point in $Q\sm B_{2r}(\xx)$. 	The exploration algorithms for constructing $\XX_{\xx}''' $ and $ \check \XX_{\xx}'''$ visit the points of $\PP\sm \YY = \check\PP\sm \YY$ in the same order, but different points in $\YY$ may be visited in between. If the thinning of $\PP$ explores points of $\PP\sm \YY$ in step $k+1$, there is a $\check k \ge 1$ such that the thinning of $\check \PP$ considers the same points in step $\check k + 1$. 
	For this to happen, we must be under one of the following  two events, which are illustrated in Figure \ref{fig:E1E2}:
	\begin{enumerate}
	\item $E_{1,k,\check k} = E_{1, k}^{(1)} \cap E_{1, k, \check k}^{(2)} \cap E_{1, k, \check k}^{(3)}$, where 
	\begin{align*}
	E_{1, k}^{(1)} =& \big\{\text{in step } k+1 \text{ of the exploration of } \PP, \text{ we continue exploring a cluster contained in } B_{2r}(\xx)\};\\
	E_{1, k, \check k}^{(2)} =& \big\{\text{in step } \check k + 1 \text{ of the exploration of } \check \PP, \text{ we continue exploring a cluster which up to step $\check k$ looks}\\
	& \text{the same as the one considered up to step $k$ in } \PP\big\};\\
	E_{1, k, \check k}^{(3)} =& \big\{\text{the points explored in the next step agree, i.e. } \PP_{S_{k+1}\sm S_k}=\check{\PP}_{\check S_{\check k+1}\sm \check S_{\check k}} \big\}.
	\end{align*}
	\item $E_{2,k,\check k} = E_{2, k}^{(1)} \cap E_{2,  \check k}^{(2)} \cap E_{2, k, \check k}^{(3)}$, where
	\begin{align*}
	E_{2, k}^{(1)} =& \big\{\text{in step } k+1 \text{ of the exploration of } \PP, \text{ we start a new cluster}\};\\
	E_{2, \check k}^{(2)} =& \big\{\text{in step } \check k + 1 \text{ of the			 exploration of } \check \PP, \text{ we start a new cluster}\big\};\\
	E_{2, k, \check k}^{(3)} =& \big\{\text{the smallest unexplored Poisson point of the exploration at step $k$ + 1 in } \PP \\
	& \text{ and of the exploration at step } \check k + 1 \text{ in } \check \PP \text{ agree}\big\}.
	\end{align*}
	\end{enumerate} 


 	Let $\XX'''_{\xx,\le i}$ and $\check \XX'''_{\xx,\le i}$ be the parts constructed before step $i+1$, and $\XX'''_{\xx,i+1}$ and $\check \XX'''_{\xx,i+1}$ be the part of the processes constructed in step $i+1$. 
 	Then, 
 	\begin{align}\nonumber
 		& \P\big(\XX'''_{\xx} \ne \check\XX'''_{\xx}\big)
 		\nonumber
 		\le 
 		\E\bigg[\sum_{k=1}^{4\la(Q)}\sum_{\check k =1}^{4\la(Q)} \P\big(\XX'''_{\xx,{k+1}}\ne \check\XX'''_{\xx,\check k +1},  \XX'''_{\xx,\le k}= \check\XX'''_{\xx,\le \check k}, E_{1,k,\check k}\cup E_{2,k,\check k} \mid \PP_{S_k}^*\cup \check \PP^*_{\check S_{\check k}}\big) \bigg]\\ \nonumber
 		& \quad + 2\P(\Ed)+\P(\YY\cap B_{r+3s}(\xx)\ne \es)\\ \nonumber
 		&\label{eq:P'''}
 		\le 16\la(Q)^2 
 		\sup_{k,\check k\ge 0}\E\big[\P\big(\XX'''_{\xx,{k+1}}\ne \check\XX'''_{\xx,\check k +1},  \XX'''_{\xx,\le k}= \check\XX'''_{\xx,\le \check k}, E_{1,k,\check k}\cup E_{2,k,\check k} \mid \PP_{S_k}^*\cup \check \PP^*_{\check S_{\check k}}\big) \big] \\
 		&\quad + 2\P(\Ed)+\P(\YY\cap B_{r+3s}(\xx)\ne \es). 
 	\end{align}
 	We bound the expectation under each of the two events $E_{1, k ,\check k}$ and $E_{2, k ,\check k}$ separately.
	\medskip
 	

\noindent {\bf Event $E_{1,k,\check k}$.} 	For $E_{1,k,\check k}$,  the so far explored points up to step $k$ and $\check k$, respectively, consists of the union $\CC_{k,\check k}$ of some clusters  in $\PP_{B_{2r}(\xx)}$ that agree for both thinnings (including the current cluster) and some extra points, which we denote by $A_k\su \PP$ and $\check A_{\check k}\su \check \PP$, respectively. See Figure \ref{fig:E1E2} for illustration. Note that $\CC_{k,\check k} \cap \YY =\es$ and $A_k\cup \check A_{\check k}\su \YY$. Recall the neighborhood notation $N(B)$ from Section~\ref{ss:Gibbs}. The unexplored window before step $k+1$ is $Q^k = Q_{(\xx_0,\ff)}\sm N(\CC_{k-1,\check k -1}\cup A_k)$ and the unexplored window before step $\check k +1 $ is  $\check Q^{\check k} =Q_{(\xx_0,\ff)}\sm N(\CC_{k-1,\check k -1}\cup\check A_{\check k})$, respectively, where $\xx_0$ is the starting point of the current cluster. 

 	Note that  $S_k\cup \check S_{\check k} \cup ( Q\sm B_{2r}(\xx))$ is a stopping set (strictly speaking, $S_k \ti \{1\}\cup \check S_{\check k} \ti \{1\}\cup ( Q\sm B_{2r}(\xx))\ti\{1,2\}$ is a stopping set with respect to the Poisson process $\PP\ti \{1\} \cup \check \PP_{Q\sm B_{2r}(\xx)}\ti \{2\}$ on $Q\ti \{1,2\}$). Hence $\PP_{B_{2r}(\xx)\sm (S_k \cup \check S_{\check k})} $ is a Poisson process on $B_{2r}(\xx)\sm ({S_k}\cup \check S_{\check k})$ conditionally on $\sigma(\PP_{S_k},\check \PP_{\check S_{\check k}})$.  	  
 	We further have $B_{2r}(\xx)\sm ({S_k}\cup \check S_{\check k}) = (Q_{(\xx_0,\ff)}\cap B_{2r}(\xx))\sm N(\CC_{k-1,\check k -1}\cup A_k\cup \check A_{\check k}) $.
 	 Take a Poisson process $\PP''$ independent of $\s (\PP_{S_k}, \check\PP_{\check S_{\check k}})$ on $ N(A_k\cup  \check A_{\check k}) \cap \big((Q_{(\xx_0,\ff)} \cap B_{r+3s}(\xx))\sm N(\CC_{k-1,\check k -1})\big)$    and let $\PP'= \PP_{B_{r+3s}(\xx)\sm (S_k \cup \check S_{\check k})} \cup \PP''$. Thus, $\PP'$ is the union of two independent Poisson processes on the disjoint sets $$B_{r+3s}(\xx)\sm (S_k \cup \check S_{\check k})=(Q_{(\xx_0,\ff)}\cap B_{r+3s}(\xx))\sm N(\CC_{k-1,\check k -1}\cup A_k\cup \check A_{\check k})$$ and $$N(A_k\cup  \check A_{\check k}) \cap \big((Q_{(\xx_0,\ff)} \cap B_{r+3s}(\xx))\sm N(\CC_{k-1,\check k -1})\big)$$ whose union is $(Q_{(\xx_0,\ff)} \cap B_{r+3s}(\xx)) \sm N(\CC_{k-1,\check k -1})$ and we may think of $\PP_{B_{r+3s}(\xx)\sm (S_k \cup \check S_{\check k})}$ as a restriction of $\PP'$.  
 	
 When a PI on $Q$ has the form $\k(\ww,\vp) = \k(\ww,\vp\cap B')\one_B$ for some set $B \su B'\su Q$, the formula for the retention  probabilities shows that
 	$$p\big(\ww,Q_{(\ww,\ff)},\phi\big) = p\big(\ww,B_{(\ww,\ff)},\phi\cap B'\big)\mathds{1}\{\ww\in B\}.$$
 	Thus, the Poisson embedding $T_{Q,\psi}(\vp^*_Q)$  equals $T_{B,\psi\cap B'}(\vp^*_B)$.
 	Taking $B= B' = B_{r+3s}(\xx)$, this means that the part of the Gibbs process generated in step ${k+1}$ by a Poisson embedding on $Q^k$ equals the Poisson embedding 
	$$\XX_{k+1}^{\prime\prime\prime}=T_{Q^k\cap B_{r+3s}(\xx),\phi_k\cap B_{r+3s}(\xx)}\big( \PP_{Q^k\cap B_{r+3s}(\xx)}^*\big) \cap S_{k+1},$$ where $\phi_k$ are the boundary conditions for this step. 
 	Moreover, $\PP_{Q^k\cap B_{r+3s}(\xx)} \cap S_{k+1} = \PP_{ B_{r+3s}(\xx)\sm S_k}\cap S_{k+1}= \PP_{ B_{r+3s}(\xx)\sm (S_k\cup \check S_{\check k})}\cap S_{k+1} $ since $(\PP_{\check S_{\check k}} \sm \PP_{ S_{ k}})\cap B_{r+3s}(\xx) \su \check A_{\check k} $ and $\PP_{S_{k+1}\sm S_k} \cap \check A_{\check k} = \PP_{\check S_{\check k+1}\sm \check S_{\check k}} \cap \check A_{\check k} = \es$ under $E_{1,k,\check k}$.   
 	On the event $\PP''=\es$, we thus have $\PP_{Q^k\cap B_{r+3s}(\xx)}\cap S_{k+1} = \PP' \cap S_{k+1}$ and hence
 	\begin{align*}
 		\XX'''_{\xx,k+1} 
 		&= T_{(Q_{(\xx_0,\ff)}\sm N(\CC_{k,\check k})) \cap B_{r+3s}(\xx),\phi_k\cap B_{r+3s}(\xx)}^{(4)}\big( (\PP')^* \big) \cap S_{k+1}. 
 	\end{align*}
 	Here, the  thinning $T^{(4)}$ used the PI $\k^{(4)}_\xx(\ww, \vp) = \k'''_\xx(\ww, \vp)\one \big\{ \ww \in Q^k\cap B_{r+3s}(\xx) \big\}$.
 	By similar arguments, if $\PP''=\es$, then
 	\begin{align*}\check \XX'''_{\xx,k+1}
 		= {T}_{(Q_{(\xx_0,\ff)}\sm N(\CC_{k,\check k})) \cap B_{r+3s}(\xx),\check \phi_{\check k}\cap B_{r+3s}(\xx)}^{(5)}\big( (\PP')^*\big) \cap \check S_{\check k+1}, 
 	\end{align*}
 	where the thinning $ T^{(5)}$ used the PI $\check \k^{(5)}_\xx(\ww, \vp) = \k'''_\xx(\ww, \vp)\one \big\{ \ww \in \check Q^{\check k} \cap B_{r+3s}(\xx)\big\}$. If the boundary conditions agree, i.e., $\phi_{ k} \cap B_{r+3s}(\xx) = \check \phi_{\check k} \cap B_{r+3s}(\xx)$, the only difference is that the  PIs $\k^{(4)}$ and $ \k^{(5)}$ are set to vanish on slightly different sets.
 	Since  
 	$$ \la\big( (Q^k\cap B_{r+3s}(\xx)) \De (\check Q^{\check k} \cap B_{r+3s}(\xx)))\le \la(N(A_k\cup \check A_{\check k}) \cap B_{r+3s}(\xx)) \le \la(N(\YY)\cap B_{r+3s}(\xx)\big),  $$
 	taking $P=(Q_{(\xx_0,\ff)}\sm N(\CC_{k,\check k})) \cap B_{r+3s}(\xx)$, $P_1=Q^k\cap B_{r+3s}(\xx) $ and $P_2=\check Q^{\check k} \cap B_{r+3s}(\xx)$ in Lemma \ref{lem:Delta}  with $M=4\la(Q)$ shows
 	\begin{align}&\P\big(\XX'''_{\xx,k+1}\ne \check\XX'''_{\xx,\check k+1},  \XX'''_{\xx,\le k}= \check\XX'''_{\xx,\le \check k}, E_{1,k,\check k} \mid \PP_{S_k}^*,\check \PP^*_{\check S_{\check k }}\big) \nonumber 
 		\le \P\big(\PP'' \ne \es \mid \PP_{S_k},\check \PP_{\check S_{\check k }}\big) \\
		&\qquad +4\la\big(Q\big)\big(2+\la(B_{r+3s}(\xx))\big)\la(N(\YY) \cap B_{r+3s}(\xx))  + \P\big(\PP(B_{2r}(\xx))>4\la(Q)\mid \PP_{S_k},\check \PP_{\check S_{\check k }}\big).\label{eq:YY}
 	\end{align}
 	Note that $\P\big(\PP'' \ne \es \mid \PP_{S_k},\check \PP_{\check S_{\check k }} \big)\le \la\big(N(A_k\cup \check A_{\check k})\cap B_{r+3s}(\xx)\big)$ and 
 	\begin{align*}
 		&\E\big[\la\big(N(\YY)\cap B_{r+3s}(\xx)\big)\big] \le 2\E \big[\la\big(N(\PP_{Q\sm B_{r+4s}(\xx)})\cap B_{r+3s}(\xx)\big)\big]  + \la(Q)\P(\YY\cap B_{r+4s}(\xx)\ne  \es)\\
 		&\le 2\E \bigg[ \sum_{\ww\in \PP_{Q\sm B_{r+4s}(\xx)}} \la(N(\ww)\cap B_{r+3s}(\xx) )\bigg] + \la(Q)\P(\YY\cap B_{r+4s}(\xx)\ne  \es)\\
 		&=2 \int_{Q\sm B_{r+4s}(\xx)}  \int_{B_{r+3s}(\xx)}\one\{\ww\sim \zz\} \la(\d\zz)\la(\d \ww)+ \la(Q)\P(\YY\cap B_{r+4s}(\xx)\ne  \es) . 
 	\end{align*}
 	Hence, taking expectations in \eqref{eq:YY}, we obtain
 	\begin{align*}
 		&\E\big[\P(\XX'''_{\xx,k+1}\ne \check\XX'''_{\xx,\check k+1},  \XX'''_{\xx,\le k}= \check\XX'''_{\xx,\le \check k}, E_{1,k,\check k} \mid \PP_{S_k}^*,\check \PP^*_{\check S_{\check k }} )\big] \\
 		&\le C\la(Q)^{3}\bigg( \int_{Q\sm B_{r+4s}(\xx)}  \int_{B_{r+3s}(\xx)}\one\{\ww\sim \zz\} \la(\d\zz)\la(\d \ww)  + \P(\YY \cap B_{r+4s}(\xx)\ne \es) \bigg)+\P(\Ed).
 	\end{align*}
 	By \eqref{cond:sharp} and \eqref{cond:exponential_sim},   this  is of order $\la(Q)^{-l-2}$ if $c_{\ms{St}}$ is chosen large enough.
 	Inserting in \eqref{eq:P'''} yields the bound under $E_{1,k,\check k}$.
	\medskip

\noindent  {\bf Event $E_{2,k, \check k}$.}	Under $E_{2,k,\check k}$, recall that we use the cube-wise ordering introduced in Section \ref{ss:pert}.

Let $\ZZ= T_{Q^k, \es}(\PP_{Q^k}^*) $ and $\check \ZZ= T_{\check Q^{\check k}, \es}(\check \PP_{Q^{\check k}}^*) $ be the Poisson embeddings in step $k+1$ and $\check k +1$, respectively. It is enough to check whether $\ZZ\cap C_j\ne  \check\ZZ \cap C_j$ for all cubes  $C_j$ intersecting $ B_{r+3s}(\xx)$ such that $\PP \cap C_p = \check \PP\cap C_p =\es$ for all $p<j$. For such a $j$, we have 
 	\begin{equation*}
 		\ZZ\cap C_j = T_{Q^{k,j},\phi_{k}}( \PP^*_{Q^{k,j}} )\cap C_j= T_{Q_j\sm N(\CC_{k,\check k}) ,\phi_{k}}^{(6)}\big(\PP^*_{Q^{k,j}(\xx)} \big), 
 	\end{equation*}
 	where $Q_j = Q_{(x_k,\ff)}\cap B_{r+3s}(\xx)  \cap (\cup_{p\ge j} C_p)  $,  $Q^{k,j}=Q^k \cap Q_j =Q_j \sm N(\CC_{k,\check k}\cup A_k)$, and $\k^{(6)}_\xx(\ww, \vp) = \k'''_\xx(\ww, \vp)\one\big\{ \ww \in Q^{k,j}\big\}$. Here $Q^k$ was the unexplored window after step $k$ and $Q_j$ are the unexplored cubes before $C_j$. The analogous statement holds for $\check \ZZ\cap C_j$. Again, $\PP^*_{Q^{k,j} \sm N(\check A_k)} =  \check \PP^*_{\check Q^{\check k,j}\sm N( A_k)} $ is a Poisson process on $Q_j \sm N(\CC_{k,\check k}\cup A_k\cup A_{\check k})=Q^{k,j} \sm N(\check A_k)=\check Q^{\check k,j}\sm N( A_k)$ given $\PP_{S_k \cup (\cup_{p<j} C_p)}$ and $\check \PP_{\check S_{\check k} \cup (\cup_{p<j} C_p)}$. Hence, we may argue as under $E_{1,k,\check k}$ to see that 
 	\begin{align*}
 		&\E\big[\P\big(\XX'''_{\xx,k+1}\cap C_j\ne \check\XX'''_{\xx,\check k+1}\cap C_j, \PP\cap C_p = \check\PP\cap C_p=\es, p<j, \XX'''_{\xx,\le k}= \check\XX'''_{\xx,\le \check k},  \\& \qquad E_{2,k,\check k} \mid \Pds_{S_k\cup(\cup_{p<j} {C_p})}, \check \PP^*_{\check S_{\check k}\cup(\cup_{p<j} {C_p})}\big)\big] \\
 		&\le C\la(Q)^{3}\bigg(\int_{Q\sm B_{r+4s}(\xx)} \int_{B_{r+3s}(\xx)} \one\{ \zz\sim \ww\} \la(\d\zz) \la(\d\ww) + \P(\YY \cap B_{r+4s}(\xx)\ne \es)\bigg)+\P(\Ed).
 	\end{align*}
 	The number of cubes intersecting $B_{r+3s}(\xx)$ is of order $\la(B_{r+3s+1}(\xx))\le\la(Q)$, which leads to the bound
 	\begin{align*}
 		&\P\big(\XX'''_{\xx,k+1}\ne \check\XX'''_{\xx,\check k+1}, \XX'''_{\xx,\le k}\ne \check\XX'''_{\xx,\le \check k}, E_{2,k,\check k}\big) \\
 		&\le C\la(Q)^{4}\bigg(\int_{Q\sm B_{r+4s}(\xx)}  \int_{B_{r+3s}(\xx)}\one\{\ww\sim \zz\} \la(\d\zz)\la(\d \ww)  + \P(\YY \cap B_{r+4s}(\xx)\ne \es)+\P(\Ed)\bigg).
 	\end{align*}
 	This is again of order $\la(Q)^{-l-2}$ when $c_{\ms{St}}$ is large enough.
	\medskip

 
\noindent{\bf Model (P)} 	For Model (P), we let $\G=\GGc(\Phi_{Q,\es}^*)$ and $\check\G=\GGc(\check \Phi_{Q,\es}^*)$ and define   
 \begin{align*}	\YY_0 &=\big\{ \ww \in \PP \cup \check\PP  \mid \ww \lrsa_{\G } \PP_{Q\sm B_{2r}(\xx)} \text{ or } \ww \lrsa_{\check \G } \check \PP_{Q\sm B_{2r}(\xx)}\big\}\\
 		\YY &=\big\{ \ww \in \PP \cup \check\PP  \mid \ww \lrsa_{\G } \YY_0\cap \PP \text{ or } \ww \lrsa_{\check \G } \YY_0 \cap \check \PP\big\}.
 		\end{align*}
Here we construct $\YY$ in two steps to ensure that it is a union of clusters in both $\G$ and $\check \G$. This is has to do with the exploration based construction of the RCMs, which makes it possible to have an edge between two points in $\PP_{B_{2r}(\xx)}$ in $\G$ but not in $\check \G$. This was not possible when the graphs were based on a symmetric relation.
 	Note that $\P(\YY\cap B_{r+4s}(\xx)\ne \es) \le \P(\YY_0\cap B_{r+4.5s}(\xx)\ne \es) + 2\P(\PP_{B_{r+4s}(\xx)} \lrsa_\G \PP_{Q\sm B_{r+4.5s}(\xx)}),$
 	which we may bound using the assumption \eqref{cond:sharp_vv}.
 		
 	The two graphs $\G$ and $\check \G$ have the same clusters on $\PP\sm \YY = \check \PP\sm \YY$, and the points are visited in the same order. To see this, consider the construction of $\G$ and $\check \G$ as thinnings of directed graphs $\overrightarrow{\G}$ and $\overrightarrow{\check{\G}}$ from Remark \ref{rem:directed}. Recall that a cluster in $\G$  starts at a point $\xx_0\in\PP$ that has no directed path in $\overrightarrow{\G}$ to an $\iota$-smaller point.  If $\xx_0\notin \YY$, there can be no directed path from $\xx_0$ to a smaller point in $\overrightarrow{\check{\G}}$ either since this would have to pass through $\check \PP_{Q\sm B_{2r}(\xx)}$ and hence $\xx_0$ would be in $\YY_0$. Thus, $\xx_0$ must also be a starting point of a cluster in $\check \G$. If a point $\zz$ is in the cluster of $\xx_0 $ in $\G$, then  there is a directed path from $\zz$ to $\xx_0$ in $\overrightarrow{\G}$, but not to any smaller point. This path cannot go via any points of $\YY_0$, otherwise $\xx_0$ would be in $\YY$. Hence, the path would also be a path in $\overrightarrow{\check{\G}}$. There cannot be any path from $\zz$ to a smaller point than $\xx_0$ in $\overrightarrow{\check{\G}}$ since such a path would have to go via $\check \PP_{Q\sm B_{2r}(\xx)}\su \YY_0$, which would again imply $\zz\in \YY_0$ and hence  $\xx_0\in \YY$. Thus, the clusters starting at an $\xx_0\notin \YY$ agree, and they are explored in the same order. In particular, at each steps $k,\check k$ in the two explorations, the clusters of $\GGc_k$ and $\check \GG_{\check k}^{\ms{cl}}$ consist of a union of clusters $\CC_{k,\check k}$ that agree (including the one being explored, if it agrees so far) and some extra points $A_k\su \PP$ and $\check A_{\check k }\su \check \PP$,  just as under Model (L).

 	We define events $E_{1,k,\check k}$ and $ E_{2,k,\check k}$ as for Model (L). Then,  \eqref{eq:P'''}  holds again, where the expectation for given $k,\check k$ can be computed as
 	$$\lim_{\de \to 0}\E[\P(\XX^{\de,\prime\prime\prime}_{\xx,{k+1}}\ne \check\XX^{\de,\prime\prime\prime}_{\xx,\check k +1},  \XX'''_{\xx,\le k}= \check\XX'''_{\xx,\le \check k}, E_{1,k,\check k}\cup E_{2,k,\check k} \mid \Phi_{Q,\es}^*\cap S_k , \check \Phi^*_{Q,\es} \cap \check S_{\check k}) ]. $$
 	Here $\XX^{\de,\prime\prime\prime}_{\xx,{k+1}}$ and $\check\XX^{\de,\prime\prime\prime}_{\xx,\check k +1}$ are the processes generated on $Q^\de$ in step $k+1$ and $\check k +1$, respectively, and $\XX'''_{\xx,\le k}$ and  $\check\XX'''_{\xx,\le \check k}$ are the parts of $\XX'''_{\xx}$ and  $\check\XX'''_{\xx}$ generated before step $k+1$ and $\check k + 1$. 
 	
 	Inductive application of Lemma \ref{lem:RCM_it} shows that conditionally on $\Phi_{Q,\es}^*\cap S_k$, we have that $\bar B^*_{r+3s}(\xx)\cap (\Phi_{Q,\es}^* \sm \bar S_k^*)$ is a Poisson process of intensity 
	$$ \rho_k(\zz)=\one\big\{\zz\in Q_{(\xx_k,\ff)}\cap B_{r+3s}(\xx)\big\} \prod_{\yy\in S_{k-1}} (1-\pi(\zz,\yy)).$$ 
 	 Moreover, conditionally  on $\Phi_{Q,\es}^*\cap \bar S_k^*$ and $\check \Phi_{Q,\es}^*\cap \Bar{\Check{S}}_{\check k}^*$, we have that $\bar B^*_{r+3s}(\xx)\cap \big(\Phi_{Q,\es}^* \sm (\bar S_k^*\cup \Bar{\Check{S}}^*_{\check k})\big)$ is a Poisson process of reduced intensity $\rho_{k,\check k}\la$, where 
	 $$\rho_{k,\check k}(\zz)=\one\big\{\zz\in Q_{(\xx_k,\ff)}\cap Q_{( \check \xx_{\check k},\ff)}\cap  B_{r+3s}(\xx)\big\}\prod_{\yy\in S_{k-1}\cup \check S_{\check k-1} }(1-\pi(\zz,\yy)).$$    
	 It follows that
 	\begin{align}\nonumber
 		&\P(\bar B^*_{r+3s}(\xx)\cap (\Phi_{Q,\es}^* \sm \bar S_k^* )\cap \Bar{\Check{S}}^*_{\check k} \ne \es, E_{1,k,\check k}\cup E_{2,k,\check k}) \\ \nonumber
 		&\le \E\big[ \big(\big| \bar B^*_{r+3s}(\xx)\cap (\Phi_{Q,\es}^* \sm \bar S_k^*) \big|- \big| \bar B^*_{r+3s}(\xx)\cap (\Phi_{Q,\es}^* \sm (\bar S_k^* \cup \Bar{\Check{S}}^*_{\check k}))\big|\big)\one\big\{ E_{1,k,\check k}\cup E_{2,k,\check k} \big\}\big]\\ \nonumber
 		&=\E \bigg[\bigg(\int_{B_{r+3s}(\xx)} \rho_{k}(\zz)\d \zz - \int_{B_{r+3s}(\xx)} \rho_{k,\check k}(\zz)\d \zz\bigg) \one\big\{ E_{1,k,\check k}\cup  E_{2,k,\check k} \big\} \bigg] \\ \nonumber
 		&\le \E \bigg[\int_{B_{r+3s}(\xx)} \bigg(1-\prod_{\yy \in \check A_{\check k}}(1-\pi(\zz,\yy))\bigg)\d \zz  \one\big\{ E_{1,k,\check k}\cup  E_{2,k,\check k}  \big\} \bigg]  \\ \nonumber
 		&\le \E\bigg[\one\big\{\YY \cap  B_{r+4s}(\xx) =\es, \check A_{\check k} \su \YY \big\}\sum_{\yy\in \check A_{\check k}}\int_{ B_{r+3s}(\xx)} \pi(\zz,\yy) \d \zz\bigg] + \la( B_{r+4s}(\xx))\P(\YY\cap  B_{r+4s}(\xx) \ne \es)\\
 		& = 2\int_{Q\sm  B_{r+4s}(\xx)} \int_{ B_{r+3s}(\xx)} \pi(\zz,\yy) \d \zz\d\yy + \la( B_{r+4s}(\xx))\P(\YY\cap  B_{r+4s}(\xx) \ne \es). \label{eq:pi}
 	\end{align}
  	Then,  $\P(\bar  B^*_{r+3s}(\xx)\cap (\Phi_{Q,\es}^* \sm \bar S_k^*)  \ne  \bar B^*_{r+3s}(\xx)\cap (\check{\Phi}_{Q,\es}^* \sm  \Bar{\Check{S}}^*_{\check k}))$ is bounded by twice the bound in \eqref{eq:pi}. 
 	
 	We may extend the Poisson process $\PP_{B_{r+3s}(\xx)} \sm (S_k\cup \check S_{\check k})$ to a Poisson process $\PP^\prime$ of intensity $\rho_{k,\check k}(\zz)=\prod_{\yy \in \CC_{k,\check k}}(1-\pi(\zz,\yy))$ on $B_{r+3s}(\xx)$ and attach arbitrary marks to the additional points to form the process $\Phi'$. We can bound the probability $\PP' \ne \PP_{B_{r+3s}(\xx)} \sm (S_k\cup \check S_{\check k})$ by computing as in \eqref{eq:pi}. Thus we may assume that we are thinning $\Phi'$. Lifting $\Phi'$ to $Q^\de$ and arguing as in  Model (L), we find that
 	\begin{align*}
 		&\P(\XX_{\xx,k+1}^{\de,\prime\prime\prime} \ne \check \XX_{\xx,\check k+1}^{\de,\prime\prime\prime}, \XX_{\xx,\le k}^{\de,\prime\prime\prime} = \check \XX_{\xx,\le \check k}^{\de,\prime\prime\prime},E_{1,k,\check k}\cup  E_{2,k,\check k}  )      \le  C\la(Q)^{4}\bigg(  \P(\YY \cap B_{r+4s}(\xx) \ne \es) + \P(\Ed) +\P(E^\de)\\& +\int_{Q\sm  B_{r+4s}(\xx)} \int_{ B_{r+3s}(\xx)} \pi(\xx,\yy) \d \xx\d\yy \\
& 		+\int_{B_{r+3s}(\xx)\ti \Uu^\de} \int_{Q^{\de} \sm (B_{r+4s}(\xx)\ti \Uu^\de)}\one\{ (\zz,\rr)\approx^\de (\yy,\ss)\}  (\la\ti \Leb_{[0,1]}^\de)^2(\d(\yy,\ss,\zz,\rr))\bigg)
 	\end{align*} 
 	where $\approx^\de$ is the symmetric relation corresponding to $\pi$, see  Section \ref{ss:pert_RCM}.
 	Integrating over $\Uu^\de\ti \Uu^\de$ and letting $\de \to \infty$, the second integral converges to
 	$$\int_{B_{r+3s}(\xx)} \int_{Q \sm B_{r+4s}(\xx)}\pi(\zz,\ww)\la(\d \yy)\la(\d \zz).$$
 	The claim thus follows from the assumptions \eqref{cond:sharp_vv} and \eqref{cond:sharp_v}.
 \enp

	\section*{Acknowledgements}
	C. Hirsch was supported by a research grant (VIL69126) from VILLUM FONDEN. M. Otto was supported by the NWO Gravitation project NETWORKS under grant agreement no. 024.002.003 and a grant (W253098-1-035) from the Drs.~J.R.D. Kuikenga Fonds voor Mathematici.

	\bibliographystyle{plainnat}
	\bibliography{./lit}
\end{document}